\documentclass[amstex,11pt,a4paper]{article}

\usepackage[cp1251]{inputenc}
\usepackage[russian]{babel}
\usepackage{amsfonts,amssymb}
\usepackage{graphicx}%[pdftex]
\usepackage{latexsym,exscale}
\usepackage{amsmath}
\usepackage{amsthm}
\usepackage{blindtext}

\usepackage{amsmath}
\usepackage{relsize}

\graphicspath{{pictures/}}
\DeclareGraphicsExtensions{.pdf,.png,.jpg}

\usepackage[version=0.96]{pgf}
\usepackage{tikz}

\usetikzlibrary{intersections}

\numberwithin{equation}{section}

\newtheorem{Lemma}{Лемма}

\newtheorem{Corollary}{Следствие}

\newtheorem*{Observation}{Замечание}

\newtheorem*{Main Theorem}{Основная теорема}
\newtheorem*{HomStatement}{Утверждение о гомоморфизмах}

\newtheorem*{WLemma}{Лемма о словах}
\newtheorem*{KurTh}{Теорема Куроша
{\rm\cite{[Kur67]}}}

\newtheorem*{FiWiTh}{Теорема Файна-Вильфа
{\rm\cite{[FiWi65]}}}

\newtheorem*{CSVCGrTh}{Теорема
{\rm\cite{[Mazh18]},\cite{[KlMaMi]}}}

\newtheorem*{ThKM}{Теорема
{\rm\cite{[KM17]}}}

\newtheorem*{BSTh}{Теорема Басса-Серра
{\rm\cite{[Ser77]},\cite{[Har00]}}}

\newenvironment{ProofMT} % имя окружения
{\par\noindent{\bf{Доказательство Основной теоремы.}}} % команды для \begin
{\hfill$\scriptstyle\square \\$}

\newenvironment{ProofHomStatement}
{\par\noindent{\bf{Доказательство Утверждения о гомоморфизмах.}}}
{\hfill$\scriptstyle\square \\$}

\newenvironment{ProofWL} % имя окружения
{\par\noindent{\bf{Доказательство Леммы о словах.}}} % команды для \begin
{\hfill$\scriptstyle\square \\$}

\newenvironment{ProofMLIII}
{\par\noindent{\bf{Доказательство леммы \ref{LemMLPartIII}.}}}
{\hfill$\scriptstyle\square \\$}

\newenvironment{ProofMLII}
{\par\noindent{\bf{Доказательство леммы \ref{LemMLPartII}.}}}
{\hfill$\scriptstyle\square \\$}

{\par\noindent{\bf{Доказательство.}}} % команды для \begin
{\hfill$\scriptstyle\blacksquare \\$}

\newcounter{Remark}

\font\MMMM=msbm10

\def\ltimes{\mathop{\raise0.5pt\hbox{\MMMM\char"6E}}}
\def\rtimes{\mathop{\raise0.5pt\hbox{\MMMM\char"6F}}}

\let\tmp\H
\def\H#1{\ifmmode{{\fam\trafam H}}#1\else\tmp#1\fi}
.pk scaled 1440 %2000
\font\ZV=cmr10.pk scaled 2000
\def\zvezda{\mathop{\,\vrule width0pt depth2pt height8pt
            \smash{\lower7pt\hbox{\ZV *}}\,}\limits}
\font\bigcmsy=cmsy10.pk scaled 2000
\def\bigtimes{\mathop{\,\vrule width0pt depth2pt height8pt
            \smash{\lower2pt\hbox{\bigcmsy\char'002}}\,}\limits}
.pk scaled 1200

\def\krest{\mathop{\,\vrule width0pt depth2pt height8pt
            \smash{\lower0pt\hbox{\ZV $\times$}}\,}\limits}

\def\acts{\!\curvearrowright\!}

\newbox\tmpbox
\newdimen\tmpdim
\def\narrow[#1]#2\par %?????? #2 ?????? ?? ?????? #1, ???????? ? \tmpbox
{%
\setbox\tmpbox\hbox{#2}%
\ifdim\wd\tmpbox>#1
\setbox\tmpbox\vbox{\hsize=#1 #2}%
\tmpdim=\ht\tmpbox%
\loop%
\setbox\tmpbox\vbox{\hsize=\wd\tmpbox \advance\hsize by -1pt
#2}%
\ifdim\ht\tmpbox=\tmpdim%
\relax%
\repeat%
\setbox\tmpbox\vbox{\hsize=\wd\tmpbox \advance\hsize by 1pt #2}%
\fi%
\box\tmpbox
}
\def\disp#1{
                $$
                \setbox\tmpbox\vbox{\narrow[\hsize]\noindent#1\par}
                \box\tmpbox
                $$
          }

\usepackage{authblk}
\title{Свободные произведения групп сильно вербально замкнуты}
\author{Андрей Михайлович Мажуга}
\affil{Механико-математический факультет МГУ им.~М.В.Ломоносова, аспирант}
\affil{Всероссийская академия внешней торговли, младший научный сотрудник Института макроэкономических исследований}
\date{\small \begin{flushleft} УДК 512.543.72 \end{flushleft}}
\makeatletter
\def\@xfootnote[#1]{%
  \protected@xdef\@thefnmark{#1}%
  \@footnotemark\@footnotetext}
\makeatother

\begin{document}
\maketitle

\begin{abstract}
В ряде недавних работ было установлено, что многие почти свободные группы, фундаментальные группы почти всех связных поверхностей и все группы, являющиеся нетривиальными свободными произведениями групп с тождествами, алгебраически замкнуты в \emph{любой} группе, в которой они вербально замкнуты. В данной работе мы установим, что любая группа, являющаяся нетривиальным свободным произведением групп, алгебраически замкнута в любой группе, в которой она вербально замкнута.\\
\\
{\bf Ключевые слова:} вербально замкнутые подгруппы, алгебраически замкнутые подгруппы, ретракты групп.
\end{abstract}
{\def\thefootnote{}\addtocounter{footnote}{-1}%
%\footnote{%
%Работа выполнена при поддержке Российского фонда фундаментальных
%исследований, грант 15-01-05823} }

\section{Введение}

Подгруппа $H$ группы $G$ называется \emph{вербально замкнутой} (в $G$)~\cite{[MR14]} (см. также \cite{[Rom12]}, \cite{[RKh13]}, \cite{[KM17]}, \cite{[Mazh17]}, \cite{[Mazh18]}, \cite{[KlMaMi]}), если любое уравнение вида $w(x_1,\dots,x_n)=h$, где $w(x_1,\dots,x_n)\in F_n(x_1,\dots,x_n)$ и $h\in H$, имеющее решение в $G$ имеет решение и в $H$.

Подгруппа $H$ группы $G$ называется \emph{алгебраически замкнутой} (в $G$), если любая система уравнений вида $\left \{ w_1\left (x_1,\dots,x_n, H\right)=1,\dots,w_m\left ( x_1,\dots,x_n, H\right)=1 \right \}$, где $w_i\left (x_1,\dots,x_n,H\right) \in F_n(x_1,\dots,x_n)\ast H$, имеющая решение в $G$ имеет решение и в $H$.

Подгруппа $H$ группы $G$ называется \emph{ретрактом} (группы $G$), если $G$ представима как полупрямое произведение нормальной подгруппы $N$ и подгруппы $H$ (т.е. $G=N\rtimes H$).

Из приведенных определений легко следует, что любая подгруппа, являющаяся ретрактом, является и алгебраически замкнутой подгруппой, а любая алгебраически замкнутая подгруппа является и вербально замкнутой. В общем случае обратные импликации неверны в обоих случаях (см., например,~\cite{[KM17]}). Следовательно, возникает вопрос: при каких ограничениях на подгруппу~$H$ и группу~$G$ выполняются обратные импликации? В работе~\cite{[MR14]} были установлены следующие достаточно общие результаты:
{\it
\begin{enumerate}
    \item[\sffamily $\mathrm{1)}$] если группа $G$ конечно представима и $H$ есть конечно порожденная алгебраически замкнутая подгруппа в $G$, то $H$ является ретрактом группы $G$;

    \item[\sffamily $\mathrm{2)}$] если группа $G$ конечно порождена над $H$\footnote[1]{Группа $G$ \emph{конечно порождена над $H$}, если $G=\left \langle H,X \right \rangle$ для некоторого конечного подмножества $X\subseteq G$.} и $H$ есть нетерова по уравнениям\footnote[2]{Группа $H$ называется \emph{нетеровой по уравнениям}, если любая система уравнений с коэфффициентами из $H$ и конечным числом неизвестных эквивалентна некоторой своей конечной подсистеме.} алгебраически замкнутая подгруппа в $G$, то $H$ является ретрактом группы $G$.
\end{enumerate}
}

Для случая вербально замкнутых подгрупп подобных общих результатов не известно, а известные результаты относятся к определенным типам подгрупп. Так в конечно порожденных свободных группах (см.~\cite{[MR14]}) и конечно порожденных свободных нильпотентных группах (см.~\cite{[RKh13]}) вербально замкнутая подгруппа всегда является ретрактом. Следующий результат был установлен в работе~\cite{[KM17]}:
\begin{ThKM}
Пусть~$G$ есть произвольная группа и~$H$ есть её
вербально замкнутая
%конечно порождённая  !!!!!!!!
почти свободная
бесконечная
недиэдральная\footnote[3]{Для бесконечной группы
\emph{недиэдральная} означает неизоморфная свободному произведению
двух групп порядка два.}
подгруппа,
любая бесконечная абелева подгруппа которой является циклической.
Тогда
\begin{enumerate}
\item[\rm1)]
%если $H$ нециклическая подгруппа, то
$H$ алгебраически замкнута в $G$;
\item[\rm2)]
если группа $G$
конечно порождена
над $H$, то $H$ является
ретрактом группы $G$.
\end{enumerate}
\end{ThKM}
Мы будем называть группу $H$ \emph{сильно вербально замкнутой}~\cite{[Mazh18]}, если $H$ является алгебраически замкнутой подгруппой в любой группе, в которой она является вербально замкнутой подгруппой (таким образом, вербальная замкнутость --- это свойство подгруппы, а сильная вербальная замкнутость --- это свойство абстрактной группы). Заметим, что первый пункт приведенной выше теоремы описывает определенный класс сильно вербально замкнутых групп (в частности, все нетривиальные свободные группы принадлежат этому классу). Из результатов работ~\cite{[Mazh18]} и~\cite{[KlMaMi]} следует:

\begin{CSVCGrTh}Сильно вербально замкнутыми являются:
\begin{enumerate}
 \item[\sffamily $\mathrm{\bullet}$]
все абелевы группы;
 \item[\sffamily $\mathrm{\bullet}$]
фундаментальные группы всех связных поверхностей, за исключением, возможно, бутылки Клейна;
 \item[\sffamily $\mathrm{\bullet}$]
все почти свободные группы, не имеющие неединичных нормальных конечных подгрупп;
 \item[\sffamily $\mathrm{\bullet}$]

все свободные произведения $\zvezda\limits_{i\in I}\!H_i$, где $I$ есть любое множество, содержащие хотя бы два элемента, и $H_i$ суть нетривиальные группы, удовлетворяющие нетривиальным тождествам\footnote[4]{Если~$I(x_1,\dots,x_r)\ne1$ есть элемент свободной группы~$F_r(x_1,\dots,x_r)$ с базисом~$x_1,\dots,x_r$, то \emph{группа~$G$ удовлетворяет нетривиальному тождеству $I$}, если для всех~$g_1,\dots,g_r\in G$ выполнено~$I(g_1,\dots,g_r)=1$ в~$G$.} (возможно, разным).
\end{enumerate}
\end{CSVCGrTh}

Мы будем называть группу $H=\zvezda\limits_{\mathfrak{a}\in \mathfrak{A}}\!H_\mathfrak{a}$ \emph{нетривиальным свободным произведением} (групп~$H_\mathfrak{a}$, $\mathfrak{a}\in \mathfrak{A}$), если множество индексов $\mathfrak{A}$ содержит хотя бы два элемента и все группы $H_\mathfrak{a}$ нетривиальны. Мы будем называть элемент $w\in H$ \emph{гиперболическим}, если он не принадлежит никакой подгруппе (группы $H$) вида $H_\mathfrak{a}^g$, где $\mathfrak{a}\in \mathfrak{A}$ и $g\in H$.

Введем некоторые обозначения. Запись $\left \langle X \right \rangle$ означает подгруппу группы $G$, порожденную множеством $X\subseteq G$. Символы $\left \langle x \right \rangle_\infty$ и $\left \langle x \right \rangle_n$ означают бесконечную циклическую группу и циклическую группу порядка $n$, порожденную элементом $x$, соответственно. $F_n(z_1,\dots,z_n)$ есть свободная группа ранга $n$ (с базисом $z_1,\dots,z_n$). Коммутант группы $G$ мы обозначаем как $G'$.

Основным результатом данной работы является доказательство следующей теоремы.
\begin{Main Theorem}\label{TMTm} Группа, являющаяся нетривиальным свободным произведением групп, сильно вербально замкнута.
\end{Main Theorem}

Если рассматриваемая в данной теореме группа изоморфна группе $\langle a\rangle_2\ast\langle b\rangle_2$, то утверждение непосредственно следует из последнего пункта предыдущей теоремы. Для случая, когда рассматриваемая группа не изоморфна $\langle a\rangle_2\ast\langle b\rangle_2$, основная теорема легко выводится (см. раздел~5) из пункта~$a)$ следующей леммы.

\begin{WLemma}\label{LemMLPartIII} Пусть группа $H=\zvezda\limits_{\mathfrak{a}\in \mathfrak{A}}\!H_\mathfrak{a}$ есть нетривиальное свободное произведение групп~$H_\mathfrak{a}$ и пусть $\underline{x}=(x_0,\dots,x_{n-1})$ есть такой кортеж элементов из $H$, что подгруппа $\left \langle x_0,\dots,x_{n-1} \right \rangle$ содержит хотя бы два некоммутирующих гиперболических элемента группы $H$. Тогда:
\begin{enumerate}
     \item[\sffamily $\mathrm{a)}$] существует слово $M_{\underline{x}}(z_0,\dots,z_{n-1})\in F_n(z_0,\dots,z_{n-1})'$, вообще говоря, зависящее от кортежа $\underline{x}$, такое, что если для некоторого кортежа $(y_0,\dots,y_{n-1})$ элементов группы $H$ выполнено равенство:
\begin{equation}\label{LemMLPartIIIMEQff}
M_{\underline{x}}(x_0,\dots,x_{n-1})=M_{\underline{x}}(y_0,\dots,y_{n-1}),
\end{equation}
то найдется такой элемент $v\in \left \langle M_{\underline{x}}(x_0,\dots,x_{n-1}) \right \rangle_\infty$, что $x_i = v^{-1}y_iv$ при $i=0,\ \dots,\ n-1$;

     \item[\sffamily $\mathrm{b)}$] существуют слова $M_{\underline{x}}'(z_0,\dots,z_{n-1}), M_{\underline{x}}''(z_0,\dots,z_{n-1})\in F_n(z_0,\dots,z_{n-1})'$, вообще говоря, зависящие от кортежа $\underline{x}$, такие, что если для некоторого кортежа $(y_0,\dots,y_{n-1})$ элементов группы $H$ одновременно выполнены равенства:
\begin{equation}\label{LemMLPartIIIMEQpartb}
M_{\underline{x}}'(x_0,\dots,x_{n-1})=M_{\underline{x}}'(y_0,\dots,y_{n-1})\ \ \text{и}\ \ M_{\underline{x}}''(x_0,\dots,x_{n-1})=M_{\underline{x}}''(y_0,\dots,y_{n-1}),
\end{equation}
то $x_i = y_i$ при $i=0,\ \dots,\ n-1$.
 \end{enumerate}
\end{WLemma}
%\begin{WLemma}\label{TFML} Пусть группа $H=\zvezda\limits_{\mathfrak{a}\in \mathfrak{A}}\!H_\mathfrak{a}$ есть нетривиальное свободное произведение групп $H_\mathfrak{a}$ и пусть %$\underline{x}=(x_0,\dots,x_{n-1})$ есть такой кортеж элементов из $H$, что подгруппа $\left \langle x_0,\dots,x_{n-1} \right \rangle$ содержит хотя бы два некоммутирующих %гиперболических элемента группы $H$. Тогда существует слово $M_{\underline{x}}(z_0,\dots,z_{n-1})\in F_n(z_0,\dots,z_{n-1})'$, вообще говоря, зависящее от кортежа $\underline{x}$, %такое, что если для некоторого кортежа $(y_0,\dots,y_{n-1})$ элементов $y_i\in H$ выполнено равенство:
%\begin{equation}\label{AddPartIVEQ1}
%M_{\underline{x}}(x_0,\dots,x_{n-1})=M_{\underline{x}}(y_0,\dots,y_{n-1}),
%\end{equation}
%то найдется такой элемент $v\in \left \langle M_{\underline{x}}(x_0,\dots,x_{n-1}) \right \rangle_\infty$, что $x_i = v^{-1}y_iv$ при $i=0,\ \dots,\ n-1.$
%\end{WLemma}

Отметим, что из пункта~$b)$ сформулированной леммы легко выводится (см. раздел 5) следующее утверждение.
\begin{HomStatement} Пусть $G$ есть конечно порожденная группа, группа $H=\zvezda\limits_{\mathfrak{a}\in \mathfrak{A}}\!H_\mathfrak{a}$ есть нетривиальное свободное произведение групп $H_\mathfrak{a}$ и {\rm $\text{\textrm{Hom}}(G,H)$} есть множество гомоморфизмов из $G$ в $H$. Если гомоморфизм {\rm$\psi \in \text{\textrm{Hom}}(G,H)$} таков, что его образ
%{\rm$\textrm{Im}(\psi)$}
содержит хотя бы два некоммутирующих гиперболических элемента группы $H$, то найдутся такие элементы $w_1,w_2\in G$, вообще говоря, зависящие от $\psi$, что для любого гомоморфизма {\rm$\varphi \in \text{\textrm{Hom}}(G,H)$} из равенств $\varphi(w_1)=\psi(w_1)$ и $\varphi(w_2)=\psi(w_2)$ следует, что $\varphi=\psi$.
\end{HomStatement}

Практически вся оставшаяся часть работы посвящена доказательству Леммы о словах, которое состоит из трех частей. Первая и вторая части посвящены соответственно доказательству лемм~\ref{LemMLPartII} и~\ref{LemMLPartIII} (см. разделы~3 и~4), являющихся частными случаями Леммы о словах. В разделе~5 мы завершаем доказательство Леммы о словах и выводим из нее Основную теорему и Утверждение о гомоморфизмах.

%Суть нашего доказательства, вероятно, наилучшим образом иллюстрирует следующее элементарное утверждение. Над не алгебраически замкнутым полем $K$ для любого $n\in\mathbb{N}$ %существует многочлен $T_n(z_0,\dots,z_{n-1})\in K[z_0,\dots,z_{n-1}]$, такой, что если для некоторых $k_0,\ \dots,\ k_{n-1}\in K$ выполнено равенство %$T_n(k_0,\dots,k_{n-1})=T_n(0,\dots,0)$, то $k_0=\dots=k_{n-1}=0$. В силу существования таких многочленов, над $K$ любая система уравнений вида $\{ f_i(x_1,x_2,\dots)=0\ |\ i=0,\ %\dots,\ n-1\}$ эквивалентна одному уравнению вида $T_n(f_0(x_1,x_2,\dots),\dots,f_{n-1}(x_1,x_2,\dots))=T_n(0,\dots,0)$. Для того, чтобы доказать основную теорему мы показываем, что %в определении алгебраически замкнутой подгруппы можно рассматривать более узкий класс систем уравнений (см. лемму~\ref{ProofMT2TL}) и для систем из этого класса мы строим слово %$L_n(z_0,\dots,z_{n-1})\in F_n(z_0,\dots,z_{n-1})$ (см. определения (\ref{L2Def}) и (\ref{MW})), играющее в нашем доказательстве роль схожую с ролью многочленом %$T_n(z_0,\dots,z_{n-1})$ из вышеприведенного утверждения (см. пункт $L3.2)$ леммы \ref{LemTML1}).
\

{\bf Обозначения}, которыми мы пользуемся, в целом стандартны. Если $k\in \mathbb{Z}$, $x$ и $y$ суть элементы группы, то записи $x^y$, $x^{ky}$, $x^{-y}$ и $x^{-ky}$ обозначают элементы $y^{-1}xy$, $y^{-1}x^ky$, $y^{-1}x^{-1}y$ и $y^{-1}x^{-k}y$ соответственно. Коммутатор элементов $x$ и $y$ группы $G$, $[x,y]$, мы определяем как $x^{-1}y^{-1}xy.$
%Запись $\left \langle X \right \rangle$ означает подгруппу группы $G$, порожденную множеством $X\subseteq G$.
Централизатор подмножества $X$ группы $G$ мы обозначаем как $C_G(X)$.
%Запись $\left \langle x \right \rangle_\infty$ означает бесконечную циклическую группу, порожденную элементом $x$.
Иногда вместо записи вида $x_1,\dots,x_n$ мы используем сокращенную запись вида $\underline{x}$ и записываем, к примеру, $w(x_1,\dots,x_n)$ как $w(\underline{x})$ или $F_n(x_1,\dots,x_n)$ как $F(\underline{x})$. Мы определяем функции $\lfloor\cdot\rfloor:x \mapsto \lfloor x\rfloor$ и $\lceil\cdot\rceil:x \mapsto \lceil x\rceil$ как $\lfloor x\rfloor=\max\{n\in\mathbb{Z}\ |\ n\le x\}$ и $\lceil x\rceil=\min\{n\in\mathbb{Z}\ |\ n\ge x\}$ соответственно.

\section{Вспомогательные определения и утверждения}

Если группа $H=\zvezda\limits_{\mathfrak{a}\in \mathfrak{A}}\!H_\mathfrak{a}$ есть нетривиальное свободное произведение групп $H_\mathfrak{a}$ и если $w_1,w_2$ суть слова из $H$, то запись $w_1\equiv w_2$ будет означать, что эти слова равны посложно (т.е. у них равны первые, вторые и т.д. вплоть до последнего слоги).

Мы будем говорить, что слово $w\equiv a_1\cdots a_m\in H$:
\begin{itemize}
\item[$\bullet$]
\emph{приведено}, если или $m=0$ (в этом случае мы пишем $w\equiv1$) или если $m\ge1$, $a_j\in H_{\mathfrak{a}_j}\!\setminus\!\{1\}$, $j=1,\ \dots,\ m$  и слоги $a_j,a_{j+1}$ при $j=1,\ \dots,\ m-1$ принадлежат различным $H_\mathfrak{a}$;

\nobreak
\item[$\bullet$]
\emph{циклически приведено}, если или $m=1$ и $a_1\ne 1$ или если $m\ge 2$, $w$ приведено и слоги $a_1$, $a_m$ принадлежат различным $H_\mathfrak{a}$;

\nobreak
\item[$\bullet$]
\emph{простое}, если $m\ge 2$, $w$ циклически приведено и не является собственной степенью\footnote[$\dagger$]{Мы говорим, что элемент $g\in G$ есть \emph{собственная степень} в группе $G$, если $g=h^n$ для некоторого $h\in G$ и целого числа $n\ge 2$.} в $H$.
\end{itemize}

Если $w_1,w_2\in H$ суть приведенные слова, то запись $w_1\cdot w_2$ будет означать, что или $w_1\equiv1$ или $w_2\equiv1$ или $w_1,w_2\not\equiv1$ и последний слог слова $w_1$ и первый слог слова $w_2$ принадлежат различным~$H_\mathfrak{a}$; запись $w_1\circ w_2$ будет означать, что или $w_1\equiv1$ или $w_2\equiv1$ или $w_1,w_2\not\equiv1$ и между последним слогом слова $w_1$ и первым слогом слова $w_2$ возможно лишь слияние, но не сокращение. Например, в группе $\left \langle a \right \rangle_5\ast\left \langle b \right \rangle_2$ мы можем написать $ba^2b=b\cdot a^2b=ba^4\circ a^3b.$

Множество всех гиперболических элементов группы $H$ мы будем обозначать как $M_H$. Если $w\in M_H$, то легко видеть, что слово $w$ может быть единственным образом представлено в виде
\begin{equation}\label{A2107_1}
w=f^{-1}\cdot A^k\circ f,
\end{equation}
где $A$ --- простое слово, $f$ --- приведенное слово и $k\in\mathbb{N}.$

Мы будем говорить, что слово $w_1$ \emph{сопряжено} слову $w_2$ элементом $f$ (и обозначать это как $w_1=w_2^{f}$), если $w_1=f^{-1}w_2f_1$.

Несложно убедиться в том, что следующие понятия определены корректно:

\begin{itemize}
\item[$\bullet$]
\emph{длиной} слова $w\in H\!\setminus\!\{1\}$, $|w|$, мы будем называть число слогов в его приведенной форме;

\item[$\bullet$]
\emph{центральной длиной} слова $w\in H\!\setminus\!\{1\}$, $|w|_c$, мы будем называть длину циклически приведенного слова, которому сопряжено $w$;

\item[$\bullet$]
если слово $w\in M_H$ представлено в виде~(\ref{A2107_1}), то его \emph{радикальной длиной}, $|w|_r$, мы будем называть длину простого слова $A$.
\end{itemize}
Например, если слово $w$ имеет вид~(\ref{A2107_1}), то $|w|=k|A|+2|f|$ или $|w|=k|A|+2|f|-1$ (в зависимости от того, произошло ли слияние или нет), $|w|_c=k|A|$ и $|w|_r=|A|$. Легко понять, что для любого слова $w\in M_H$ и любого $k\in \mathbb{N}$ мы имеем: $|w|_r\le|w|_c\le|w|$ (причем первое неравенство превращается в равенство тогда и только тогда, когда $w$ не есть истинная степень в $H$, а второе тогда и только тогда, когда слово $w$ циклически приведено), $|w^k|\le k|w|$, $|w^k|_c=k|w|_c$ и $|w^k|_r=|w|_r$. Заметим также, что $w\in M_H$ тогда и только тогда, когда $|w|_c\ge2$, и что радикальную длину мы определяем лишь для гиперболических элементов группы~$H.$

Мы будем говорить, что слово $w_1$ есть \emph{циклическая перестановка} циклически приведенного слова $w_2$ при помощи элемента $f$ если $w_1=f^{-1}w_2f_1$ и $|w_1|=|f^{-1}w_2f_1|.$ Легко видеть, что если слова $w_1, w_2\in H$ циклически приведены и $w_1$ есть циклическая перестановка $w_2$, то существуют такие приведенные слова $C_1$ и $C_2$, что $w_1\equiv C_1\cdot C_2$ и $w_2\equiv C_2\cdot C_1$.

Если $w_1,w_2\in H\!\setminus\!\{1\}$ суть приведенные слова, то мы будем говорить, что $w_2$ есть \emph{подслово} в слове $w_1$ если $w_1\equiv f_1\cdot w_2\cdot f_2$ для некоторых приведенных слов $f_1,f_2\in H$. Если $A$ есть простое слово и $U$ есть подслово слова $A^k$, $k\in\mathbb{N}$, то слово $U$ мы будем называть  \emph{$A$-периодическим}.

Следующую теорему легко вывести из теоремы Файна-Вильфа (см. теорему 1 в работе~\cite{[FiWi65]}):

\begin{FiWiTh}Пусть группа $H=\zvezda\limits_{\mathfrak{a}\in \mathfrak{A}}\!H_\mathfrak{a}$ есть нетривиальное свободное произведение групп $H_\mathfrak{a}$ и пусть $w_1,w_2\in H$ суть циклически приведенные слова, причем $|w_1|_c,|w_2|_c\ge2$. Если $U$ есть подслово как в $w_1^{k_1}$, так и в $w_2^{k_2}$ для некоторых $k_1,k_2\in\mathbb{N}$ и $|U|\ge |w_1|_r+|w_2|_r-\text{{\rm НОД}}(|w_1|_r,|w_2|_r)$, то существуют такие приведенные слова $C_1,C_2\in H$, что $w_1 \equiv (C_1\cdot C_2)^{m_1}$ и $w_2 \equiv (C_2\cdot C_1)^{m_2}$ для некоторых $m_1,m_2\in\mathbb{N}$.
\end{FiWiTh}

Все утверждения следующей леммы известны.

\begin{Lemma}\label{MWL} Пусть группа $H=\zvezda\limits_{\mathfrak{a}\in \mathfrak{A}}\!H_\mathfrak{a}$ есть нетривиальное свободное произведение групп $H_\mathfrak{a}$ и пусть $A,B,T,U\in H$, причем $A$ и $B$ суть простые слова. Тогда:

\begin{enumerate}
\item[P1)]
если слово $U\equiv B$ есть подслово в $B^k$, $k\in\mathbb{N}$, то $B^k\equiv B^m\cdot U\cdot B^{k-m-1}$ для некоторого $m\in\{0,\dots,k-1\}$;

\item[P2)]
если $U$ есть приведенное слово и слово $T\equiv B\cdot U\cdot B$ есть подслово в $B^k$, $k\in\mathbb{N}$, то $U\equiv B^m$ для некоторого $m\in\{0,\dots,k-2\}$;

\item[P3)]
если слово $B^{-1}$ является $B$-периодическим, то существуют однозначно определенные приведенные слова $C_1,C_2\in H$ такие, что $B\equiv C_1\cdot C_2$, $B^{-1}\equiv C_2\cdot C_1$ и $C_1^2=C_2^2=1$;

\item[P4)]
если приведенное слово $U$ одновременно $A$-периодическое и $B$-периодическое, причем $|U|\ge |A|+|B|-1$, то существуют однозначно определенные приведенные слова $C_1,C_2\in H$ такие, что $A\equiv C_1\cdot C_2$ и $B\equiv C_2\cdot C_1$. Кроме того, если $U$ одновременно начинается и с подслова $A$ и с подслова $B$, то $A\equiv B$.
\end{enumerate}
\end{Lemma}

\begin{proof} $P1)$ Полагая $w_1\equiv w_2\equiv U \equiv B$ в условии теоремы Файна-Вильфа, мы заключаем, что найдутся такие приведенные слова $C_1,C_2\in H$, что $B\equiv C_1\cdot C_2$ и $B\equiv C_2\cdot C_1$. Если $C_1\equiv1$ или $C_2\equiv1$, то легко понять, что утверждение доказано. Пусть $C_1,C_2\not\equiv 1$, тогда т.к. $[C_1,C_2]=1$, то или $C_1,C_2\in H_\mathfrak{a}^g\!\setminus\!\{1\}$ для некоторого $\mathfrak{a}\in\mathfrak{A}$ и $g\in H$ или $C_1= f^{k_1}$, $C_2= f^{k_2}$ для некоторого гиперболического элемента $f\in H$ и $k_1,k_2\in\mathbb{N}$, причем $k_1$ и $k_2$ не могут иметь разные знаки, т.к. иначе между последним слогом $C_1$ и первым слогом $C_2$ будет происходить сокращение. В обоих случаях мы получаем, что $B$ не является простым словом.

$P2)$ Утверждение немедленно следует из $P1)$.

$P3)$ Полагая $w_1\equiv B$ и $w_2\equiv U \equiv B^{-1}$ в условии теоремы Файна-Вильфа, мы заключаем, что найдутся такие приведенные слова $C_1$ и $C_2$, что $B\equiv C_1\cdot C_2$ и $B^{-1}\equiv C_2\cdot C_1$. Значит $C_1\cdot C_2\equiv C_1^{-1}\cdot C_2^{-1}$, из чего следует (т.к. $|C_i|=|C_i^{-1}|$, $i=1,2$), что $C_i\equiv C_i^{-1}$, $i=1,2$.

Пусть $B\equiv C_1'\cdot C_2'$ и $B^{-1}\equiv C_2'\cdot C_1'$ для некоторых приведенных слов $C_1'$ и $C_2'$, отличных от $C_1$ и $C_2$. Без ограничения общности мы можем считать, что $|C_1'|>|C_1|$. Значит $C_1'\equiv C_1\cdot D$ для некоторого приведенного слова $D\not\equiv 1$. Т.к. $B^{-1} = B^{C_1}=D\cdot C_2'\cdot C_1$, то $B^{-D} = B^{-1}$, что невозможно (например, в силу того, что $C_H(B^{-1})=\left \langle B \right \rangle_\infty$ и $0<|D|_c<|B|_r$).

$P4)$ Существование слов $C_1$ и $C_2$ немедленно следует из теоремы Файна-Вильфа. Единственность доказывается также, как и в пункте~$P3)$, а доказательство последнего утверждения аналогично доказательству пункта~$P1)$.
\end{proof}

\begin{Lemma}\label{F2TLem1} Пусть группа $H=\zvezda\limits_{\mathfrak{a}\in \mathfrak{A}}\!H_\mathfrak{a}$ есть нетривиальное свободное произведение групп $H_\mathfrak{a}$, $A,B\in H$ суть простые слова, $f\in H$ есть приведенное слово и натуральные числа $k$ и $m$ таковы, что $k\ge \frac{|f|}{|A|} + \frac{|B|}{|A|} +1$ и $m\ge \frac{|f|}{|B|} + \frac{|A|}{|B|} +1$. Тогда, если при приведении слова $A^kfB^m$ затрагивается его первый или последний слог, то:
\begin{enumerate}
    \item[\sffamily $\mathrm{1)}$] $A\equiv C_1^{-1}\cdot C_2^{-1}$ и $B\equiv C_1\cdot C_2$ для некоторых приведенных слов $C_1,C_2\in H$;
    \item[\sffamily $\mathrm{2)}$] существуют такие $\alpha,\beta\in\mathbb{Z}$, что выполнены равенства:
    \begin{equation}\label{F2TLem1EQ1}
        f=A^{-\alpha}C_2B^{-\beta}=C_2B^{\alpha-\beta}=A^{-\alpha+\beta-1}C_1^{-1}.
    \end{equation}
 \end{enumerate}
\end{Lemma}
\begin{proof} Положим $\beta=\left \lceil\! \frac{|f|}{|B|}\! \right \rceil$, тогда ясно, что при приведении слова $fB^{\beta}$ его последний слог не затрагивается.
%Положим $\beta$ равным такому минимальному натуральному числу, что при приведении слова $fB^{\beta}$ его последний слог не затрагивался (ясно, что $\beta\le\left \lceil \frac{|f|}{|B|} \right \rceil\le m-2$).
Пусть $\hat{f}=fB^{\beta}$ есть приведенное слово и пусть при приведении слова $A^k\hat{f}B^{m-\beta}$ затрагивается его первый или последний слог. В таком случае несложно видеть, что (в силу выбора констант $k$ и $m$) слова $B^{m-\beta}$ и $A^{-k}$ имеют общее подслово длины не меньшей, чем $|A|+|B|-1$. Значит, в силу пункта~$P4)$ леммы~\ref{MWL}, найдутся такие приведенные слова $C_1$ и $C_2$, что $A^{-1}\equiv C_2\cdot C_1$ и $B\equiv C_1\cdot C_2$, более того, в таком случае $\hat{f}\equiv A^{-\alpha}\cdot C_2$, где $\alpha=\left \lfloor\! \frac{|\hat{f}|}{|A|}\! \right \rfloor$ (чтобы легче это понять, см. рис.~\ref{fig2}). Следовательно, мы имеем $f=A^{-\alpha}C_2B^{-\beta},$ откуда (учитывая, что $A\equiv C_1^{-1}\cdot C_2^{-1}$ и $B\equiv C_1\cdot C_2$) не составляет труда вывести оставшиеся два равенства из~(\ref{F2TLem1EQ1}).
\begin{figure}[h]
\begin{center}
\begin{tikzpicture}
%\draw[thick] (0,0) --++(2,0) node[midway,xshift=0ex,yshift=2ex] {$\hat{f}$}--+(0,0.15)--+(0,0)--++(1.5,0)node[midway,xshift=0ex,yshift=2ex]{$B$}--+(0,0.15)--+(0,0)--++(1.5,0)node[midway,xshift=0ex,yshift=2ex]{$B$};

\draw[rounded corners=10pt,thick] (1.65,0.15) .. controls (1.55,0.09) and (1.75,-0.09) .. (1.65,-0.15);
\draw[rounded corners=10pt,thick] (1.72,0.15) .. controls (1.62,0.09) and (1.82,-0.09) .. (1.72,-0.15);
\draw[rounded corners=10pt,thick] (7.8,0.15) .. controls (7.7,0.09) and (7.9,-0.09) .. (7.8,-0.15);
\draw[rounded corners=10pt,thick] (7.8,-1.55) .. controls (7.7,-1.61) and (7.9,-1.79) .. (7.8,-1.85);

\draw[very thick] (0,0)--++(0,0.15)--++(0,-0.15)--++(1.65,0);

\draw[very thick] (1.72,0)--++(1.88,0)node[very near start,xshift=0ex,yshift=2ex]{$\hat{f}$}
--++(0,0.15)--++(0,-0.15)--++(1.5,0)node[midway,xshift=0ex,yshift=2ex]{$B$}
--++(0,0.15)--++(0,-0.15)--++(1.5,0)node[midway,xshift=0ex,yshift=2ex]{$B$}
--++(0,0.15)--++(0,-0.15)--++(1.2,0)node[midway,xshift=0ex,yshift=2ex]{$B$};

\draw[] (1.72,0)--++(1.88,0)
--++(0,-0.1)--++(0,0.1)--++(0.9,0)node[midway,xshift=0ex,yshift=-2ex]{$C_1$}
--++(0,-0.15)--++(0,0.15)--++(0.6,0)node[midway,xshift=0ex,yshift=-2ex]{$C_2$}
--++(0,-0.15)--++(0,0.15)--++(0.9,0)node[midway,xshift=0ex,yshift=-2ex]{$C_1$}
--++(0,-0.15)--++(0,0.15)--++(0.6,0)node[midway,xshift=0ex,yshift=-2ex]{$C_2$}
--++(0,-0.15)--++(0,0.15)--++(0.9,0)
--++(0,-0.15);

%\draw[very thick] (0,0)--++(0,0.15)--++(0,-0.15)--++(1.65,0);

\draw[very thick] (0,-1.7)--++(0,0.15)--++(0,-0.15)--++(1.5,0)node[midway,xshift=0ex,yshift=2ex]{$A^{-1}$}
--++(0,0.15)--++(0,-0.15)--++(0.15,0);

\draw[rounded corners=10pt,thick] (1.65,-1.55) .. controls (1.55,-1.61) and (1.75,-1.79) .. (1.65,-1.85);
\draw[rounded corners=10pt,thick] (1.72,-1.55) .. controls (1.62,-1.61) and (1.82,-1.79) .. (1.72,-1.85);

\draw[very thick] (1.72,-1.7)--++(1.28,0)node[midway,xshift=0ex,yshift=2ex]{$A^{-1}$}
--++(0,0.15)--++(0,-0.15)--++(1.5,0)node[midway,xshift=0ex,yshift=2ex]{$A^{-1}$}
--++(0,0.15)--++(0,-0.15)--++(1.5,0)node[midway,xshift=0ex,yshift=2ex]{$A^{-1}$}
--++(0,0.15)--++(0,-0.15)--++(1.5,0)node[midway,xshift=0ex,yshift=2ex]{$A^{-1}$}
--++(0,0.15)--++(0,-0.15)--++(0.3,0);

\draw[] (1.72,-1.7)--++(1.88,0)
--++(-0.6,0)node[midway,xshift=0ex,yshift=-2ex]{$C_2$}
--++(0,-0.1)--++(0,0.1)--++(0.6,0)
--++(0,-0.1)--++(0,0.1)--++(0.9,0)node[midway,xshift=0ex,yshift=-2ex]{$C_1$}
--++(0,-0.15)--++(0,0.15)--++(0.6,0)node[midway,xshift=0ex,yshift=-2ex]{$C_2$}
--++(0,-0.15)--++(0,0.15)--++(0.9,0)node[midway,xshift=0ex,yshift=-2ex]{$C_1$}
--++(0,-0.15)--++(0,0.15)--++(0.6,0)node[midway,xshift=0ex,yshift=-2ex]{$C_2$}
--++(0,-0.15)--++(0,0.15)--++(0.9,0)
--++(0,-0.15);

\foreach \x in {0,3.6,4.5,5.1,6,6.6}
\draw[xshift=\x cm,loosely dotted] (0,0.5)--(0,-2.2);

%\foreach \x in {5,7}
%\draw[rounded corners=10pt,thick] (0,0.15) .. controls (-0.1,0.09) and (0.1,-0.09) .. (0,-0.15);

\end{tikzpicture}
\caption{Сверху --- приведенное слово $\hat{f}\cdot B^{m-\beta}\equiv\hat{f}\cdot(C_1\cdot C_2)^{m-\beta}$, снизу --- приведенное слово $A^{-k}\equiv(C_2\cdot C_1)^{k}.$}
%Если при приведении слова $A^{k}fB^m$ изменяется его первый или последний слог, то (в силу выбора констант $k$ и $m$) слова $A^{-k}$ и $\hat{f}\cdot B^{m-\beta}$ начинаются на общее %подслово длины не меньшей, чем $|\hat{f}|+|B|+|A|.$}
\label{fig2}
\end{center}
\end{figure}
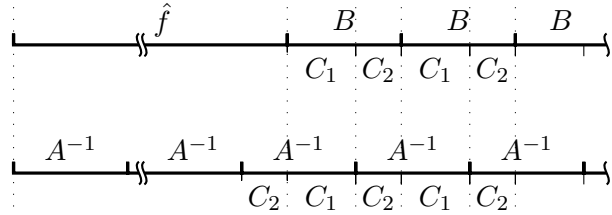

\end{proof}

\begin{Corollary}\label{LemBFLL} Пусть группа $H=\zvezda\limits_{\mathfrak{a}\in \mathfrak{A}}\!H_\mathfrak{a}$ есть нетривиальное свободное произведение групп $H_\mathfrak{a}$, $a\in H$ есть приведенное слово, $B\in H$ есть простое слово и натуральное число $k$ таково, что $k\ge \frac{|a|}{|B|}+2$. Тогда:
\begin{enumerate}
    \item[\sffamily $\mathrm{1)}$] если $a$ есть гиперболический элемент, то при приведении слова $B^kaB^k$ его первый и последний слог не затрагивается;

    \item[\sffamily $\mathrm{2)}$] если $[a,B]\ne1$, то при приведении слова $B^{-k}aB^k$ его первый и последний слог не затрагивается.
\end{enumerate}

\end{Corollary}
\begin{proof} 1) Если при приведении слова $B^kaB^k$ затрагивается его первый или последний слог, то из леммы~\ref{F2TLem1} легко следует, что $B\equiv C_1^{-1}\cdot C_2^{-1}$ и $B\equiv C_1\cdot C_2$ для некоторых приведенных слов $C_1,C_2\in H$ и что $a=C_2B^m$ для некоторого $m\in\mathbb{Z}$. Из равенства $C_1^{-1}\cdot C_2^{-1}\equiv C_1\cdot C_2$ следует, что $C_1^2=C_2^2=1$, но в таком случае $a= C_2(C_1C_2)^m$ есть элемент второго порядка, что противоречит условию.

2) Если при приведении слова $B^{-k}aB^k$ затрагивается его первый или последний слог, то из леммы~\ref{F2TLem1} легко следует, что $B^{-1}\equiv C_1^{-1}\cdot C_2^{-1}$ и $B\equiv C_1\cdot C_2$ для некоторых приведенных слов $C_1,C_2\in H$ и что $a=C_2B^m$ для некоторого $m\in\mathbb{Z}$. Учитывая, что $B$ есть просто слово, из равенств $B\equiv C_1\cdot C_2 \equiv C_2\cdot C_1$ следует (см. доказательство пункта~$P1)$ леммы~\ref{MWL}), что или $C_1\equiv 1$ и $C_2\equiv B$, или $C_1\equiv B$ и $C_2\equiv 1$. В обоих случаях мы получаем противоречие с условием (т.к. $[a,B]\ne1$ по условию).
\end{proof}

\begin{Lemma}\label{LemIvanAdd} Пусть группа $H=\zvezda\limits_{\mathfrak{a}\in \mathfrak{A}}\!H_\mathfrak{a}$ есть нетривиальное свободное произведение групп $H_\mathfrak{a}$, $A,B\in H$ суть простые слова, $f\in H$ есть приведенное слово, $|B|>|A|$ и $k\in\mathbb{N}$ таково, что $k|B|>|f|$. Тогда для любого $l\in\mathbb{N}$ при приведении слова $B^{6k+8}f^{-1}A^lfB^{6k+8}$ его первый и последний слог не затрагивается.
\end{Lemma}

\begin{proof}Положим $m=\left \lceil\! \frac{|B|}{|A|}\! \right \rceil$ и рассмотрим два случая.

Пусть $l\ge 2(mk+m+1)$. В силу того, что $A$ есть простое слово достаточно доказать, что при приведении слов $D_1\equiv B^{6k+8}f^{-1}A^{mk+m+1}$ и $D_2\equiv A^{mk+m+1}fB^{6k+8}$ их первые и последние слоги не затрагиваются. Докажем это для $D_2$ (для $D_1$ доказательство аналогичное).

Пусть при приведении слова $D_2$ затрагивается его первый или последний слог. Тогда, т.к. выполнены оценки:
\begin{equation}
mk+m+1>\frac{|f|}{|A|}+\frac{|B|}{|A|}+1\ \ \text{и} \ \ 6k+8> \frac{|f|}{|B|}+\frac{|A|}{|B|}+1,
\end{equation}
то из пункта~$1)$ леммы~\ref{F2TLem1} следует, что существуют такие приведенные слова $C_1,C_2\in H$, что $A\equiv C_1^{-1}\cdot C_2^{-1}$ и $B\equiv C_1\cdot C_2$. Следовательно, $|A| = |B|$, а это противоречит условию леммы.

Пусть $l< 2(mk+m+1)$. В таком случае мы имеем оценку:
\begin{align}\label{AddIEq1}
|f^{-1}A^{l}f|&<2|f|+2(m(k+1)+1)|A|\nonumber\\
&<2k|B|+2\biggl(\biggl(\frac{|B|}{|A|}+1\biggr)(k+1)+1\biggr)|A|\nonumber\\
&=2k|B|+2(|B|+|A|)(k+1)+2|A|\nonumber\\
&<(6k+6)|B|,
\end{align}
из которой немедленно следует, что $6k+8>|f^{-1}A^{l}f|/|B|+2$. Т.к. $f^{-1}A^{l}f$ есть гиперболический элемент, то утверждение следует из последнего неравенства и пункта~1) следствия~\ref{LemBFLL}.
\end{proof}

Все пункты следующей теоремы, на которую мы будем в дальнейшем ссылаться как на теорему Басса-Серра, хорошо известны (мы приводим доказательство последнего пункта, т.к. автору не удалось найти на него явной ссылки).

\begin{BSTh} Пусть группа $H=\zvezda\limits_{\mathfrak{a}\in \mathfrak{A}}\!H_\mathfrak{a}$ есть нетривиальное свободное произведение групп $H_\mathfrak{a}$, тогда группа $H$ может быть вложена в группу автоморфизмов некоторого ориентированного дерева $T$ таким образом, что при действии $H\acts T$:
\begin{enumerate}
    \item[\sffamily $\mathrm{1)}$] стабилизаторы вершин дерева $T$ имеют вид $H_\mathfrak{a}^g$ для некоторого $\mathfrak{a}\in \mathfrak{A}$ и $g\in H$;

    \item[\sffamily $\mathrm{2)}$] стабилизаторы ребер дерева $T$ тривиальны;

    \item[\sffamily $\mathrm{3)}$] автоморфизмы, соответствующие гиперболическим элементам группы $H$, имеют единственную инвариантную прямую;

    \item[\sffamily $\mathrm{4)}$] расстояние (относительно естественной метрики дерева $T$), на которое гиперболический элемент $h\in H$ сдвигает свою инвариантную прямую равно $|h|_c$.
\end{enumerate}
\end{BSTh}

\begin{proof} 4) В качестве множества вершин дерева $T$ возьмем множество левых смежных классов вида $gH_\mathfrak{a}$, $\mathfrak{a}\in \mathfrak{A}$, $g\in H$ и будем считать, что две (различные) вершины соединены ребром тогда и только тогда, когда соответствующие им смежные классы имеют нетривиальное пересечение. Действие элемента $h\in H$ на вершине $gH_\mathfrak{a}$ дерева $T$ определим как $h(gH_\mathfrak{a}) = hgH_\mathfrak{a}.$ Мы опускаем все необходимые здесь проверки корректности, т.к. оно относительно просты.

Гиперболический элемент $h\in H$ может быть единственным образом представлен как $h=g^{-1}\cdot A\circ g$, где $A$ есть циклически приведенное слово и $g$ есть приведенное слово. Пусть $A\equiv s_1s_2\cdots s_{n-1}s_n$, $n\ge2$ где $s_i\in H_{\mathfrak{a}_i}$ суть слоги слова $A$. В частности, $|h|_c=|A|=n.$

Определим вершины $v_{k,i}$, где $k\in\mathbb{Z}$ и $i=0,\ \dots,\ n-1$ дерева $T$ следующим образом:
$$v_{k,i}=\begin{cases}
g^{-1}A^k H_{\mathfrak{a}_n} & \text{ если } i=0;\\
g^{-1}A^ks_1\cdots s_i H_{\mathfrak{a}_i} & \text{ если } i=1,\ \dots,\ n-1.
\end{cases}$$
Легко видеть, что $v_{k,0}\cap v_{k,1}=\{g^{-1}A^k\},$ $v_{k,n-1}\cap v_{k+1,0}=\{g^{-1}A^ks_1\cdots s_{n-1}\}$ и $v_{k,i}\cap v_{k,i+1}=\{g^{-1}A^ks_1\cdots s_i\}$ при $i=1,\ \dots,\ n-2.$ Следовательно, вершины $v_{k,i}$ образуют прямую в дереве $T$. Т.к. очевидно, что $h(v_{k,i})=v_{k+1,i}$, то эта прямая инвариантна при действии элемента $h$. Осталось заметить, что расстояние (относительно естественной метрики дерева $T$) между вершинами $v_{k,i}$ и $v_{k+1,i}$ равно $n$ (действительно, отрезок, соединяющий эти вершины состоит из вершин $v_{k,i},\ v_{k,i+1},\ \dots,\ v_{k,n-1},\ v_{k+1,0},\ \dots,\ v_{k+1,i-1},\ v_{k+1,i}$).
\end{proof}

Далее в работе, когда мы говорим о действии группы $H=\zvezda\limits_{\mathfrak{a}\in \mathfrak{A}}\!H_\mathfrak{a}$ на дереве, мы всегда имеем в виду действие, описанное в теореме Басса-Серра, мы всегда рассматриваем левые действия и мы всегда вычисляем расстояния между вершинами дерева $T$ относительно естественной метрики этого дерева (для вершин $v_1$ и $v_2$ мы обозначаем это расстояние как $d(v_1,v_2)$).

\begin{Lemma}\label{GaoT1} Пусть группа $H=\zvezda\limits_{\mathfrak{a}\in \mathfrak{A}}\!H_\mathfrak{a}$, являющаяся нетривиальным свободным произведением групп $H_\mathfrak{a}$, действует на дереве $T$ и пусть $X_1,X_2\in H$ суть гиперболические элементы, $\alpha_1$ и $\alpha_2$ соответствующие этим элементам инвариантные прямые и $I = \alpha_1\cap \alpha_2$ --- общая часть этих прямых. Тогда, если $[X_1,X_2]\ne 1$, то или прямые $\alpha_1$ и $\alpha_2$ не пересекаются или $I$ есть отрезок и справедлива оценка:
\begin{equation}\label{GaoT1Eq1}
|I|\le |X_1|_r+|X_2|_r-1,
\end{equation}
где $|I|$ есть длина отрезка $I$.
\end{Lemma}

\begin{proof} Пусть $X_1= f_1^{k_1}$ и $X_2= f_2^{k_2}$ для некоторых элементов $f_1,f_2\in H$, не являющихся собственной степенью в $H$, и $k_1,k_2\in\mathbb{N}$. Из теоремы Басса-Серра несложно вывести, что два гиперболических элемента группы $H$ коммутируют тогда и только тогда, когда их инвариантные прямые совпадают. В частности, из этого следует, что $[f_1,f_2]\ne 1$ и инвариантные прямые элементов $f_1$ и $f_2$ совпадают с инвариантными прямыми элементов $X_1$ и $X_2$ соответственно. Значит, т.к. $T$ --- дерево, $I$ есть или пустое множество или отрезок (одну вершину мы считаем отрезком) или луч.

Если длина $I$ не меньше чем $|f_1|_c+|f_2|_c=|X_1|_r+|X_2|_r$, то ясно, что найдутся такие вершины $v_1$ и $v_2$ дерева $T$, что $v_1,v_2\in I$ и $d(v_1,v_2)=|f_1|_c+|f_2|_c$. Пусть элементы $f_1$ и $f_2$ осуществляют сдвиги своих инвариантных прямых по направлению от $v_1$ к $v_2$. В таком случае легко видеть, что $f_1f_2(v_1)=v_2$ и $f_1^{-1}f_2^{-1}(v_2)=v_1$. Значит, $[f_1,f_2](v_1)=f_1^{-1}f_2^{-1}f_1f_2(v_1)=v_1$, т.е. коммутатор $[f_1,f_2]$ принадлежит стабилизатору вершины $v_1$, следовательно, $[f_1,f_2]\in H_\mathfrak{a}^g\!\setminus\!\{1\}$ для некоторого $\mathfrak{a}\in\mathfrak{A}$ и $g\in H$. Но, в соответствии с пунктом~$C1)$ леммы~\ref{MCL}, из этого следует, что $f_1,f_2\in H_\mathfrak{a}^g\!\setminus\!\{1\}$, а это противоречит тому, что элементы $X_1= f_1^{k_1}$ и $X_2= f_2^{k_2}$ гиперболические. Оставшиеся три случая разбираются аналогично, например, если элемент $f_2$ осуществляет сдвиг своей инвариантной прямой по направлению от $v_2$ к $v_1$ (а элемент $f_1$ от $v_1$ к $v_2$), то вместо действия элемента $[f_1,f_2]$ мы рассмотрим действие элемента $[f_1,f_2^{-1}]$.
\end{proof}

{\Observation {\rm Опираясь на теорему Файна-Вильфа можно доказать, что в условиях леммы~\ref{GaoT1} выполнена более точная оценка: $$|I|\le|X_1|_r+|X_2|_r - \text{НОД}(|X_1|_r,|X_2|_r).$$Мы не приводим это доказательство, т.к. оно несколько сложнее и для целей данной работы достаточно оценки~(\ref{GaoT1Eq1}).}}

\begin{Lemma}\label{AddIIICILL}
Пусть группа $H=\zvezda\limits_{\mathfrak{a}\in \mathfrak{A}}\!H_\mathfrak{a}$, являющаяся нетривиальным свободным произведением групп $H_\mathfrak{a}$, действует на дереве $T$ и пусть $X_1,X_2\in H$ суть некоммутирующие гиперболические элементы. Тогда для любого $k\in\mathbb{N}$, $k\ge 2$ справедлива оценка:
\begin{equation}\label{AddIIICILLEq1}
|[X_1^k,X_2^k]|_c\ge 2(k|X_1|_c+k|X_2|_c-|I|),
\end{equation}
где $|I|$ есть длина общего отрезка инвариантных прямых элементов $X_1$ и $X_2$ (если прямые не пересекаются, то мы полагаем $|I|=0$).
\end{Lemma}

Перед тем как перейти к доказательству, отметим, что приведенная оценка выполняется и при $k=1$, но мы опускаем этот случай в силу того, что он требует отдельного рассмотрения и не используется в наших последующих рассуждениях.

\begin{proof} Наше доказательство будет основано на оценке снизу расстояния $d\big(v,[X_1^k,X_2^k](v)\big)$, на которое элемент $[X_1^k,X_2^k]$ может сдвинуть вершину $v$ дерева $T$.

Пусть $\alpha$ и $\beta$ суть инвариантные прямые элементов $X_1$ и $X_2$ соответственно. Мы опускаем доказательство случая, при котором прямы $\alpha$ и $\beta$ не пересекаются в силу того, что, во-первых, это доказательство в достаточной степени аналогично разобранному ниже случаю и, во-вторых, в силу того, что для этого случая оценка~(\ref{AddIIICILLEq1}) может быть улучшена до оценки: $$|[X_1^k,X_2^k]|_c\ge 2(k|X_1|_c+k|X_2|_c+2|J|),$$ где $|J|$ есть длина отрезка, соединяющего прямые $\alpha$ и $\beta$.

Итак, пусть $|I|\ge 1$. Покажем, что без потери общности при рассмотрение действия элемента $[X_1^k,X_2^k]=X_1^{-k}X_2^{-k}X_1^kX_2^k$ на дереве $T$ мы можем считать, что $|X_2|_c\ge |X_1|_c$ и что элементы $X_1$ и $X_2$ осуществляют сдвиги своих инвариантных прямых разнонаправленно (относительно отрезка $I$). Действительно, если это не так, то вместо слова $X_1^{-k}X_2^{-k}X_1^kX_2^k$ рассмотрим слово:
\begin{align*}
&X_2^{-k}X_1^{k}X_2^{k}X_1^{-k}\ \ \text{если сдвиги сонаправлены и}\ |X_2|_c\le |X_1|_c;\\
&X_1^{k}X_2^{-k}X_1^{-k}X_2^{k}\ \ \text{если сдвиги сонаправлены и}\ |X_2|_c\ge |X_1|_c;\\
&X_2^{-k}X_1^{-k}X_2^kX_1^{k}\ \ \text{если сдвиги разнонаправлены и}\ |X_2|_c\le |X_1|_c.
\end{align*}
Т.к. данные слова имеют одинаковую центральную длину (каждое из этих слов может быть получено из $X_1^{-k}X_2^{-k}X_1^kX_2^k$ подходящим сопряжением и/или взятием обратного), то достаточно оценить центральную дину любого из них.

Несложно показать, что минимум величины $d\big(v,[X_1^k,X_2^k](v)\big)$ достигается на вершинах, лежащих на инвариантной прямой элемента $[X_1^k,X_2^k]$. Следовательно, достаточно оценить расстояние, на которое элемент $[X_1^k,X_2^k]$ сдвигает свою инвариантную прямую. Т.к. $|X_2|_c\ge|X_1|_c$, $k\ge2$ и выполнена оценка $|I|\le |X_1|_r+|X_2|_r-1$ (см. лемму~\ref{GaoT1}), то элемент $X_2^k$ сдвигает свою инвариантную прямую $\beta$ на расстояние $k|X_2|_c$, большее $|I|$. Возможны два случая: $k|X_1|_c<|I|$ и $k|X_1|_c\ge|I|$. Мы рассмотрим лишь первый из них, т.к. второй, по сути, аналогичен.

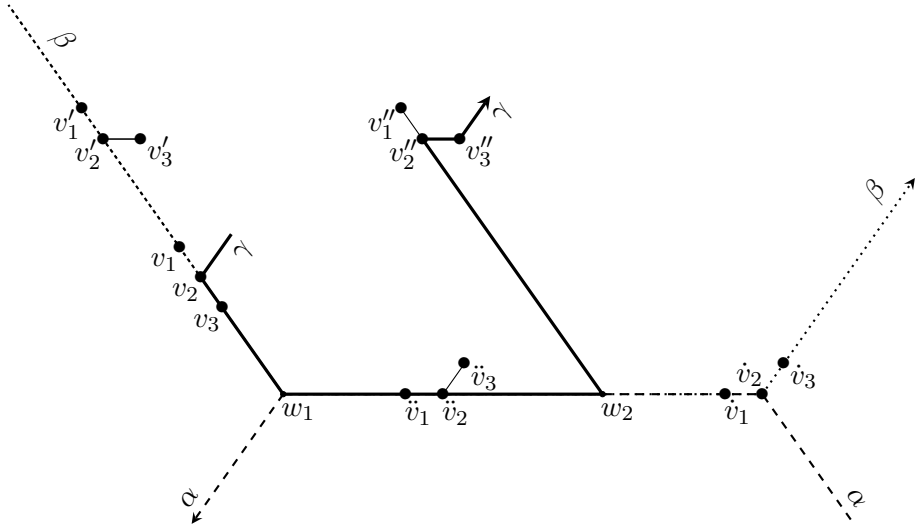
\begin{figure}[h]
\begin{center}
\begin{tikzpicture}[scale=0.7,>=stealth]

\draw[dotted,thick,->] (0,0)--++(125:9)--++(-55:9)node[very near start,sloped,above] {$\beta$}--++(9,0)--++(55:5)node[very near end,sloped,above]{$\beta$};

\draw (0,0)++(125:2)node[xshift=-1.2ex,yshift=-1.2ex]{$v_3$}++(125:0.7)node[xshift=-1.2ex,yshift=-1.2ex] {$v_2$}++(125:0.7)node[xshift=-1.2ex,yshift=-1.2ex] {$v_1$};
\draw (0,0)++(125:2)node{\textbullet}++(125:0.7)node {\textbullet}++(125:0.7)node {\textbullet};

\draw (0,0)++(8.3,0)node[xshift=1.0ex,yshift=-1.4ex]{$\dot{v}_1$}++(0.7,0)node[xshift=-1.0ex,yshift=1.6ex] {$\dot{v}_2$}++(55:0.7)node[xshift=1.6ex,yshift=-1ex] {$\dot{v}_3$};
\draw (0,0)++(8.3,0)node[]{\textbullet}++(0.7,0)node[]{\textbullet}++(55:0.7)node[]{\textbullet};

\draw (0,0)++(2.3,0)node[xshift=1.0ex,yshift=-1.5ex]{$\ddot{v}_1$}++(0.7,0)node[xshift=1.0ex,yshift=-1.5ex]{$\ddot{v}_2$}++(55:0.7)node[xshift=1.6ex,yshift=-1.1ex]{$\ddot{v}_3$};
\draw (0,0)++(2.3,0)node[]{\textbullet}++(0.7,0)node[]{\textbullet}--++(55:0.7)node[]{\textbullet};
%\draw (0,0)++(2.6,0)node[]{}++(0.7,0)node[]{}--++(55:0.7)node[]{};

\draw (0,0)++(125:6.6)node[xshift=-1.2ex,yshift=-1.2ex]{$v_1'$}++(-55:0.7)node[xshift=-1.2ex,yshift=-1.2ex]{$v_2'$}--++(0.7,0)node[xshift=1.6ex,yshift=-1ex]{$v_3'$};
\draw (0,0)++(125:6.6)node[]{\textbullet}++(-55:0.7)node[]{\textbullet}--++(0.7,0)node[]{\textbullet};

\draw (0,0)++(6,0)--++(125:6.6)node[xshift=-1.5ex,yshift=-1.2ex]{$v_1''$}++(-55:0.7)node[xshift=-1.5ex,yshift=-1.2ex]{$v_2''$}--++(0.7,0)node[xshift=1.6ex,yshift=-1ex]{$v_3''$};
\draw (0,0)++(6,0)++(125:6.6)node[]{\textbullet}++(-55:0.7)node[]{\textbullet}++(0.7,0)node[]{\textbullet};

\draw[dashed,thick] (0,0)++(-125:3)++(55:3)--++(9,0)--++(-55:3)node[very near end,sloped,above] {$\alpha$};
\draw[dashed,thick,->] (0,0)--++(-125:3)node[very near end,sloped,above] {$\alpha$};

\draw[very thick,->] (0,0)--++(125:2.7)++(55:1)--++(-125:1)node[very near start,sloped,below] {$\gamma$}++(-55:2.7)--++(6,0)--++(125:5.9)--++(0.7,0)--++(55:1)node[very near end,sloped,below] {$\gamma$};

\filldraw [black] (0,0) circle (1.4pt) node[xshift=1.2ex,yshift=-1.4ex] {$w_1$};
\filldraw [black] (6,0) circle (1.4pt) node[xshift=1.2ex,yshift=-1.4ex] {$w_2$};
\end{tikzpicture}
\caption{$[w_1,\dot{v}_2]=I$ --- общий отрезок инвариантных прямых $\alpha$ и $\beta$ элементов $X_2$ и $X_1$ соответственно. Элемент $X_2^k$ сдвигает прямую $\beta$ слева на право; элемент $X_1^k$ сдвигает прямую $\alpha$ справа на лево, причем $k|X_1|_c\le |I|$.}
\label{fig1}
\end{center}
\end{figure}

Пусть вершины $v_1$, $v_2$ и $v_3$ лежат на прямой $\beta$ и находятся на расстоянии $k|X_2|_c-|I|+1$, $k|X_2|_c-|I|$ и $k|X_2|_c-|I|-1$ от отрезка $I$ соответственно, причем при действии $X_2^k$ эти вершины смещаются по направлению к $I$ (см. рис.~\ref{fig1}). При действии элемента $X_2^k$ вершины $v_1$, $v_2$ и $v_3$ переходят соответственно в вершины $\dot{v}_1$, $\dot{v}_2$ и $\dot{v}_3$, причем $\dot{v}_1,\dot{v}_2\in I$ и $\dot{v}_3\not\in \alpha$. Далее, при действии элемента $X_1^{k}$ вершины $\dot{v}_1$, $\dot{v}_2$ и $\dot{v}_3$ переходят соответственно в вершины $\ddot{v}_1$, $\ddot{v}_2$ и $\ddot{v}_3$, причем $\ddot{v}_1,\ddot{v}_2\in I$ (здесь используется то, что $k|X_1|_c<|I|$) и $\ddot{v}_3\not\in \alpha\cup\beta$. Далее, при действии элемента $X_2^{-k}$ вершины $\ddot{v}_1$, $\ddot{v}_2$ и $\ddot{v}_3$ переходят соответственно в вершины $v_1'$, $v_2'$ и $v_3'$, причем $v_1',v_2'\in\beta$, $v_1',v_2'\not\in\alpha$ и $v_3'\not\in\alpha\cup\beta$. Наконец, элемент $X_1^{-k}$ переводит отрезок $[v_1',w_1]$ в отрезок $[v_1'',w_2]$, а вершины $v_1'$, $v_2'$ и $v_3'$ соответственно в вершины $v_1''$, $v_2''$ и $v_3''$. Вычислим расстояние $d\big(v_2,[X_1^k,X_2^k](v_2)\big)=d(v_2,v_2'')$:
\begin{align}\label{InLiPi1}
d(v_2,v_2'')&=d(v_2,w_1)+d(w_1,w_2)+d(w_2,v_2'')\nonumber\\
&=k|X_2|_c-|I|+k|X_1|_c+d(w_1,v_2')\nonumber\\
&=k|X_2|_c-|I|+k|X_1|_c+k|X_2|_c-(|I|-k|X_1|_c)=2\big(k|X_1|_c+k|X_2|_c-|I|\big).
\end{align}
Аналогичные вычисления показывают, что $d\big(v_3,[X_1^k,X_2^k](v_3)\big)=d\big(v_2,[X_1^k,X_2^k](v_2)\big)$, тогда как $d\big(v_1,[X_1^k,X_2^k](v_1)\big)=d\big(v_2,[X_1^k,X_2^k](v_2)\big)+2$. Следовательно, несложно видеть, что инвариантная прямая $\gamma$ элемента $[X_1^k,X_2^k]$ (отмечена на рис.~\ref{fig1} сплошной линией) проходит через вершины $v_2$, $v_3$, $v_2''$, $v_3''$ и не проходит через вершины $v_1$ и $v_1''$ (отметим также, что $\gamma\cap\alpha =[w_1,w_2]$ и $\gamma\cap\beta =[v_2,w_2]$). Значит из~(\ref{InLiPi1}) следует необходимая нам оценка.
\end{proof}

\begin{Lemma}\label{MCL} Пусть группа $H=\zvezda\limits_{\mathfrak{a}\in \mathfrak{A}}\!H_\mathfrak{a}$ есть нетривиальное свободное произведение групп $H_\mathfrak{a}$ и $M_H$ есть множество всех гиперболических элементов группы $H$, тогда:
\begin{enumerate}
\item[C1)]
если $X_1,X_2\in H$ и $[X_1,X_2]\in H_\mathfrak{a}^g\!\setminus\!\{1\}$ для некоторого $\mathfrak{a}\in \mathfrak{A}$ и $g\in H$, то $X_1,X_2\in H_\mathfrak{a}^g\!\setminus\!\{1\};$

\item[C2)]
если $X_1,X_2\in M_H$ и $[X_1,X_2]=1$, то $|X_1|_r=|X_2|_r$;

\item[C3)]
если $X_1,X_2\in M_H$, $[X_1,X_2]\ne1$ и $n\ge 2$, то выполнена оценка:
\disp{\it
$|[X_1^n,X_2^n]|_c>2(n-1)(|X_1|_c+|X_2|_c);$
}
\item[C4)]
если $X_1,X_2\in M_H$, $[X_1,X_2]\ne1$ и $n\ge 2$, то существует такое $N\in\mathbb{N}$, что при $m>N$ слово $[X_1^n,X_2^{nm}]$ не является в $H$ степенью, большей, чем два;

\item[C5)]
если $X_1,X_2\in M_H$, $[X_1,X_2]\ne1$ и $n\ge 2$, то:
\disp{\it
для любого $C\in\mathbb{N}$ существует такое $N\in\mathbb{N}$, что $|[X_1^n,X_2^{nm}]|_r>C$ при $m\ge N$.
}
\end{enumerate}
\end{Lemma}

Прежде чем приступить к доказательству данной леммы сформулируем следующую версию теоремы Куроша о строении подгрупп свободных произведений групп (далее в тексте мы будем ссылаться на эту теорему просто как на теорему Куроша):

\begin{KurTh} Если $G=\zvezda\limits_{\alpha}\!A_\alpha$ и $H$ есть произвольная подгруппа группы $G$, то существует разложение подгруппы $H$,
$$H=F\ast\!\zvezda\limits_{\beta}\!B_\beta,$$
такое, что:
 \begin{enumerate}
    \item[\sffamily $\mathrm{1)}$] множитель $F$ есть свободная подгруппа, а каждый множитель $B_\beta$ сопряжен в $G$ с некоторой подгруппой одного из множителей $A_\alpha$;
    \item[\sffamily $\mathrm{2)}$] если $g$ есть произвольный элемент из $G$, $A_\alpha$ --- произвольный свободный сомножитель исходного разложения группы $G$ и если пересечение $D=H\cap A_\alpha^g$ нетривиально, то в указанном разложении $H$ найдется свободный множитель, сопряженный с $D$ внутри $H$.
 \end{enumerate}
\end{KurTh}

\begin{proof} $C1)$ Рассмотрим подгруппу $S=\left \langle X_1,X_2 \right \rangle$ группы $H$. Из условия ясно, что группа $S$ не циклическая, значит, опираясь на теорему Куроша и теорему Грушко~\cite{[Gru40]}, мы можем заключить, что или $S$ есть подгруппа в $H_\mathfrak{a}^g$ (и утверждение доказано) или $S=\left \langle f_1 \right \rangle_\infty\ast\left \langle f_2 \right \rangle_\infty$ или $S=\left \langle f_1 \right \rangle_\infty\ast\left \langle h_2 \right \rangle$ или $S=\left \langle h_1 \right \rangle\ast\left \langle h_2 \right \rangle,$ где $f_1,f_2\in M_H$ и $h_1,h_2\in H\!\setminus\!M_H$. Положим $u = [X_1,X_2]$, тогда из пункта~2) теоремы Куроша следует, что если $S=\left \langle f_1 \right \rangle_\infty\ast\left \langle f_2 \right \rangle_\infty$, то ситуация $u\in H\!\setminus\! M_H$ невозможна, если $S=\left \langle f_1 \right \rangle_\infty\ast\left \langle h_2 \right \rangle$, то $u^s\in\left \langle h_2 \right \rangle$ для некоторого $s\in S$, если же $S=\left \langle h_1 \right \rangle\ast\left \langle h_2 \right \rangle$, то $u^s\in\left \langle h_1 \right \rangle$ или $u^s\in\left \langle h_2 \right \rangle$ для некоторого $s\in S$. Заметим, что если $S=\left \langle f_1 \right \rangle_\infty\ast\left \langle h_2 \right \rangle$ или $S=\left \langle h_1 \right \rangle\ast\left \langle h_2 \right \rangle,$ то коммутант $S'$ совпадает с декартовой подгруппой\footnote[$\dagger$]{\emph{Декартовой подгруппой} группы вида $H=\zvezda\limits_{\mathfrak{a}\in \mathfrak{A}}\!H_\mathfrak{a}$ мы называем ядро естественного эпиморфизма $\zvezda\limits_{\mathfrak{a}\in \mathfrak{A}}\!H_\mathfrak{a}\twoheadrightarrow\! \krest\limits_{\mathfrak{a}\in \mathfrak{A}}\!\!H_\mathfrak{a}.$} группы $S$. Хорошо известно, что декартова подгруппа тривиально пересекает свободные сомножители, следовательно, ситуация $u^s\in\left \langle h_i \right \rangle$, $i=1,2$ невозможна, т.к. $u^s\in S'\!\setminus\!\{1\}$.

$C2)$ Если $X_1,X_2\in M_H$ и $[X_1,X_2]=1$, то, в силу теоремы Куроша, $\left \langle X_1,X_2 \right \rangle=\left \langle f \right \rangle_\infty$ для некоторого $f\in M_H$. Значит $X_1=f^{k_1}$, $X_2=f^{k_2}$ для некоторых $k_1,k_2\in\mathbb{Z}\!\setminus\!\{0\}$ и мы имеем $|X_i|_r=|f^{k_i}|_r=|f|_r$, $i=1,2.$

$C3)$ Нужная оценка немедленно следует из оценок $|I|\le |X_1|_r+|X_2|_r-1$ и $|[X_1^k,X_2^k]|_c\ge 2(k|X_1|_c+k|X_2|_c-|I|)$, установленных в леммах~\ref{GaoT1} и~\ref{AddIIICILL}, соответственно.

$C4)$ %Учитывая пункт $C3)$, достаточно показать, что при всех достаточно больших значениях $m$ выполнено неравенство $2[X_1^n,X_2^{nm}]_r\ge[X_1^n,X_2^{nm}]_c$.
Без ограничения общности мы можем считать, что $X_2$ не является собственной степенью в $H$, т.е. $X_2=f^{-1}\cdot B \circ f$ для некоторого простого слова $B$ и приведенного слова $f$. Пусть $a = f X_1^{n}f^{-1}$ есть приведенное слово, положим $k=\left \lceil\! \frac{|a|}{|B|}\! \right \rceil+2$ и будем считать, что $m>2k$. Сопрягая слово $[X_1^n,X_2^{nm}]$ элементом $f^{-1}B^{-k}$, мы имеем:
\begin{equation}\label{AddIIIMCLC4EQ1}
[X_1^n,X_2^{nm}]^{f^{-1}B^{-k}}=E_1B^{-nm+2k}E_2B^{nm-2k},
\end{equation}
где $E_1=B^{k}a^{-1}B^{-k}$ и $E_2= B^{-k}aB^k$ суть приведенные слова. Учитывая наш выбор $k$, условие $[a,B]\ne1$, пункт~2) следствия~\ref{LemBFLL}, а также принимая во внимание то, что $nm-2k>0$ (т.к. $n\in\mathbb{N}$ и $m>2k$), несложно видеть, что слово $g\equiv E_1 B^{-nm+2k} E_2 B^{nm-2k}$ циклически приведено. Значит:
%
%В силу выбора $k$ и т.к. $[a,B]\ne1$ и $m-2k>0$ из пункта 2) следствия~\ref{LemBFLL} легко вывести, что
\begin{equation}\label{AddIIIMCLC4EQ2}
|[X_1^n,X_2^{nm}]|_c=|g|=|E_1\cdot B^{-nm+2k}\cdot E_2\cdot B^{nm-2k}|=2(mn-2k)|B|+|E_1|+|E_2|.
\end{equation}
Т.к. по построению величины $|E_1|$ и $|E_2|$ не зависят от выбора $m$, то из~(\ref{AddIIIMCLC4EQ2}) следует, что при достаточно большом значение $m$ (конкретно, при $m>N= \big\lceil  \big(2k+3+(|E_1|+|E_2|)/|B|\big)/n \big\rceil$), слово $g$ заканчивается на подслово $B^{nm-2k}$, длина которого больше, чем $1/3|g|+|B|$. Фиксируем такое $m$.

Пусть оказалось, что слово $[X_1^n,X_2^{nm}]$ есть собственная степень, большая двух. В таком случае найдется такое простое слово $C$, что $g\equiv C^l$ и $l\ge 3$. Значит слово $g$ заканчивается на одновременно $B$- и $C$-периодическое подслово $B^{nm-2k}$, длина которого больше, чем $|C|+|B|$. Следовательно, в силу пункта~$P4)$ леммы~\ref{MWL}, мы имеем $B\equiv C$, т.е. мы имеем равенство:
\begin{equation}\label{AddIIIMCLC4EQ3}
E_1\cdot B^{-nm+2k}\cdot E_2\cdot B^{nm-2k}\equiv B^l.
\end{equation}
Т.к. $nm-2k>0$, то из этого равенства следует, что слово $B^{-1}$ является $B$-периодическим, значит, в силу пункта~$P3)$ леммы~\ref{MWL}, существуют такие приведенные слова $C_1,C_2\in H$, что $B\equiv C_1\cdot C_2$, $B^{-1}\equiv C_2\cdot C_1$ и $C_1^2=C_2^2=1$. Т.е. равенство~(\ref{AddIIIMCLC4EQ3}) имеет вид:
$$E_1\cdot(C_2\cdot C_1)^{nm-2k}\cdot E_2\cdot(C_1\cdot C_2)^{nm-2k}\equiv (C_1\cdot C_2)^{l},$$
из которого легко следует, что $E_1\equiv B^{l_1}C_1$ и $E_2\equiv C_2B^{l_2}$ для некоторых $l_1,l_2\in\mathbb{N}\cup\{0\}$, в частности $E_1^2=E_2^2=1$. Но это означает, что мы имеем $X_1^{2n}=a^{2f}=E_2^{2B^{-k}f}=1$, что противоречит условию $X_1\in M_H$.

$C5)$ Пусть $N$ настолько велико, что для $m>N$ выполнен пункт~$C4)$, тогда:
\begin{align*}
|[X_1^n,X_2^{nm}]|_r&\ge\frac{1}{2}|[X_1^n,X_2^{nm}]|_c\\
&>(n-1)\big(|X_1|_c+|X_2^m|_c\big)\\
&>(n-1)m|X_2|_c\to\infty\ \ \text{при}\ \ m\to\infty,
\end{align*}
где при первом и втором переходе мы воспользовались $C4)$ и $C3)$, соответственно.
\end{proof}

{\Observation {\rm В утверждении пункта~$C4)$ случай, когда слово $[X_1^n,X_2^{nm}]$ является квадратом в $H$ не может быть исключен. Действительно, для любых $k_1,k_2\in\mathbb{N}$ в группе $H=\left \langle a \right \rangle_2\ast\left \langle b \right \rangle_2\ast\left \langle c \right \rangle_2$ выполнено:$$\big[(ab)^{k_1},(bc)^{k_2}\big]=\big((ba)^{k_1}(cb)^{k_2-1}c \big )^2.$$}}

\section{Доказательство Леммы о словах, часть I}

Основной целью данного раздела является доказательство следующей леммы.
\begin{Lemma}\label{LemMLPartII} Пусть группа $H=\zvezda\limits_{\mathfrak{a}\in \mathfrak{A}}\!H_\mathfrak{a}$ есть нетривиальное свободное произведение групп $H_\mathfrak{a}$ и пусть $\underline{x}=(x_0,\dots,x_{n-1})$, $n\ge2$ есть такой кортеж гиперболических элементов группы $H$, что $[x_i,x_{i+1}]\ne1$ при $i=0,\ \dots,\ n-2$. Тогда:
\begin{enumerate}
     \item[\sffamily $\mathrm{a)}$] существует слово $T_{\underline{x}}(z_0,\dots,z_{n-1})\in F_n(z_0,\dots,z_{n-1})'$, вообще говоря, зависящее от кортежа $\underline{x}$ такое, что если для некоторого кортежа $(y_0,\dots,y_{n-1})$ гиперболических элементов группы $H$ выполнено равенство:
\begin{equation*}
T_{\underline{x}}(x_0,\dots,x_{n-1})=T_{\underline{x}}(y_0,\dots,y_{n-1}),
\end{equation*}
то найдется такой элемент $v\in \left \langle T_{\underline{x}}(x_0,\dots,x_{n-1}) \right \rangle_\infty$, что $x_i = y_i^v$ при $i=0,\ \dots,\ n-1$;

     \item[\sffamily $\mathrm{b)}$] существуют слова $T_{\underline{x}}'(z_0,\dots,z_{n-1}), T_{\underline{x}}''(z_0,\dots,z_{n-1})\in F_n(z_0,\dots,z_{n-1})'$, вообще говоря, зависящие от кортежа $\underline{x}$ такие, что если для некоторого кортежа $(y_0,\dots,y_{n-1})$ гиперболических элементов группы $H$ одновременно выполнены равенства:
\begin{equation}\label{LemMLPartIIaddTSS}
T_{\underline{x}}'(x_0,\dots,x_{n-1})=T_{\underline{x}}'(y_0,\dots,y_{n-1})\ \ \text{и}\ \ T_{\underline{x}}''(x_0,\dots,x_{n-1})=T_{\underline{x}}''(y_0,\dots,y_{n-1}),
\end{equation}
то $x_i = y_i$ при $i=0,\ \dots,\ n-1$.
 \end{enumerate}
\end{Lemma}
%\begin{Lemma}\label{LemMLPartII} Пусть группа $H=\zvezda\limits_{\mathfrak{a}\in \mathfrak{A}}\!H_\mathfrak{a}$ есть нетривиальное свободное произведение групп $H_\mathfrak{a}$ и %пусть $\underline{x}=(x_0,\dots,x_{n-1})$ есть такой кортеж гиперболических элементов группы $H$, что $[x_i,x_{i+1}]\ne1$ при $i=0,\ \dots,\ n-2$. Тогда существует слово %$T_{\underline{x}}(z_0,\dots,z_{n-1})\in F_n(z_0,\dots,z_{n-1})'$, вообще говоря, зависящее от кортежа $\underline{x}$ такое, что если для некоторого кортежа $(y_0,\dots,y_{n-1})$ %гиперболических элементов $y_i\in H$ выполнено равенство:
%\begin{equation*}
%T_{\underline{x}}(x_0,\dots,x_{n-1})=T_{\underline{x}}(y_0,\dots,y_{n-1}),
%\end{equation*}
%то найдется такой элемент $v\in \left \langle T_{\underline{x}}(x_0,\dots,x_{n-1}) \right \rangle_\infty$, что $x_i = y_i^v$ при $i=0,\ \dots,\ n-1.$
%\end{Lemma}

Предварительно мы установим ряд вспомогательных лемм. Доказательства леммы~\ref{LemIvanML1} и пункта~$V2)$ леммы~\ref{LemIvanPL} во многом основаны на идеях доказательств соответственно леммы 4 и леммы 1, изложенных в работе С.В. Иванова~\cite{[Ivanov]}.

\begin{Lemma}\label{LemIvanPL} Пусть группа $H=\zvezda\limits_{\mathfrak{a}\in \mathfrak{A}}\!H_\mathfrak{a}$ есть нетривиальное свободное произведение групп $H_\mathfrak{a}$ и пусть $X_1$, $X_2$ суть некоммутирующие гиперболические элементы группы $H$. Тогда:
%$M_H$ есть множество всех гиперболических элементов группы $H$, $X_1,X_2\in M_H$ и $[X_1,X_2]\ne 1$, то:
\begin{enumerate}
\item[V1)]
для любого $k\in\mathbb{N}$ слово $[X_1^{10},X_2^{10k}]$ сопряжено в $H$ слову вида $B(X_1,X_2,k)^{l(X_1,X_2,k)}$, где $l(X_1,X_2,k)\in\mathbb{N}$ и $B(X_1,X_2,k)$ есть простое слово\footnote[$\dagger$]{Мы используем в формулировке данной леммы громоздкие обозначения типа $B(X_1,X_2,k)$ и $l(X_1,X_2,k)$, явно указывающие на слова $X_1$, $X_2$ и константу $k$, для которых эти объекты определяются, чтобы избежать путаницы в дальнейшем.};

\item[V2)]
для любого $k\in\mathbb{N}$ существует такой элемент $s(X_1,X_2,k)\in H$, что если положить {\rm$Y_1(X_1,X_2,k)\overset{\textrm{def}}{=}X_1^{s(X_1,X_2,k)}$} и {\rm$Y_2(X_1,X_2,k)\overset{\textrm{def}}{=}X_2^{s(X_1,X_2,k)},$} то будет выполнено:
\begin{enumerate}
\item[1)]
$\big[Y_1(X_1,X_2,k)^{10},Y_2(X_1,X_2,k)^{10k}\big]=B(X_1,X_2,k)^{l(X_1,X_2,k)};$
\item[2)]
$4|\big[X_1^{10},X_2^{10k}\big]|_c=4|B(X_1,X_2,k)^{l(X_1,X_2,k)}|>\max\big(|Y_1(X_1,X_2,k)|,|Y_2(X_1,X_2,k)^k|\big)$;
\item[3)]
$Y_1(X_1,X_2,k),\ Y_2(X_1,X_2,k)\not\in\left \langle B(X_1,X_2,k) \right \rangle_\infty;$
\end{enumerate}

\item[V3)]
для $X_1$ и $X_2$ может быть определена неубывающая функция вида $\kappa \big(\underline{\ \ },X_1,X_2 \big)\!\colon\!\mathbb{N}\to\mathbb{N}$ такая, что для любого $N\in\mathbb{N}$ из неравенства $k>\kappa\big(N,X_1,X_2 \big)$ следует, что одновременно выполнено:
\begin{enumerate}
\item[1)]
$l(X_1,X_2,k)\in\{1,2\};$
\item[2)]
$|B(X_1,X_2,k)|>N$;
\item[3)]
$4|B(X_1,X_2,k)|>|s(X_1,X_2,k)|.$
\end{enumerate}
%$$l(X_1,X_2,k)\in\{1,2\} \ \ \text{и}\ \ |B(X_1,X_2,k)^{l(X_1,X_2,k)}|>N\ \text{и}\ 2|B(X_1,X_2,k)^{l(X_1,X_2,k)}|>|s(X_1,X_2,k)|.$$
\end{enumerate}
\end{Lemma}

{\Observation {\rm В дальнейшем, когда это не будет вызывать путаницы, мы будем использовать обозначения $B$, $l$, $Y_1$, $Y_2$ и $s$ вместо $B(X_1,X_2,k)$, $l(X_1,X_2,k)$, $Y_1(X_1,X_2,k)$, $Y_2(X_1,X_2,k)$ и $s(X_1,X_2,k)$ соответственно.
}}

\begin{proof} Т.к. при доказательстве пунктов~$V1)$ и~$V2)$ величины $X_1$, $X_2$ и $k$ можно считать постоянными, то мы будем использовать упрощенные обозначения: $B$, $l$, $Y_1$, $Y_2$ и $s$.

$V1)$ По сути данный пункт утверждает, что коммутатор двух некоммутирующих гиперболических элементов $X_1^{10}$ и $X_2^{10k}$ есть гиперболический элемент, а это немедленно следует из пункта~$C1)$ леммы~\ref{MCL}.

$V2)$ Пусть $X_1=f_1^{-1}\cdot A_1^{k_1}\circ f_1=A_1^{k_1f_1}$ и $X_2=f_2^{-1}\cdot A_2^{k_2}\circ f_2=A_2^{k_2f_2}$, где $A_1,A_2\in H$ --- простые слова, $f_1,f_2\in H$ --- приведенные слова и $k_1,k_2\in\mathbb{N}$. Тогда для некоторого элемента $u\in H$ мы имеем:
$$[X_1^{10},X_2^{10k}]= A_1^{-10k_1f_1}A_2^{-10k_2kf_2}A_1^{10k_1f_1}A_2^{10k_2kf_2}=B^{lu}.$$
Сопрягая это равенство элементом $f_1^{-1}$ и полагая $f_3=f_2f_1^{-1}$, мы получаем:
\begin{equation}\label{ICEQ1}
A_1^{-10k_1}A_2^{-10k_2kf_3}A_1^{10k_1}A_2^{10k_2kf_3}=B^{luf_1^{-1}}.
\end{equation}
Пусть числа $\alpha,\beta\in\mathbb{Z}$ таковы, что приведенное слово $f_4=A_2^{\beta}f_3A_1^{\alpha}$ имеет минимальную длину в двойном смежном классе $\left \langle A_2 \right \rangle f_3\left \langle A_1 \right \rangle.$ Отметим, что выбор $\alpha$ и $\beta$ при определении элемента $f_4$ не зависит от $k$. Сопрягая равенство~(\ref{ICEQ1}) элементом $A_1^{\alpha}$, мы получаем:
\begin{equation}\label{ICEQ2}
P_0Q_0P_1Q_1P_2Q_2P_3Q_3\equiv A_1^{-10k_1}f_4^{-1}A_2^{-10k_2k}f_4A_1^{10k_1}f_4^{-1}A_2^{10k_2k}f_4=B^{luf_1^{-1}A_1^\alpha},
\end{equation}
где $P_0\equiv A_1^{-10k_1}$, $P_1\equiv A_2^{-10k_2k}$, $P_2\equiv A_1^{10k_1}$, $P_3\equiv A_2^{10k_2k}$, $Q_0\equiv Q_2\equiv f_4^{-1}$ и $Q_1\equiv Q_3\equiv f_4$.

Пусть $|f_4|\ge |A_1|+|A_2|$. В таком случае, т.к. $10k_1,10k_2k> 2$ и в силу выбора элемента $f_4$, при приведении слова вида $P_iQ_iP_{i+1}$, $i\in\mathbb{Z}/4\mathbb{Z}$ не может произойти больше чем $1/2(|A_1|+|A_2|)$ сокращений и двух слияний (иначе $f_4$ не будет иметь минимальную длину в $\left \langle A_2 \right \rangle f_3\left \langle A_1 \right \rangle$). Значит при циклическом приведении слова из левой части равенства~(\ref{ICEQ2}) некоторые слоги каждого из слов $P_i$, $Q_i$ останутся незатронутыми. Следовательно, справедлива оценка (напомним, что $|h^k|_c=k|h|_c$ для любого гиперболического элемента $h\in H$ и $k\in\mathbb{N}$):
\begin{align*}
|[X_1^{10},X_2^{10k}]|_c=|B^l|&>2|A_1^{10k_1}|+2|A_2^{10k_2k}|+4|f_4|-4(|A_1|+|A_2|)-8\nonumber\\
&=|A_1^{10k_1}|+|A_2^{10k_2k}|+4|f_4|+(10k_1|A_1|+10k_2k|A_2|-4|A_1|-4|A_2|-8)\nonumber\\
&>|A_1^{10k_1}|+|A_2^{10k_2k}|+2|f_4|.
\end{align*}

Пусть $|f_4|< |A_1|+|A_2|$. В таком случае справедлива оценка:
\begin{align*}
|[X_1^{10},X_2^{10k}]|_c=|B^l|&\ge 2|X_1^{10}|_c+2|X_2^{10k}|_c-2|X_1|_c-2|X_2^k|_c\nonumber\\
&=2|A_1^{10k_1}|+2|A_2^{10k_2k}|-2|A_1^{k_1}|-2|A_2^{k_2k}|\nonumber\\
&=|A_1^{10k_1}|+|A_2^{10k_2k}|+2|A_1|+2|A_2|+(8|A_1^{k_1}|+8k|A_2^{k_2}|-2|A_1|-2|A_2|)\nonumber\\
&>|A_1^{10k_1}|+|A_2^{10k_2k}|+2|f_4|,
\end{align*}
где на первом шаге мы воспользовались пунктом~$C3)$ леммы~\ref{MCL}, а на четвертом шаге тем, что $|f_4|< |A_1|+|A_2|$. Таким образом, вне зависимости от длины слова $f_4$, выполнена оценка:
\begin{equation}\label{ICEQ3}
|B^l|>|A_1^{10k_1}|+|A_2^{10k_2k}|+2|f_4|.
\end{equation}

Понятно, что существует элемент $s_1\in H$, вообще говоря, зависящий от $X_1$, $X_2$ и $k$ такой, что:
\begin{equation}\label{ICEQ5}
B^{l}=[A_1^{10k_1},A_2^{10k_2kf_4}]^{s_1}\ \text{и}\  |s_1|\le \frac{1}{2}\left(|B|+|[A_1^{10k_1},A_2^{10k_2kf_4}]|\right).
\end{equation}
Заметим, что выполнена оценка:
\begin{align}\label{ICEQ4}
|s_1|&\le \frac{1}{2}\left(|B|+|[A_1^{10k_1},A_2^{10k_2kf_4}]|\right)\nonumber\\
&\le \frac{1}{2}\left(|B^l|+2|A_1^{10k_1}|+2|A_2^{10k_2k}|+4|f_4|\right)\nonumber\\
&<\frac{3}{2}|B^l|,
\end{align}
где на третьем шаге мы воспользовались неравенством~(\ref{ICEQ3}).
Положим $s\overset{\textrm{def}}{=}f_1^{-1}A_1^{\alpha}s_1$, тогда:
\begin{align*}
Y_1&=X_1^s=A_1^{k_1f_1s}=A_1^{k_1f_1f_1^{-1}A_1^{\alpha}s_1}=A_1^{k_1s_1},\\
Y_2&=X_2^s=A_2^{k_2f_2s}=A_2^{k_2f_2f_1^{-1}A_1^{\alpha}s_1}=A_2^{k_2f_3A_1^{\alpha}s_1}=A_2^{k_2A_2^{\beta}f_3A_1^{\alpha}s_1}=A_2^{k_2f_4s_1},
\end{align*}
откуда, в соответствии с~(\ref{ICEQ5}), мы имеем $[Y_1^{10},Y_2^{10k}]=B^l$.

Оценим $|Y_1|$ и $|Y_2^k|$:
\begin{align*}
|Y_1|&=|A_1^{k_1s_1}|\le |A_1^{k_1}|+2|s_1|<4|B^l|;\\
|Y_2^k|&=|A_2^{k_2kf_4s_1}|\le |A_2^{k_2k}|+2|f_4|+2|s_1|<4|B^l|,
\end{align*}
где мы воспользовались~(\ref{ICEQ3}) и~(\ref{ICEQ4}).

Чтобы завершить доказательство пункта~$V2)$ осталось показать, что $Y_1,Y_2\not\in\left \langle B \right \rangle_\infty.$ Пусть $Y_1=B^m$ для некоторого $m\in\mathbb{Z}\!\setminus\!\{0\}$ ($m\ne0$, т.к. $Y_1$ есть гиперболический элемент).
%Далее вместо $Y_2^k$ будем писать просто $Y_2$.
В таком случае $[B^{10m},Y_2^{10k}]=B^l$, откуда $B^{10mY_2^{10k}}=B^{l+10m}.$ Т.к. $B$ есть простое слово и центральные длины сопряженных слов равны, то из последнего равенства следует, что $|10m|=|l+10m|$, т.е. или $B^{10mY_2^{10k}} = B^{10m}$ и $l=0$ (и мы сразу получаем противоречие с условием), или же $B^{10mY_2^{10k}} = B^{-10m}$ и $l=-20m$. Во втором случае, учитывая, что в $H$ корень из гиперболического элемента извлекается однозначно, мы имеем равенство $B^{Y_2^{10k}}=B^{-1}$, из которого следует, что $B^{Y_2^{20k}}=B$. Учитывая, что $B$ есть простое слово, $Y_2$ есть гиперболический элемент и опираясь на теорему Куроша, из последнего равенства легко получить, что $Y_2 = B^{m_1}$ для некоторого $m_1\in\mathbb{Z}\!\setminus\!\{0\}$. Но это влечет равенство $[Y_1^{10},Y_2^{10k}]=1$, противоречащие условию. Случай $Y_2\in\left \langle B \right \rangle_\infty$ рассматривается аналогично.

$V3)$. Из пунктов~$C4)$ и~$C5)$ леммы~\ref{MCL} немедленно следует, что соответственно условия $V3.1)$ и $V3.2)$ выполняются при всех достаточно больших $k.$ Покажем, что условие $V3.3)$ тоже выполнено при всех достаточно больших $k$.

Пусть $N_1\in\mathbb{N}$ настолько велико, что для всех $k>N_1$ выполнено условие $V3.1)$ и неравенство $|B(X_1,X_2,k)|>|f_1|+|A_1^{\alpha}|$ (здесь мы пользуемся тем, что величины $|f_1|$, $|A_1|$ и $\alpha$ не зависят от $k$). Тогда справедлива оценка:
\begin{align*}
|s(X_1,X_2,k)|&=|f_1^{-1}A_1^{\alpha}s_1(X_1,X_2,k)|\nonumber\\
&\le |f_1|+|A_1^{\alpha}|+|s_1(X_1,X_2,k)|\nonumber\\
&<|B(X_1,X_2,k)|+\frac{3}{2}|B(X_1,X_2,k)^{l(X_1,X_2,k)}|\\
&\le4|B(X_1,X_2,k)|,
\end{align*}
где на третьем шаге мы воспользовались неравенствами~(\ref{ICEQ4}) и $|B(X_1,X_2,k)|>|f_1|+\alpha|A_1|$, а на четвертом шаге тем, что $l(X_1,X_2,k)\in\{1,2\}.$
%Т.к. $|A_2|\ge 2$, то из неравенства~(\ref{ICEQ3}) следует, что $|B(X_1,X_2,k)^{l(X_1,X_2,k)}|\to\infty$ при $k\to\infty.$ Значит для всех достаточно больших $k$ справедлива %оценка:
%где на третьем шаге мы воспользовались неравенством~(\ref{ICEQ4}) а также тем, что величины $|f_1|$, $|A_1|$ и $\alpha$ не зависят от $k$.

Теперь легко видеть, что нужная нам функция $\kappa\big(\underline{\ \ },X_1,X_2 \big)$ существует. Действительно, при $N=1$ достаточно положить $\kappa\big(1,X_1,X_2 \big)=N_1$ (условия $V3.1)$ и $V3.3)$ выполнены в силу выбора $N_1$, а условие $|B(X_1,X_2,k)|>1$ выполнено для любого $k\in\mathbb{N}$). Если функция $\kappa\big(\underline{\ \ },X_1,X_2 \big)$ уже определена для всех значений, меньших $N$, то число $\kappa\big(N,X_1,X_2 \big)$ достаточно выбрать так, чтобы одновременно выполнялись неравенства $\kappa\big(N,X_1,X_2 \big)>\kappa\big(N-1,X_1,X_2 \big)$ и $|B\big(X_1,X_2,\kappa\big(N,X_1,X_2 \big)\big)|>N$.
\end{proof}

Определим слово $L_2(z_0,z_1)\in F_2(z_0,z_1)'$ следующим образом:
\begin{equation}\label{L2Def}
L_2(z_0,z_1)\overset{\textrm{def}}{=}[z_0^{10},z_1^{10}]^{5000}z_0[z_0^{10},z_1^{10}]^{200}z_1[z_0^{10},z_1^{10}]^{400}z_0^{-1}[z_0^{10},z_1^{10}]^{600}z_1^{-1}[z_0^{10},z_1^{10}]^{5000}.
\end{equation}

\begin{Lemma}\label{LemIvanSL} Пусть группа $H=\zvezda\limits_{\mathfrak{a}\in \mathfrak{A}}\!H_\mathfrak{a}$ есть нетривиальное свободное произведение групп $H_\mathfrak{a}$, $X_1$ и $X_2$ суть некоммутирующие гиперболические элементы группы $H$, а $k$ есть произвольное натуральное число. Если элементы $Y_1$, $Y_2$, $s$ и $B^l$ определены для слов $X_1$, $X_2$ и константы $k$ как в замечании к лемме~\ref{LemIvanPL}, а простое слово $B'$ есть циклическая перестановка слова $B$ при помощи элемента $f$ (т.е. $B' = B^f$), причем $|f|\le 1/2|B|$, то:
\begin{enumerate}
\item[W1)]
слово $L_2(X_1,X_2^k)$ сопряжено элементом $sB^{4990l}f$ циклически приведенному слову:
\begin{equation}\label{WDef}
W\equiv T_1R_1T_2R_2T_3R_3T_4R_4,
\end{equation}
где $T_1=B'^{10l}Y_1^fB'^{10l}$, $T_2=B'^{10l}Y_2^{kf}B'^{10l}$, $T_3=B'^{10l}Y_1^{-f}B'^{10l}$, $T_4=B'^{10l}Y_2^{-kf}B'^{10l}$ --- приведенные слова, $R_i\equiv B'^{(200i-20)l}$, $i=1,2,3$ и $R_4\equiv B'^{9980l}.$ Кроме того, каждое из слов $T_i$ начинается и заканчивается на подслово $B'$ и выполнены оценки:
\begin{align}
4|B^l|<|&T_i|< 25|B^l|\ \ \text{при}\ \ i=1,2,3,4\label{LemIvanSLEst1};\\
11136|B^l|<|&W|<11220|B^l|;\label{LemIvanSLEst2}
\end{align}

\item[W2)]
слово $L_2(X_1,X_2^k)$ является гиперболическим элементом и не является собственной степенью в $H$.
%$не является собственной степенью в $H$ и $L_2(X_1,X_2^k)\in M_H$;

%\item[W3)]
%для любого $C\in\mathbb{N}$ существует такое $N\in\mathbb{N}$, что $|L_2(X_1,X_2)|_r>C$ при $|X_2|_c\ge N.$
\end{enumerate}
\end{Lemma}

{\Observation {\rm В дальнейшем (при доказательстве леммы~\ref{LemIvanML1}) нам придется <<подправить>> простое слово $B$, перейдя к сопряженному ему простому слову $B'$. Элементы $B'$ и $f$ (<<поправка>> для $B$) введены в условие леммы~\ref{LemIvanSL} для того, чтобы в пункте~$W1)$ получить оценки, не зависящие от небольшой ($|f|\le 1/2|B|$) поправки $f$.}}

\begin{proof} $W1)$ Учитывая пункты~$V1)$ и~$V2.1)$ леммы~\ref{LemIvanPL} и определение~(\ref{L2Def}) несложно видеть, что элемент $L_2(X_1,X_2^k)$ сопряжен элементом $sB^{4990l}f$ (не обязательно приведенному) слову:
\begin{equation}\label{W0Def}
W_0\equiv\underbrace{B'^{10l}Y_1^fB'^{10l}}_{=T_1}B'^{180l}\underbrace{B'^{10l}Y_2^{kf}B'^{10l}}_{=T_2}B'^{380l}\underbrace{B'^{10l}Y_1^{-f}B'^{10l}}_{=T_3}B'^{580l}\underbrace{B'^{10l}Y_2^{-kf}B'^{10l}}_{=T_4}B'^{9980l}.
\end{equation}

Положим $\alpha_1=1$ и $\alpha_2=k$, тогда из пункта~$V2.2)$ леммы~\ref{LemIvanPL} следует, что справедлива оценка $|Y_i^{\alpha_i}|<4|B^l|$, $i=1,2$. Значит, $|Y_i^{\alpha_if}|\le |Y_i^{\alpha_i}|+2|f|<5|B^l|=5|B'^l|$ и мы имеем $8l\ge 5l+3>|Y_i^{\alpha_if}|/|B'|+3$. Учитывая последнюю оценку, а также то, что $Y_i^{\alpha_if}$ есть гиперболический элемент, в силу пункта~1) следствия~\ref{LemBFLL}, мы заключаем, что при приведении слова вида $B'^{8l}Y_i^{\pm\alpha_if}B'^{8l}$ его первый и последний слог не изменяется.
Учитывая, что $B'$ есть простое слово, из сказанного легко следует оценка $|T_i|>4|B^l|$, а также то, что слова $T_i$ начинаются и заканчиваются на подслово $B'$ (в частности, слова $T_i$ циклически приведены). Учитывая вид слова $W_0$ и то, что $B'$ есть простое слово, теперь легко видеть, что слово $W$ циклически приведено. Оценка $|T_i|< 25|B^l|$ немедленно следует из оценки $|Y_i^{\alpha_if}|<5|B^l|$, а оценки~(\ref{LemIvanSLEst2}), учитывая что слово $W$ приведено, легко получаются из оценок~(\ref{LemIvanSLEst1}).

$W2)$ Положим $f\equiv1$. Т.к. слово $L_2(X_1,X_2^k)$ сопряжено в $H$ циклически приведенному слову $W$, то из оценки~(\ref{LemIvanSLEst2}) следует, что $L_2(X_1,X_2^k)$ есть гиперболический элемент. Осталось показать, что $W$ не является собственной степенью в $H$.

Положим $a\equiv T_1R_1T_2R_2T_3R_3T_4$, тогда $W\equiv aB^{9980l}$. Т.к. $R_i\equiv B^{(200i-20)l}$, $i=1,2,3$, то из оценки $|T_i|< 25|B^l|$ легко следует, что:
$$|a|<4\cdot 25|B^l|+(180+380+580)|B^l|=1240l|B|.$$
Т.к. слово $W$ циклически приведено, то $W\equiv aB^{9980l}\equiv C^m$ для некоторого простого слова $C$ и $m\in\mathbb{N}$. Если $m\ge2$, то, т.к. $|a|<1240l|B|$ и $|W|>11136|B^l|$, слово $W$ заканчивается на одновременно $B$- и $C$-периодическое подслово, длина которого больше, чем $|B|+|C|.$ Следовательно, в силу пункта~$P4)$ леммы~\ref{MWL}, мы имеем $C\equiv B$, т.е. $T_1R_1T_2R_2T_3R_3T_4\equiv B^{m-9980l}$ и $T_1$ есть $B$-периодическое слово. Т.к. слово $T_1=B^{10l}Y_1B^{10l}$ начинается и заканчивается на подслово $B$, то из пункта~$P2)$ леммы~\ref{MWL} следует, что $Y_1\in \left \langle B \right \rangle_\infty$, а это противоречит пункту~$V2.3)$ леммы~\ref{LemIvanPL}.
\end{proof}

\begin{Lemma}\label{LemICAdd2} Пусть группа $H=\zvezda\limits_{\mathfrak{a}\in \mathfrak{A}}\!H_\mathfrak{a}$ есть нетривиальное свободное произведение групп $H_\mathfrak{a}$, $\underline{X}=(X_0,X_1,X_2)$ есть кортеж гиперболических элементов группы $H$ и выполнены условия:
\begin{enumerate}
    \item[\sffamily $\mathrm{1)}$] $[X_i,X_{i+1}]\ne1$ при $i=0,1$;

    \item[\sffamily $\mathrm{2)}$] $m>\kappa \big(|L_2(X_0,X_1)|,X_1,X_2)|\big)$, где функция $\kappa \big(\underline{\ \ },X_{1},X_{2} \big)\!\colon\!\mathbb{N}\to\mathbb{N}$ определена для слов $X_{1}$ и $X_{2}$ как в лемме~\ref{LemIvanPL}.
  \end{enumerate}
Тогда для любых целых чисел $l_1$ и $l_2$, одновременно не равных нулю, из равенства:
        \begin{equation}\label{LemICAdd2MEQ}
        \left ( L_2(X_0,X_1)^{l_1}L_2(X_1,X_2^m)^{l_2} \right )^{l_3}=X_1^{l_4}
        \end{equation}
следует, что $l_3=0$.
\end{Lemma}
\begin{proof} %Т.к. $X_1$ есть гиперболический элемент, то из~(\ref{LemICAdd2MEQ}) ясно, что равенство $l_3=0$ влечет равенство $l_4=0$, т.е. достаточно доказать, что $l_3=0$.
 Рассмотрим несколько случаев.

1) Пусть $l_2=0$ (в таком случае $l_1\ne 0$). Если $l_3\ne 0$, то из равенства~(\ref{LemICAdd2MEQ}) следует, что:
$$|L_2(X_0,X_1)|_r=|L_2(X_0,X_1)^{l_1l_3}|_r=|X_1^{l_4}|_r=|X_1|_r.$$
Но это невозможно, т.к.:
\begin{equation}\label{LemICAdd2MEQ4}
|X_1|_r<4|[X_0^{10},X_1^{10}]|_c<\frac{4}{11136}|L_2(X_0,X_1)|_c<|L_2(X_0,X_1)|_r,
\end{equation}
где первое, второе и третье неравенство выполнены в силу пункта~$V2.2)$ леммы~\ref{LemIvanPL}, оценки~(\ref{LemIvanSLEst2}) и пункта~$W2)$ леммы~\ref{LemIvanSL}, соответственно.

2) Случай $l_1=0$ (и $l_2\ne 0$) аналогичен случаю 1).

3) Пусть $l_1\ne 0$ и $l_2>0$. Предположим, что $l_3\ne 0$ (без ограничения общности мы будем считать, что $l_3>0$). В силу пункта~$W2)$ леммы~\ref{LemIvanSL} мы имеем $L_2(X_0,X_1)=g^{-1}\cdot A\circ g$ для некоторого простого слова $A$ и приведенного слова $g$. Полагая $f\equiv1$ в условии леммы~\ref{LemIvanSL} и учитывая ее пункт~$W1)$, несложно понять, что после сопряжения равенства~(\ref{LemICAdd2MEQ}) элементом $sB^{-38l}$ (здесь мы используем обозначения $s$ и $B^l$ вместо более аккуратных $s(X_1,X_2,m)$ и $B(X_1,X_2,m)^{l(X_1,X_2,m)}$), его левая часть будет иметь вид:
\begin{equation}\label{LemICAdd2MEQ2}
\Big( \underbrace{B^{38l}f^{-1}A^{l_1}fB^{38l}}_{=U}B^{-38l}\underbrace{B^{4990l}\overbrace{T_1B^{180l}T_2B^{380l}T_3B^{580l}T_4B^{4990l}B^{4990l}}^{\equiv V}T_1\cdots T_4B^{4990l}}_{l_2\ \text{раз}\ B^{4990l}VB^{-4990l}}B^{-38l} \Big)^{l_3},
\end{equation}
где $f=gs$ есть приведенное слово, элементы $T_i$, $i=1,2,3,4$ определены для слов $X_1$, $X_2$ и константы $k$ как в пункте~$W1)$ леммы~\ref{LemIvanSL}, а $V$ --- циклически приведенное слово (т.к., в силу пункта~$W1)$ леммы~\ref{LemIvanSL}, слова $T_i$ начинаются и заканчиваются на подслово $B$).

Покажем, что при приведении слова $B^{38l}f^{-1}A^{l_1}fB^{38l}$ его первый и последний слог не затрагивается. Для этого проверим, что выполнены все условия леммы~\ref{LemIvanAdd}.
Из условия 2) данной леммы, в силу пунктов~$V3.2)$ и~$V3.3)$ леммы~\ref{LemIvanPL}, мы имеем: $|B|>|L_2(X_0,X_1)|$ и $4|B|>|s|$. Т.к. $|L_2(X_0,X_1)|\ge |A|+2|g|-1$, то легко видеть, что выполнены неравенства: $|B|>|A|$ и $|B|>g$. Значит $|f|\le |g|+|s|<|B|+4|B|=5|B|$ и полагая (в условии леммы~\ref{LemIvanAdd}) $k=5$, мы имеем $38l\ge 38=6k+8$. Следовательно, все условия леммы~\ref{LemIvanAdd} выполнены и в приведенном слове $U=B^{38l}f^{-1}A^{l_1}fB^{38l}$ первый и последний слог совпадает с первым и последним слогом слова $B$, соответственно (в частности, из этого следует, что слово $U$ циклически приведено).
Значит после приведения слова~(\ref{LemICAdd2MEQ2}) мы будем иметь:
\begin{equation}\label{LemICAdd2MEQ3}
E\equiv\Big(U\cdot B^{4952l}\cdot V^{l_2-1}\cdot T_1\cdot B^{180l}\cdot T_2\cdot B^{380l}\cdot T_3\cdot B^{580l} \cdot T_4\cdot B^{4952l} \Big)^{l_3}.
\end{equation}
Заметим, что слово $E$ циклически приведено, т.к. его первый и последний слог совпадает с первым и последним слогом слова $B$, соответственно. Значит из равенства $E=X_1^{l_4sB^{-38l}}$ следует, что $E\equiv C^{kl_4}$ для некоторого $k\in\mathbb{N}$ и простого слова $C$, причем $|X_1|_r=|C|.$
%$X_1^{l_4sB^{-38l}}=C^{kl_4}$ для некоторого простого слова $C$ и $k\in\mathbb{N}$, т.е. $W\equiv C^{kl_4}$ и $|X_1|_r=|C|$.
Учитывая~(\ref{LemICAdd2MEQ4}) мы имеем:
\begin{equation}\label{LemICAdd2MEQ5}
|C|=|X_1|_r<|L_2(X_0,X_1)|_r=|A|<|B|.
\end{equation}
Т.к. $|C|<|B|$, то легко видеть, что слово $E$ оканчивается на одновременно $C$- и $B$-периодическое подслово $U$, длина которого больше, чем $|B|+|C|$ (например, в качестве $U$ можно взять $B^{4952l}$). Значит, в силу пункта~$P4)$ леммы~\ref{MWL}, мы имеем $B\equiv C$, следовательно, $|C|=|B|$, что противоречит~(\ref{LemICAdd2MEQ5}).

4) Случай $l_1\ne 0$ и $l_2<0$ аналогичен случаю 3).
\end{proof}

\begin{Lemma}\label{LemIvanML1} Пусть группа $H=\zvezda\limits_{\mathfrak{a}\in \mathfrak{A}}\!H_\mathfrak{a}$ есть нетривиальное свободное произведение групп $H_\mathfrak{a}$, $X_1, X_2,\hat{X}_1,\hat{X}_2$ суть гиперболические элементы группы $H$, причем $[X_1,X_2]\ne 1$, а константа $k\in\mathbb{N}$ выбрана так, что $k\ge\kappa \big(1,X_1,X_2 \big),$ где функция $\kappa \big(\underline{\ \ },X_1,X_2 \big)$ определена для слов $X_1$ и $X_2$ как в пункте~$V3)$ леммы~\ref{LemIvanPL}. Тогда из равенства:
$$L_2(X_1,X_2^k)=L_2(\hat{X}_1,\hat{X}_2^k)$$
следует, что существует такой элемент $v\in H$, что $X_1=\hat{X}_1^v$, и $X_2=\hat{X}_2^v$.
\end{Lemma}
\begin{proof} Пусть элементы $s$, $Y_1$, $Y_2$, $B^l$, $f$, $B'$, $R_i$ и $T_i$, $i=1,2,3,4$ (соответственно $\hat{s}$, $\hat{Y}_1$, $\hat{Y}_2$, $\hat{B}^{\hat{l}}$ $\hat{f}$, $\hat{B}'$, $\hat{R}_i$ и $\hat{T}_i$, $i=1,2,3,4$) определены для $X_1$, $X_2$ и $k$ (соответственно для $\hat{X}_1$, $\hat{X}_2$ и $k$) как в условии леммы~\ref{LemIvanSL}. Положим $f\equiv1$, а элемент $\hat{f}$ определим ниже. Опираясь на пункт~$W1)$ леммы~\ref{LemIvanSL}, мы заключаем, что слова $L_2(X_1,X_2^k)$ и $L_2(\hat{X}_1,\hat{X}_2^k)$ сопряжены (при помощи элементов $sB^{4990l}$ и $\hat{s}\hat{B}^{4990\hat{l}}\hat{f}$) циклически приведенным словам $W$ и $\hat{W}$ вида~(\ref{WDef}):
$$W\equiv T_1R_1T_2R_2T_3R_3T_4R_4\ \ \text{и}\ \ \hat{W}\equiv \hat{T}_1\hat{R}_1\hat{T}_2\hat{R}_2\hat{T}_3\hat{R}_3\hat{T}_4\hat{R}_4,$$
соответственно. Т.к. $L_2(X_1,X_2^k)=L_2(\hat{X}_1,\hat{X}_2^k)$, то слово $W$ есть циклическая перестановка слова $\hat{W}$, следовательно, слова $W$ и $\hat{W}$ равны как циклические слова и $|W|=|\hat{W}|$.

Пусть $|B^l|\ge |\hat{B}^{\hat{l}}|$ (случай $|B^l|\le |\hat{B}^{\hat{l}}|$ аналогичен и мы его опускаем). В таком случае, в силу оценок~(\ref{LemIvanSLEst2}) и равенства $|W|=|\hat{W}|$, мы имеем:
\begin{equation}\label{LemIvanML1Est}
|B^l|<\frac{1}{11136}|W|=\frac{1}{11136}|\hat{W}|<\frac{11220}{11136}|\hat{B}^{\hat{l}}|<1,\!01|\hat{B}^{\hat{l}}|.
\end{equation}

Т.к. $|R_4|=9980|B^l|\ge9980|\hat{B}^{\hat{l}}|$, $|\hat{R}_4|=9980|\hat{B}^{\hat{l}}|$ и $|\hat{W}|<11220 |\hat{B}^{\hat{l}}|$, то легко видеть, что подслова $R_4$ и $\hat{R}_4$ циклического слова $\hat{W}$ имеют общее подслово $Q$ длины большей, чем $\left ( 9980-\frac{11220}{2} \right )|\hat{B}^{\hat{l}}|=4370|\hat{B}^{\hat{l}}|.$ В силу~(\ref{LemIvanML1Est}) мы имеем:
$$|Q|>4370|\hat{B}^{\hat{l}}|>\frac{4370}{1,\!01}|B^l|>4326|B^l|,$$
значит $|Q|>2163|B^l|+2185|\hat{B}^{\hat{l}}|$. Т.к. $Q$ есть одновременно $B$- и $\hat{B}'$-периодическое слово и $|Q|>|B|+|\hat{B}|$, то из пункта~$P4)$ леммы~\ref{MWL} следует, что $\hat{B}'$ есть циклическая перестановка $B$, в частности $|B|=|\hat{B}'|.$ Значит мы можем считать, что элемент $\hat{f}$ таков, что $\hat{B}^{\hat{f}}=\hat{B}'\equiv B$ и $|\hat{f}|\le 1/2|\hat{B}|$.  Из равенства $|B|=|\hat{B}|$ и неравенств $|B^l|\ge |\hat{B}^{\hat{l}}|$ и $|B^l|<1,\!01|\hat{B}^{\hat{l}}|$ следует, что $\hat{l}\le l<1,\!01\hat{l}$. Т.к. $k\ge\kappa \big(1,X_1,X_2 \big),$ то в силу пункта~$V3.1)$ леммы~\ref{LemIvanPL}, мы имеем $l\in\{1,2\}$. Значит из неравенств $\hat{l}\le l<1,\!01\hat{l}$ следует, что $l=\hat{l}$. Т.о. слово $R_i$ равно слову $\hat{R}_i$ (как элемент группы $H$), однако эти слова могут различаться как подслова циклического слова $\hat{W}$.

Мы будем говорить, что подслова $R_i$ и $\hat{R}_j$ циклического слова $\hat{W}$ перекрываются, если они имеют общее нетривиальное подслово. Покажем, что слово $R_i$ может перекрываться не более чем с одним словом $\hat{R}_j$. Действительно, если это не так, то несложно понять, что хотя бы одно из слов $\hat{T}_k$, $k=1,2,3,4$ будет $B$-периодическим. В силу пункта~$W1)$ леммы~\ref{LemIvanSL}, слово $\hat{T}_k$ начинается и заканчивается на подслово $B$, из чего, используя пункт~$P2)$ леммы~\ref{MWL} и равенство $B=\hat{B}^{\hat{f}}$, легко вывести, что хотя бы одно из слов $\hat{Y}_1$, $\hat{Y}_2$ принадлежит $\langle \hat{B} \rangle_\infty$, а это противоречит пункту~$V2.3)$ леммы~\ref{LemIvanPL}. С другой стороны, учитывая, что слова $T_k$ и $\hat{T}_k$ не могут быть $B$-периодическими и принимая во внимание оценки $|R_i|,|\hat{R}_j|\ge 180|B^l|$ и $|T_i|, |\hat{T}_i|<25|B^l|$, несложно проверить, что каждое $R_i$ перекрывается хотя бы с одним $\hat{R}_j$ по подслову длины большей, чем $(180-2\cdot25)|B^l|=130|B^l|.$ Таким образом мы показали, что в циклическом слове $\hat{W}$ каждое подслово $R_i$ перекрывается ровно с одним подсловом $\hat{R}_j$, причем длина подслова, по которому перекрываются $R_i$ и $\hat{R}_j$, больше чем $130|B^l|$. Как мы уже видели, $R_4$ перекрывается с $\hat{R}_4$, значит, учитывая соображения ориентации, несложно понять, что подслово $R_3$ перекрывается лишь с подсловом $\hat{R}_3$, $R_2$ перекрывается лишь с $\hat{R}_2$, а $R_1$ с $\hat{R}_1.$

Итак, мы показали, что в циклическом слове $\hat{W}$ подслова $R_i$  и $\hat{R}_i$, $i=1,2,3,4$ имеют общее подслово, длина которого больше, чем $130|B^l|$, а подслова $R_i$ и $\hat{R}_j$ при $i\ne j$ не перекрываются. Из этого следует, что существуют такие слова $Q_i$ и $\hat{Q}_i$ такие, что:
\begin{itemize}
\item[$\bullet$]
$Q_i\equiv \hat{Q}_i\equiv B^l$, $i=1,2,3,4$;

\item[$\bullet$]
слова $Q_i$ и $\hat{Q}_i$ являются подсловами слов $R_i$ и $\hat{R}_i$ соответственно, т.е. для некоторых $m_i,\hat{m}_i\in\mathbb{N}\cup \{0\}$ выполнено: $R_i \equiv B^{m_il}\cdot Q_i\cdot B^{(200i-21-m_i)l}$, $\hat{R}_i \equiv B^{\hat{m}_il}\cdot \hat{Q}_i\cdot B^{(200i-21-\hat{m}_i)l}$ при $i=1,2,3$ и $R_4 \equiv B^{m_4l}\cdot Q_4\cdot B^{(9979-m_4)l}$, $\hat{R}_4 \equiv B^{\hat{m}_4l}\cdot \hat{Q}_4\cdot B^{(9979-\hat{m}_4)l}$;

\nobreak
\item[$\bullet$]
в циклическом слове $\hat{W}$ слова $Q_i$ и $\hat{Q}_i$ совпадают (т.е. $Q_i$ есть подслово с слове, по которому перекрываются $R_i$ и $\hat{R}_i$ в циклическом слове $\hat{W}$).
\end{itemize}
Введем следующие обозначения: $S_{i1}=B^{m_il}$, $\hat{S}_{i1}=B^{\hat{m}_il}$ при $i=1,2,3,4$, $S_{i2}=B^{(200i-21-m_i)l}$, $\hat{S}_{i2}=B^{(200i-21-\hat{m}_i)l}$ при $i=1,2,3$ и $S_{42}=B^{(9979-m_4)l}$, $\hat{S}_{42}=B^{(9979-\hat{m}_4)l}$. Рассмотрим циклические слова:
$${ \underbrace{S_{42}\!\cdot\! B^{10l}Y_1B^{10l}\!\cdot\! S_{11}}_{=V_1}\!\cdot Q_1\!\cdot\! \underbrace{S_{12}\!\cdot\! B^{10l} Y_2^kB^{10l}\!\cdot\! S_{21}}_{=V_2}\!\cdot Q_2\!\cdot\! \underbrace{S_{22}\!\cdot\! B^{10l}Y_1^{-1}B^{10l}\!\cdot\! S_{31}}_{=V_3}\!\cdot Q_3\!\cdot\! \underbrace{S_{32}\!\cdot\! B^{10l}Y_2^{-k}B^{10l}\!\cdot\! S_{41}}_{=V_4}\!\cdot Q_4,}$$
и
$${ \underbrace{\hat{S}_{42}\!\cdot\! B^{10l}\hat{Y}_1^{\hat{f}}B^{10l}\!\cdot\! \hat{S}_{11}}_{=\hat{V}_1}\!\cdot \hat{Q}_1\!\cdot\! \underbrace{\hat{S}_{12}\!\cdot\! B^{10l} \hat{Y}_2^{k\hat{f}}B^{10l}\!\cdot\! \hat{S}_{21}}_{=\hat{V}_2}\!\cdot \hat{Q}_2\!\cdot\! \underbrace{\hat{S}_{22}\!\cdot\! B^{10l}\hat{Y}_1^{-\hat{f}}B^{10l}\!\cdot\! \hat{S}_{31}}_{=\hat{V}_3}\!\cdot \hat{Q}_3\!\cdot\! \underbrace{\hat{S}_{32}\!\cdot\! B^{10l}\hat{Y}_2^{-k\hat{f}}B^{10l}\!\cdot\! \hat{S}_{41}}_{=\hat{V}_4}\!\cdot \hat{Q}_4.}$$
При приведении этих слов сокращения или слияния могут происходить лишь в подсловах вида $B^{10l}Y_1^{\pm1}B^{10l}$, $B^{10l}\hat{Y}_1^{\pm\hat{f}}B^{10l}$, $B^{10l}Y_2^{\pm k}B^{10l}$ и $B^{10l}\hat{Y}_2^{\pm k\hat{f}}B^{10l}$, причем (в силу пункта~$W1)$ леммы~\ref{LemIvanSL}) первые и последние слоги этих слов не изменяются при приведении. Следовательно, подслова $Q_i$ и $\hat{Q}_i$ не затрагиваются при приведении. Учитывая~(\ref{W0Def}) нетрудно видеть, что после приведения первое из этих циклических слов становится посложно равно циклическому слову $W$, а второе становится посложно равно циклическому слову $\hat{W}$. Т.к. циклические слова $W$ и $\hat{W}$ равны, а $Q_i$ совпадает с $\hat{Q}_i$ в циклическом слове $\hat{W}$, то из сказанного следует, что $\hat{V}_i=V_i$ при $i=1,2,3,4$.

Из равенства $\hat{V}_1=V_1$ мы имеем:
\begin{equation}\label{EqY1}
\hat{Y}_1^{\hat{f}}=( \hat{S}_{42}B^{10l}  )^{-1}( S_{42}B^{10l} )Y_1(B^{10l}S_{11})(B^{10l}\hat{S}_{11})^{-1}=B^{(-m_4+\hat{m}_4)l}Y_1B^{(m_1-\hat{m}_1)l}.
\end{equation}
Аналогичным образом из равенств $\hat{V}_i=V_i$, $i=2,3,4$ мы имеем:
\begin{align}
\hat{Y}_2^{k\hat{f}}&=B^{(-m_1+\hat{m}_1)l}Y_2^{k}B^{(m_2-\hat{m}_2)l};\label{EqY2}\\
\hat{Y}_1^{-\hat{f}}&=B^{(-m_2+\hat{m}_2)l}Y_1^{-1}B^{(m_3-\hat{m}_3)l};\label{EqY1-1}\\
\hat{Y}_2^{-k\hat{f}}&=B^{(-m_3+\hat{m}_3)l}Y_2^{-k}B^{(m_4-\hat{m}_4)l}\label{EqY2-1}.
\end{align}
Вводя обозначения $\alpha_i=(m_i-\hat{m}_i)l$, $i=1,2,3,4$, из равенств~(\ref{EqY1}) и~(\ref{EqY1-1}) мы получаем:
\begin{equation}\label{EqBY1}
B^{\alpha_3-\alpha_4}Y_1=Y_1B^{-\alpha_1+\alpha_2},
\end{equation}
тогда как из равенств~(\ref{EqY2}) и~(\ref{EqY2-1}), учитывая однозначность извлечения корня из гиперболического элемента в $H$, мы получаем:
\begin{equation}\label{EqBY2}
B^{-\alpha_1+\alpha_4}Y_2=Y_2B^{-\alpha_2+\alpha_3}.
\end{equation}
В силу простоты слова $B$ и того, что $Y_1,Y_2\not\in\left \langle B \right \rangle_\infty$ (см. пункт~$V2.3)$ леммы~\ref{LemIvanPL}), из равенств~(\ref{EqBY1}) и~(\ref{EqBY2}) легко следует, что $\alpha_3-\alpha_4=-\alpha_1+\alpha_2=0$ и $-\alpha_1+\alpha_4=-\alpha_2+\alpha_3=0$, соответственно. Откуда мы получаем, что $\alpha_1=\alpha_2=\alpha_3=\alpha_4$. Значит, с учетом сказанного, из равенств~(\ref{EqY1}) и~(\ref{EqY2}) следует, что:
\begin{align*}
\hat{Y}_1^{\hat{f}}& = B^{-\alpha_1}Y_1B^{\alpha_1}=Y_1^{B^{\alpha_1}};\\
\hat{Y}_2^{\hat{f}}&=B^{-\alpha_1}Y_2B^{\alpha_1}=Y_2^{B^{\alpha_1}}.
\end{align*}
Осталось заметить, что:
$$X_i = Y_i^{s^{-1}}=\hat{Y}_i^{\hat{f}B^{-\alpha_1}s^{-1}}=\hat{X}_i^{\hat{s}\hat{f}B^{-\alpha_1}s^{-1}}\ \ \text{при}\ i=1,2,$$ т.е. можно положить $v\overset{\textrm{def}}{=}\hat{s}\hat{f}B^{-\alpha_1}s^{-1}.$
\end{proof}

\label{pagewith1} Введем некоторые обозначения. \emph{Подкортежем}, $\underline{a}[j,i]$, кортежа $\underline{a}=(a_0,a_1,\dots, a_{n-2},a_{n-1})$, $n\ge1$ мы будем называть кортеж $\underline{a}[j,i]=(a_i,a_{i+1}\dots,a_{i+j-2},a_{i+j-1})$, где $j=1,\ \dots,\ n$ и $i=0,\ \dots,\ n-j$. Например, если $\underline{a}=(1,1,2,0,1,23)$, то $\underline{a}[3,1]=(1,2,0)$, $\underline{a}[1,5]=(23)$ и $\underline{a}[6,0]=\underline{a}$. Если мы смотрим на $\underline{a}[j,i]$ и $\underline{a}[j',i']$ просто как на кортежи, то мы считаем их равными, тогда и только тогда, когда они совпадают покомпонентно (в частности, равные кортежи имеют одно и тоже число компонент, т.е. $j=j'$). Если же мы смотрим на $\underline{a}[j,i]$ и $\underline{a}[j',i']$ как на подкортежи кортежа $\underline{a}$, то мы считаем их равными тогда и только тогда, когда $j=j'$ и $i=i'$. Например, если вернуться к рассмотренному выше примеру, то $\underline{a}[1,1]=(1)=\underline{a}[1,4]$ как самостоятельные кортежи и $\underline{a}[1,1]\ne\underline{a}[1,4]$ как подкортежи кортежа $\underline{a}$.

\label{pagewith2} Пусть кортеж $\underline{x}=(x_0,\dots,x_{n-1})$, $n\ge2$ гиперболических элементов $x_i$ группы $H=\zvezda\limits_{\mathfrak{a}\in \mathfrak{A}}\!H_\mathfrak{a}$ таков, что $[x_i,x_{i+1}]\ne1$ при $i=0,\ \dots,\ n-2$. Для всех подкортежей $\underline{x}[j,i]$, $j=1,\ \dots,\ n$, $i=0,\ \dots,\ n-j$ кортежа $\underline{x}$ определим слова $T_{\underline{x}[i,j]}(z_i,\dots,z_{i+j-1})\in F_n(z_0,\dots,z_{n-1})'$ следующим образом:

1) на множестве упорядоченных пар $(j,i)$, где $j=1,\ \dots,\ n$ и $i = 0,\ \dots,\ n-j$ введем линейный порядок, полагая, что $(j',i')<(j,i)$ если или $j'<j$ или $j'=j$ и $i'<i$;

%2) положим $m_{j,-1}=1$ при $j=2,\ \dots,\ n$;

2) для пар $(1,i)$, $i = 0,\ \dots,\ n-1$ положим:
\begin{equation}\label{AddIII_NML_EQ0}
    T_{\underline{x}[1,i]}(z_i)\overset{\textrm{def}}{=}z_i;
\end{equation}
3) если для всех пар $(j',i')$, меньших пары $(j,i)$, $j\ge2$ слова $T_{\underline{x}[j',i']}(z_{i'},\dots,z_{i'+j'-1})$ уже определены, то положим:
        \begin{equation}\label{AddIII_NML_EQ1}
            T_{\underline{x}[j,i]}(z_i,\dots,z_{i+j-1})\overset{\textrm{def}}{=}L_2\Big(T_{\underline{x}[j-1,i]}\big(z_i,\dots,z_{i+j-2}\big)^{m_{j,i-1}},T_{\underline{x}[j-1,i+1]}\big(z_{i+1},\dots,z_{i+j-1}\big)^{m_{j,i}}\Big),
        \end{equation}
    где слово $L_2(z_0,z_1)$ определено как~(\ref{L2Def}), константа $m_{j,i-1}\in\mathbb{N}$ определена на предыдущем шаге, если $i>0$ и $m_{j,-1}=1$, если $i=0$, а константа $m_{j,i}\in\mathbb{N}$ выбирается таким образом, чтобы было выполнено условие:
        \begin{enumerate}
          \item[\sffamily $\mathrm{\bullet}$]
        если $i=0$, то:
        \begin{equation}\label{AddIII_NML_EQ2T21}
        %m_{j,i}>\kappa \big(1,X_{j-1,i},X_{j-1,i+1}\big),
        m_{j,0}>\kappa \big(1,X_{j-1,0},X_{j-1,1}\big),
         \end{equation}

        где $X_{j-1,0}=T_{\underline{x}[j-1,0]}\big(x_0,\dots,x_{j-2}\big)$, $X_{j-1,1}=T_{\underline{x}[j-1,1]}\big(x_{1},\dots,x_{j-1}\big)$ и функция $\kappa \big(\underline{\ \ },X_{j-1,0},X_{j-1,1} \big)$ определена для слов $X_{j-1,0}$ и $X_{j-1,1}$ как в пункте~$V3)$ леммы~\ref{LemIvanPL};

     \item[\sffamily $\mathrm{\bullet}$]
            если $i>0$, то:
              \begin{equation}\label{AddIII_NML_EQ2T22}
            m_{j,i}>\kappa \big(|X_{j,i-1}|,X_{j-1,i}^{m_{j,i-1}},X_{j-1,i+1}\big),
             \end{equation}
    где $|X_{j,i-1}|$ означает (как обычно) длину слова $X_{j,i-1}=T_{\underline{x}[j,i-1]}\big(x_{i-1},\dots,x_{i+j-2}\big)$, $X_{j-1,i}=T_{\underline{x}[j-1,i]}\big(x_{i},\dots,x_{i+j-2}\big)$ и $X_{j-1,i+1}=T_{\underline{x}[j-1,i+1]}\big(x_{i+1},\dots,x_{i+j-1}\big).$
    \end{enumerate}

Заметим, что при $j=1$ приведенное построение слов $T_{\underline{x}[1,i]}(z_i)$, очевидно, корректно. Тогда как при $j\ge2$ для того, чтобы данное построение было корректным необходимо и достаточно, чтобы для каждой пары $(j,i)$, $j\ge2$ была определена функция $\kappa \big(\underline{\ \ },X_{j-1,i}^{m_{j,i-1}},X_{j-1,i+1} \big)$.
%Из пункта~$V3)$ леммы~\ref{LemIvanPL} следует, что функция $\kappa \big(\underline{\ \ },X_{j-1,i},X_{j-1,i+1} \big)$ определена, если $X_{j-1,i}$ и $X_{j-1,i+1}$ суть некоммутирующие гиперболические элементы.
Корректность построения слов $T_{\underline{x}[j,i]}(z_i,\dots,z_{i+j-1})$ при $j\ge2$ будет установлена в ходе доказательства следующей леммы.

\begin{Lemma}\label{AddIIINewLemIGProp} Пусть группа $H=\zvezda\limits_{\mathfrak{a}\in \mathfrak{A}}\!H_\mathfrak{a}$ есть нетривиальное свободное произведение групп $H_\mathfrak{a}$, $M_H$ есть множество всех гиперболических элементов группы $H$ и кортеж $\underline{x}=(x_0,\dots,x_{n-1})$, $n\ge2$ гиперболических элементов $x_i$ группы $H$ таков, что $[x_i,x_{i+1}]\ne1$ при $i=0,\ \dots,\ n-2$. Тогда для всех слов $T_{\underline{x}[j,i]}(z_i,\dots,z_{i+j-1})$, где $j=2,\ \dots,\ n$ и $i=0,\ \dots,\ n-j$, определенных для кортежа $\underline{x}$ как~(\ref{AddIII_NML_EQ1}), выполнено:
\begin{enumerate}
\item[E1)] $X_{j,i}\in M_H$ и $[X_{j-1,i},X_{j-1,i+1}]\ne1,$ где $X_{j,i}=T_{\underline{x}[j,i]}\big(x_{i},\dots,x_{i+j-1}\big)$ для всех пар $(j,i)$;

\item[E2)] если $(y_0,\dots,y_{n-1})$ есть такой кортеж элементов из $H$, что $[y_k,y_{k+1}]=1$ для некоторого $k=i,\ \dots,\ i+j-2$, то $T_{\underline{x}[j,i]}(y_i,\dots,y_{i+j-1})=1$;

\item[E3)] если $(y_0,\dots,y_{n-1})$ есть кортеж элементов из $M_H$, то слово $T_{\underline{x}[j,i]}(y_i,\dots,y_{i+j-1})$ или равно единице или является гиперболическим элементом и не является собственной степенью в $H.$
%не является собственной степенью в $H$ и содержится во множестве $M_H.$
\end{enumerate}
\end{Lemma}
\begin{proof} Индукцией по множеству упорядоченных пар $(j,i)$ мы установим все пункты леммы одновременно.

При $(j,i)=(2,0)$.

В силу~(\ref{AddIII_NML_EQ0}) и~(\ref{AddIII_NML_EQ1}) мы имеем: $$T_{\underline{x}[2,0]}(z_0,z_1)=L_2(z_0,z_1^{m_{2,0}}),$$ где константа $m_{2,0}$ выбирается так, чтобы выполнялось неравенство $m_{2,0}>\kappa \big(1,x_0,x_1\big).$ Прежде всего заметим, что функция $\kappa \big(\underline{\ \ },x_0,x_1 \big):\mathbb{N}\to\mathbb{N}$ определена, т.к. по условию $x_0$ и $x_1$ суть некоммутирующие гиперболические элементы (см. пункт~$V3)$ леммы~\ref{LemIvanPL}). Следовательно, слово $T_{\underline{x}[2,0]}(z_0,z_1)$ определено корректно.

Если слова $y_0$ и $y_{1}$ коммутируют, то из определения~(\ref{L2Def}) слова $L_2(z_0,z_1)$ легко следует, что $T_{\underline{x}[2,0]}(y_0,y_{1})=L_2(y_0,y_{1}^{m_{2,0}})=1$. Если же $y_0$ и $y_{1}$ суть некоммутирующие гиперболические элементы, то из пункта~$W2)$ леммы~$\ref{LemIvanSL}$ следует, что слово $T_{\underline{x}[2,0]}(y_0,y_{1})$ есть гиперболический элемент и не является собственной степенью в $H$. В частности, мы доказали, что $X_{2,0}=T_{\underline{x}[2,0]}(x_0,x_1)\in M_H$.

Пусть все пункты леммы были установлены для всех пар $(j',i')$, меньших пары $(j,i)$. Докажем лемму для пары $(j,i)$.

Чтобы слово $T_{\underline{x}[j,i]}(z_i,\dots,z_{i+j-1})$ было определено корректно необходимо и достаточно, чтобы функция $\kappa \big(\underline{\ \ },X_{j-1,i}^{m_{j,i-1}},X_{j-1,i+1} \big):\mathbb{N}\to\mathbb{N}$ была определена. В силу пункта~$V3)$ леммы~\ref{LemIvanPL} для этого достаточно, чтобы слова $X_{j-1,i}^{m_{j,i-1}}$ и $X_{j-1,i+1}$ были гиперболическими элементами и не коммутировали. Но это следует или из условия леммы, если $j=2$ (действительно, $X_{1,i}^{m_{2,i-1}}=x_i^{m_{2,i-1}}$ и $X_{1,i+1}=x_{i+1}$) или из пункта~$E1)$ и индукционного предположения, если $j>2$ (точнее, утверждение $X_{j-1,i}\in M_H$ доказывается для пары $(j-1,i)$; утверждения $X_{j-1,i+1}\in M_H$ и $[X_{j-1,i},X_{j-1,i+1}]\ne1$ доказывается для пары $(j-1,i+1)$). Следовательно, слово слово $T_{\underline{x}[j,i]}(z_i,\dots,z_{i+j-1})$ определено корректно.

В силу определения~(\ref{AddIII_NML_EQ1}) для кортежа $(y_i,\dots,y_{i+j-1})$ элементов из $H$ мы имеем:
\begin{equation}\label{AddIII_NML_EQ5}
T_{\underline{x}[j,i]}(y_i,\dots,y_{i+j-1})=L_2\Big(T_{\underline{x}[j-1,i]}\big(y_i,\dots,y_{i+j-2}\big)^{m_{j,i-1}},T_{\underline{x}[j-1,i+1]}\big(y_{i+1},\dots,y_{i+j-1}\big)^{m_{j,i}}\Big).
\end{equation}
Если $[y_k,y_{k+1}]=1$ для некоторого $k=i,\ \dots,\ i+j-2$, то, в силу индукционного предположения, хотя бы одно из слов $T_{\underline{x}[j-1,i]}\big(y_i,\dots,y_{i+j-2}\big)$ и  $T_{\underline{x}[j-1,i+1]}\big(y_{i+1},\dots,y_{i+j-1}\big)$ равно единице в $H$. Учитывая определение~(\ref{L2Def}) слова $L_2(z_0,z_1)$ и равенство~(\ref{AddIII_NML_EQ5}), легко видеть, что это влечет $T_{\underline{x}[j,i]}(y_i,\dots,y_{i+j-1})=1.$

Если $(y_i,\dots,y_{i+j-1})$ есть кортеж гиперболических элементов и слова $T_{\underline{x}[j-1,i]}\big(y_i,\dots,y_{i+j-2}\big)$ и $T_{\underline{x}[j-1,i+1]}\big(y_{i+1},\dots,y_{i+j-1}\big)$ не коммутируют, то, в силу пункта~$E3)$ и предположения индукции, оба эти слова суть гиперболические элементы. Это значит, что, в силу пункта~$W2)$ леммы~$\ref{LemIvanSL}$, слово $T_{\underline{x}[j,i]}(y_i,\dots,y_{i+j-1})$ есть гиперболический элемент и не является собственной степенью в $H$. Если же слова $T_{\underline{x}[j-1,i]}\big(y_i,\dots,y_{i+j-2}\big)$ и $T_{\underline{x}[j-1,i+1]}\big(y_{i+1},\dots,y_{i+j-1}\big)$ коммутируют, то из определения~(\ref{L2Def}) слова $L_2(z_0,z_1)$ и равенства~(\ref{AddIII_NML_EQ5}) легко следует, что $T_{\underline{x}[j,i]}(y_i,\dots,y_{i+j-1})=1$. Т.к., в силу пункта~$E1)$ и индукционного предположения, мы имеем $X_{j-1,i},X_{j-1,i+1}\in M_H$ и $[X_{j-1,i},X_{j-1,i+1}]\ne1$, то из сказанного легко следует, что $X_{i,j}=L_2\big(X_{j-1,i}^{m_{j,i-1}},X_{j-1,i+1}^{m_{j,i}}\big)\in M_H$.

%(действительно, т.к. по определению $X_{i,j}=L_2\beg(X_{j-1,i}^{m_{j,i-1}},X_{j-1,i+1}^{m_{j,i}}\big)$, достаточно воспользоваться пунктом $E1)$ и предположением индукции).

Если $i=0$, то на этом шаг индукции заканчивается. Если же $i>0$, то покажем еще что $[X_{j,i-1},X_{j,i}]\ne1.$ В силу выбора константы $m_{j,i}$, определения функции $\kappa \big(\underline{\ \ },X_{j-1,i}^{m_{j,i-1}},X_{j-1,i+1} \big)$ и пункта~$V3.2)$ леммы~\ref{LemIvanPL} (в обозначениях леммы~\ref{LemIvanPL}) мы имеем:
$$\big|B\big(X_{j-1,i}^{m_{j,i-1}},X_{j-1,i+1},m_{j,i}\big)\big|>|X_{j,i-1}|.$$
В силу пункта~$W2)$ леммы~$\ref{LemIvanSL}$ и оценки~(\ref{LemIvanSLEst2}) мы имеем:
\begin{align*}
|X_{j,i}|_r=\big|L_2\big(X_{j-1,i}^{m_{j,i-1}},X_{j-1,i+1}^{m_{j,i}}\big)\big|_c&>11136\big|B\big(X_{j-1,i}^{m_{j,i-1}},X_{j-1,i+1},m_{j,i}\big)^{l\big(X_{j-1,i}^{m_{j,i-1}},X_{j-1,i+1},m_{j,i}\big)}\big|\\
&\ge11136\big|B\big(X_{j-1,i}^{m_{j,i-1}},X_{j-1,i+1},m_{j,i}\big)\big|.
\end{align*}
Следовательно, выполнена оценка $|X_{j,i}|_r>|X_{j,i-1}|\ge|X_{j,i-1}|_r$, которая, в силу пункта~$C2)$ леммы~$\ref{MCL}$, влечет $[X_{j,i-1},X_{j,i}]\ne1.$
\end{proof}

\begin{Lemma}\label{LemTML1} Пусть группа $H=\zvezda\limits_{\mathfrak{a}\in \mathfrak{A}}\!H_\mathfrak{a}$ есть нетривиальное свободное произведение групп $H_\mathfrak{a}$, $\underline{x}=(x_0,\dots,x_{n-1})$ есть кортеж гиперболических элементов группы $H$, причем $[x_i,x_{i+1}]\ne1$ при $i=0,\ \dots,\ n-2$ и слова $T_{\underline{x}[j,i]}(z_i,\dots,z_{i+j-1})$ определены для кортежа $\underline{x}$ как~$(\ref{AddIII_NML_EQ1})$. Тогда для любой пары индексов $(j,i)$, где $j=2,\ \dots,\ n$ и $i=0,\ \dots,\ n-j$ из равенства:
    \begin{equation}\label{MainLemMainProp}
    T_{\underline{x}[j,i]}(x_i,\dots,x_{i+j-1})=T_{\underline{x}[j,i]}(y_0,\dots,y_{j-1}),
    \end{equation}
где $(y_0,\dots,y_{j-1})$ есть кортеж гиперболических элементов группы $H$ следует, что существует такой элемент $v\in\left \langle T_{\underline{x}[j,i]}(x_i,\dots,x_{i+j-1}) \right \rangle_{\infty}$, что:
$$x_i = y_0^v,\ \ x_{i+1}=y_1^v,\ \ \dots,\ \ x_{i+j-1}=y_{j-1}^v.$$
\end{Lemma}

\begin{proof} Для всех пар $(j,i)$ положим $X_{j,i}=T_{\underline{x}[j,i]}\big(x_i,\dots,x_{i+j-1}\big).$

Если существует такое $v\in H$, что выполнены равенства $x_{i+k}=y_k^v$, $k=0,\ \dots,\ j-1$, то из~(\ref{MainLemMainProp}) следует, что $X_{j,i}=X_{j,i}^v$. Если для некоторого элемента $v\in H$ выполнено равенство $X_{j,i}=X_{j,i}^v$, то, в силу того, что слово $X_{j,i}$ есть гиперболический элемент и не является собственной степенью в $H$ (см. пункты~$E1)$ и~$E3)$ леммы~\ref{AddIIINewLemIGProp}), из теоремы Куроша легко вывести, что $v\in\left \langle X_{j,i} \right \rangle_{\infty}$.

Докажем существование элемента $v$ для каждой пары $(j,i)$ проведя индукцию по $j$.

При $j=2$ и $i=0,\ \dots,\ n-2$, в силу условия~$(\ref{MainLemMainProp})$ и определения~$(\ref{AddIII_NML_EQ1})$ мы имеем:
$$L_2\big(x_i^{m_{2,i-1}},x_{i+1}^{m_{2,i}}\big)=L_2\big(y_0^{m_{2,i-1}},y_{1}^{m_{2,i}}\big),$$
где $m_{2,-1}=1$, $m_{2,0}>\kappa \big(1,x_0,x_{1}\big)$ и $m_{2,i}>\kappa \big(|X_{2,i-1}|,x_i^{m_{2,i-1}},x_{i+1}\big)$ при $i=1,\ \dots,\ n-2$. Т.к. функция $\kappa \big(\underline{\ \ },x_i^{m_{2,i-1}},x_{i+1}\big)$ монотонно возрастает и $|X_{2,i-1}|>1$ при $i=1,\ \dots,\ n-2$
(см. пункт~$E1)$ леммы~\ref{AddIIINewLemIGProp})
, то $m_{2,i}>\kappa \big(1,x_i^{m_{2,i-1}},x_{i+1}\big)$ при $i=0,\ \dots,\ n-2$ и утверждение следует из леммы~\ref{LemIvanML1} (в условии леммы~\ref{LemIvanML1} достаточно положить: $X_1=x_i^{m_{2,i-1}}$, $X_2=x_{i+1}$, $\hat{X}_1=y_0^{m_{2,i-1}}$, $\hat{X}_2=y_{1}$ и $k=m_{2,i}$) и однозначности извлечения корня из гиперболического элемента в $H$.

При $j=3$ и $i=0,\ \dots,\ n-3$ в силу условия~$(\ref{MainLemMainProp})$ и определения~$(\ref{AddIII_NML_EQ1})$ мы имеем:
\begin{equation}\label{MainLemMainPropEQ1}
L_2\Big(T_{\underline{x}[2,i]}\big(x_i,x_{i+1}\big)^{m_{3,i-1}},T_{\underline{x}[2,i+1]}\big(x_{i+1},x_{i+2}\big)^{m_{3,i}}\Big)=L_2\Big(T_{\underline{x}[2,i]}\big(y_0,y_{1}\big)^{m_{3,i-1}},T_{\underline{x}[2,i+1]}\big(y_{1},y_{2}\big)^{m_{3,i}}\Big),
\end{equation}
где $m_{3,-1}=1$, $m_{3,0}>\kappa \big(1,X_{2,0},X_{2,1} \big)$ и $m_{3,i}>\kappa \big(|X_{3,i-1}|,X_{2,i}^{m_{3,i-1}},X_{2,i+1}\big)$ при $i=1,\ \dots,\ n-3$. Т.к. функция $\kappa \big(\underline{\ \ },X_{2,i}^{m_{3,i-1}},X_{2,i+1} \big)$ монотонно возрастает и $|X_{3,i-1}|>1$ при $i=1,\ \dots,\ n-3$ (см. пункт~$E1)$ леммы~\ref{AddIIINewLemIGProp}), то $m_{3,i}>\kappa \big(1,X_{2,i}^{m_{3,i-1}},X_{2,i+1}\big)$ при $i=0,\ \dots,\ n-3$ и лемма~\ref{LemIvanML1} применима к~(\ref{MainLemMainPropEQ1}), т.е. существует такой элемент $s\in H$, что:
$$T_{\underline{x}[2,i]}\big(x_i,x_{i+1}\big)^{m_{3,i-1}}=T_{\underline{x}[2,i]}\big(y_0^s,y_{1}^s\big)^{m_{3,i-1}}\ \ \text{и}\ \ T_{\underline{x}[2,i+1]}\big(x_{i+1},x_{i+2}\big) = T_{\underline{x}[2,i+1]}\big(y_{1}^s,y_{2}^s\big).$$
Т.к. $X_{2,i}$ есть гиперболический элемент (см. пункт~$E1)$ леммы~\ref{AddIIINewLemIGProp}) и корень из гиперболического элемента в $H$ извлекается однозначно, то из полученных равенств следует, что:
\begin{equation}\label{MainLemMainPropEQ3}
T_{\underline{x}[2,i]}\big(x_i,x_{i+1}\big)=T_{\underline{x}[2,i]}\big(y_0^s,y_{1}^s\big)\ \ \text{и}\ \ T_{\underline{x}[2,i+1]}\big(x_{i+1},x_{i+2}\big) = T_{\underline{x}[2,i+1]}\big(y_{1}^s,y_{2}^s\big).
\end{equation}
Значит, в силу базы индукции, найдутся такие элементы $u_1,u_2\in H$, что:
\begin{equation}\label{MainLemMainPropEQ2}
x_i = y_0^{su_1}\ \ \text{и}\ \ x_{i+1}=y_1^{su_1}\ \ \text{и}\ \ x_{i+1}=y_1^{su_2}\ \ \text{и}\ \ x_{i+2}=y_2^{su_2},
\end{equation}
т.е. нам осталось показать, что $u_1=u_2$ (действительно, если это так, то мы можем положить $v\overset{\textrm{def}}{=}su_1$).

Выражая $y_0$, $y_1$ и $y_2$ из~(\ref{MainLemMainPropEQ2}) и подставляя полученные выражения в~(\ref{MainLemMainPropEQ3}), мы получаем:
\begin{equation}\label{MainLemMainPropEQ4}
T_{\underline{x}[2,i]}\big(x_i,x_{i+1}\big)=T_{\underline{x}[2,i]}\big(x_i,x_{i+1}\big)^{u_1^{-1}}\ \ \text{и}\ \ T_{\underline{x}[2,i+1]}\big(x_{i+1},x_{i+2}\big) = T_{\underline{x}[2,i+1]}\big(x_{i+1},x_{i+2}\big)^{u_2^{-1}}.
\end{equation}
Т.к. (в силу пунктов~$E1)$ и~$E3)$ леммы~\ref{AddIIINewLemIGProp}) слова $X_{2,i}=T_{\underline{x}[2,i]}\big(x_i,x_{i+1}\big)$ и $X_{2,i+1}=T_{\underline{x}[2,i+1]}\big(x_{i+1},x_{i+2}\big)$ являются гиперболическими элементами и не являются собственными степенями в $H$, то из теоремы Куроша легко вывести, что из равенств~(\ref{MainLemMainPropEQ4}) следует, что:
\begin{equation}\label{MainLemMainPropEQ5}
u_1=T_{\underline{x}[2,i]}\big(x_i,x_{i+1}\big)^{-l_1}\ \ \text{и}\ \ u_2=T_{\underline{x}[2,i+1]}\big(x_{i+1},x_{i+2}\big)^{l_2}
\end{equation}
для некоторых $l_1,l_2\in\mathbb{Z}.$ Если $l_1=l_2=0$, то $u_1=u_2=1$ и утверждение доказано; пусть $l_1$ и $l_2$ одновременно не равны нулю.

Из равенств~(\ref{MainLemMainPropEQ2}) мы имеем $x_{i+1}=x_{i+1}^{u_1^{-1}u_2}$. Пусть $u_1\ne u_2$, тогда из последнего равенства следует, что найдутся такие $l_3,l_4\in\mathbb{Z}\!\setminus\!\{0\}$, что $(u_1^{-1}u_2)^{l_3}=x_{i+1}^{l_4}$. Подставляя в это равенство выражения для $u_1$ и $u_2$ из~(\ref{MainLemMainPropEQ5}) и учитывая, что $T_{\underline{x}[2,i]}\big(x_i,x_{i+1}\big)=L_2\big(x_i^{m_{2,i-1}},x_{i+1}^{m_{2,i}}\big)$ и $T_{\underline{x}[2,i+1]}\big(x_{i+1},x_{i+2}\big)=L_2\big(x_{i+1}^{m_{2,i}},x_{i+2}^{m_{2,i+1}}\big)$, мы получаем равенство:
\begin{equation}\label{MainLemMainPropEQ6}
\left ( L_2\big(x_i^{m_{2,i-1}},x_{i+1}^{m_{2,i}}\big)^{l_1}L_2\big(x_{i+1}^{m_{2,i}},x_{i+2}^{m_{2,i+1}}\big)^{l_2} \right )^{l_3}=x_{i+1}^{l_4}.
\end{equation}
Полагая $X_0 = x_i^{m_{2,i-1}}$, $X_1=x_{i+1}^{m_{2,i}}$, $X_2=x_{i+2}$, $m=m_{2,i+1}$ и возводя равенство~(\ref{MainLemMainPropEQ6}) в степень $m_{2,i}$ мы получаем:
\begin{equation}\label{MainLemMainPropEQ7}
\left ( L_2\big(X_0,X_1\big)^{l_1}L_2\big(X_1,X_2^m\big)^{l_2} \right )^{l_3m_{2,i}}=X_1^{l_4}.
\end{equation}
Не сложно проверить, что в новых обозначениях условие~(\ref{AddIII_NML_EQ2T22}) для $m_{2,i+1}$ записывается как $m>\kappa \big(|L_2(X_0,X_1)|,X_1,X_2\big)$. Значит из равенства~(\ref{MainLemMainPropEQ7}), в силу леммы~$\ref{LemICAdd2}$, следует, что $l_3m_{2,i}=0$, т.е. $l_3=0$ (т.к. $m_{2,i}\in\mathbb{N}$) и мы получили противоречие с предположением $u_1\ne u_2$.

При $j\ge4$ и $i=0,\ \dots,\ n-j$ в силу равенства~$(\ref{MainLemMainProp})$ и определения~$(\ref{AddIII_NML_EQ1})$ мы имеем:
\begin{align}\label{MainLemMainPropEQ8}
L_2\Big(T_{\underline{x}[j-1,i]}\big(x_i,\dots,x_{i+j-2}&\big)^{m_{j,i-1}},T_{\underline{x}[j-1,i+1]}\big(x_{i+1},\dots,x_{i+j-1}\big)^{m_{j,i}}\Big)=\\\nonumber
&=L_2\Big(T_{\underline{x}[j-1,i]}\big(y_0,\dots,y_{j-2}\big)^{m_{j,i-1}},T_{\underline{x}[j-1,i+1]}\big(y_{1},\dots,y_{j-1}\big)^{m_{j,i}}\Big),
\end{align}
где $m_{j,-1}=1$, $m_{j,0}>\kappa \big(1,X_{j-1,i},X_{j-1,i+1} \big)$ и $m_{j,i}>\kappa \big(|X_{j,i-1}|,X_{j-1,i}^{m_{j,i-1}},X_{j-1,i+1}\big)$ при $i=1,\ \dots,\ n-j$. Т.к. функция $\kappa \big(\underline{\ \ },X_{j-1,i}^{m_{j,i-1}},X_{j-1,i+1} \big)$ монотонно возрастает и $|X_{j,i-1}|>1$ при $i=1,\ \dots,\ n-j$ (см. пункт~$E1)$ леммы~\ref{AddIIINewLemIGProp}), то $m_{j,i}>\kappa \big(1,X_{j-1,i}^{m_{j,i-1}},X_{j-1,i+1}\big)$ при $i=0,\ \dots,\ n-j$ и лемма~\ref{LemIvanML1} применима к~(\ref{MainLemMainPropEQ8}), т.е. существует такое $s\in H$, что выполнены равенства:
\begin{align*}
T_{\underline{x}[j-1,i]}\big(x_i,\dots,x_{i+j-2}\big)^{m_{j,i-1}}&=T_{\underline{x}[j-1,i]}\big(y_0^s,\dots,y_{j-2}^s\big)^{m_{j,i-1}};\\
T_{\underline{x}[j-1,i+1]}\big(x_{i+1},\dots,x_{i+j-1}\big)&=T_{\underline{x}[j-1,i+1]}\big(y_{1}^s,\dots,y_{j-1}^s\big).
\end{align*}
Т.к. $X_{j-1,i}$ есть гиперболический элемент (см. пункт~$E1)$ леммы~\ref{AddIIINewLemIGProp}) и корень из гиперболического элемента в $H$ извлекается однозначно, то из полученных равенств, в силу индукционного предположения следует, что существуют такие элементы $u_1,u_2\in H$, что:
\begin{align*}
x_i=y_0^{su_1},\ \ x_{i+1}&=y_1^{su_1},\ \ \dots,\ \ x_{i+j-2}=y_{j-2}^{su_1};\\
x_{i+1}&=y_1^{su_2},\ \ \dots,\ \ x_{i+j-2}=y_{j-2}^{su_2},\ \ x_{i+j-1}=y_{j-1}^{su_2}.
\end{align*}
Т.е. нам осталось показать, что $u_1=u_2$ (действительно, если это так, то мы можем положить $v\overset{\textrm{def}}{=}su_1$). Т.к. $j\ge4$, то из полученных равенств мы имеем:
$$x_{i+1}^{(su_1)^{-1}}=y_1=x_{i+1}^{(su_2)^{-1}}\ \ \text{и}\ \ x_{i+2}^{(su_1)^{-1}}=y_2=x_{i+2}^{(su_2)^{-1}},$$
из чего следует, что $u_1^{-1}u_2\in C_H(x_{i+1})\cap C_H(x_{i+2}).$ Т.к. $x_{i+1}$ и $x_{i+2}$ суть некоммутирующие гиперболические элементы группы $H$, то из теоремы Куроша несложно вывести, что $C_H(x_{i+1})\cap C_H(x_{i+2})=\{1\}$.
\end{proof}

\begin{ProofMLII} $a)$ В силу леммы~\ref{LemTML1} мы можем положить:
$$T_{\underline{x}}(z_0,\dots,z_{n-1})\overset{\textrm{def}}{=}T_{\underline{x}[n,0]}(z_0,\dots,z_{n-1}),$$
где слово $T_{\underline{x}[n,0]}(z_0,\dots,z_{n-1})$ определено для кортежа $\underline{x}$ как~$(\ref{AddIII_NML_EQ1}).$
%Утверждение немедленно следует из леммы \ref{LemTML1} при $(j,i)=(n,0)$ (действительно, $\underline{x}[n,0]=\underline{x}$ в силу используемых нами обозначений).

$b)$ Кортежу $\underline{x}=(x_0,x_1,\dots,x_{n-2},x_{n-1})$ сопоставим кортеж (длины $2n-1$):
$$\underline{\hat{x}}=(x_0,x_1,\dots,x_{n-2},x_{n-1},x_{n-2},\dots,x_1,x_0).$$
Легко видеть, что кортеж $\underline{\hat{x}}$ удовлетворяет условиям леммы~\ref{LemTML1} (т.е. все его компоненты суть гиперболические элементы группы $H$ и соседние компоненты не коммутируют). Положим:
\begin{align}
T_{\underline{x}}'(z_0,\dots,z_{n-1})&\overset{\textrm{def}}{=}T_{\underline{\hat{x}}[2n-2,0]}(z_0,z_1,\dots,z_{n-2},z_{n-1},z_{n-2},\dots,z_2,z_1);\label{DefT'}\\
T_{\underline{x}}''(z_0,\dots,z_{n-1})&\overset{\textrm{def}}{=}T_{\underline{\hat{x}}[2n-2,1]}(z_1,z_2,\dots,z_{n-2},z_{n-1},z_{n-2},\dots,z_1,z_0),\label{DefT''}
\end{align}
где слова $T_{\underline{\hat{x}}[2n-2,0]}(z_0,\dots,z_{2n-3}),T_{\underline{\hat{x}}[2n-2,1]}(z_1,\dots,z_{2n-2})\in F_{2n-1}(z_0,\dots,z_{2n-2})'$ определены для кортежа $\underline{\hat{x}}$ как~$(\ref{AddIII_NML_EQ1})$. Обозначим $X'=T_{\underline{x}}'(x_0,\dots,x_{n-1})$ и $X''=T_{\underline{x}}''(x_0,\dots,x_{n-1})$. Из пункта~$E1)$ леммы~\ref{AddIIINewLemIGProp} следует, что $X'$ и $X''$ суть некоммутирующие гиперболические элементы (действительно, в обозначениях леммы~\ref{AddIIINewLemIGProp} мы имеем $X'=X_{2n-2,0}$ и $X''=X_{2n-2,1}$). Следовательно, мы имеем $C_H(X')\cap C_H(X'')=\{1\}$.

Пусть выполнены равенства~(\ref{LemMLPartIIaddTSS}), тогда из определений~(\ref{DefT'}) и~(\ref{DefT''}) и леммы~\ref{LemTML1} следует, что найдутся такие элемемнты $v_1\in\langle X'\rangle_{\infty}$ и $v_2\in\langle X''\rangle_{\infty}$, что выполнены равенства:
\begin{align*}
x_i=y_i^{v_1}\ \ \text{и}\ \ x_i=y_i^{v_2}\ \ \text{при}\ \ i=0,\ \dots,\ n-1.
\end{align*}
Т.к. по условию $x_0$ и $x_1$ суть некоммутирующие гиперболические элементы, то из полученных равенств следует, что $v_1^{-1}v_2\in C_H(x_0)\cap C_H(x_1)=\{1\}$, т.е. $v_1=v_2$. Но т.к. $v_1\in\langle X'\rangle_{\infty}$, $v_2\in\langle X''\rangle_{\infty}$ и $C_H(X')\cap C_H(X'')=\{1\}$, то из равенства $v_1=v_2$ следует, что $v_1=v_2=1$.
\end{ProofMLII}

\section{Доказательство Леммы о словах, часть II}

Основной целью данного раздела является доказательство следующей леммы, обобщающей лемму~\ref{LemMLPartII} (элементами кортежа $(y_0,\dots,y_{n-1})$ теперь могут быть произвольные элементы группы $H$, а не только гиперболические элементы, как в лемме~\ref{LemMLPartII}).
\begin{Lemma}\label{LemMLPartIII} Пусть группа $H=\zvezda\limits_{\mathfrak{a}\in \mathfrak{A}}\!H_\mathfrak{a}$ есть нетривиальное свободное произведение групп $H_\mathfrak{a}$ и пусть $\underline{x}=(x_0,\dots,x_{n-1})$, $n\ge2$ есть такой кортеж гиперболических элементов группы $H$, что $[x_i,x_{i+1}]\ne1$ при $i=0,\ \dots,\ n-2$. Тогда:
\begin{enumerate}
     \item[\sffamily $\mathrm{a)}$] существует слово $P_{\underline{x}}(z_0,\dots,z_{n-1})\in F_n(z_0,\dots,z_{n-1})'$, вообще говоря, зависящее от кортежа~$\underline{x}$, такое, что если для некоторого кортежа $(y_0,\dots,y_{n-1})$ элементов группы $H$ выполнено равенство:
\begin{equation}\label{LemMLPartIIIMEQ}
P_{\underline{x}}(x_0,\dots,x_{n-1})=P_{\underline{x}}(y_0,\dots,y_{n-1}),
\end{equation}
то найдется такой элемент $v\in \left \langle P_{\underline{x}}(x_0,\dots,x_{n-1}) \right \rangle_\infty$, что $x_i = y_i^v$ при $i=0,\ \dots,\ n-1$;

     \item[\sffamily $\mathrm{b)}$] существуют слова $P_{\underline{x}}'(z_0,\dots,z_{n-1}), P_{\underline{x}}''(z_0,\dots,z_{n-1})\in F_n(z_0,\dots,z_{n-1})'$, вообще говоря, зависящие от кортежа $\underline{x}$, такие, что если для некоторого кортежа $(y_0,\dots,y_{n-1})$ элементов группы $H$ одновременно выполнены равенства:
\begin{equation}\label{LemMLPartIIIMEQpartb}
P_{\underline{x}}'(x_0,\dots,x_{n-1})=P_{\underline{x}}'(y_0,\dots,y_{n-1})\ \ \text{и}\ \ P_{\underline{x}}''(x_0,\dots,x_{n-1})=P_{\underline{x}}''(y_0,\dots,y_{n-1}),
\end{equation}
то $x_i = y_i$ при $i=0,\ \dots,\ n-1$.
 \end{enumerate}
\end{Lemma}
%\begin{Lemma}\label{LemMLPartIII} Пусть группа $H=\zvezda\limits_{\mathfrak{a}\in \mathfrak{A}}\!H_\mathfrak{a}$ есть нетривиальное свободное произведение групп $H_\mathfrak{a}$ и %пусть $\underline{x}=(x_0,\dots,x_{n-1})$ есть такой кортеж гиперболических элементов группы $H$, что $[x_i,x_{i+1}]\ne1$ при $i=0,\ \dots,\ n-2$. Тогда существует слово %$P_{\underline{x}}(z_0,\dots,z_{n-1})\in F_n(z_0,\dots,z_{n-1})'$, вообще говоря, зависящее от кортежа $\underline{x}$ такое, что если для некоторого кортежа $(y_0,\dots,y_{n-1})$ %элементов $y_i\in H$ выполнено равенство:
%\begin{equation}\label{LemMLPartIIIMEQ}
%P_{\underline{x}}(x_0,\dots,x_{n-1})=P_{\underline{x}}(y_0,\dots,y_{n-1}),
%\end{equation}
%то найдется такой элемент $v\in \left \langle P_{\underline{x}}(x_0,\dots,x_{n-1}) \right \rangle_\infty$, что $x_i = y_i^v$ при $i=0,\ \dots,\ n-1.$
%\end{Lemma}

%Единственное отличие данной леммы от леммы \ref{LemMLPartII} заключается в отсутствии каких либо условий, наложенных на элементы $y_i$.

Прежде чем перейти к доказательству леммы~\ref{LemMLPartIII}, установим ряд вспомогательных утверждений.

Определим слова $E_n(z_0,\dots,z_{n-1}) \in F_n(z_0,\dots,z_{n-1})'$ для $n\ge 2$ следующим образом:
\begin{equation}\label{AddIIIE_NTWE_Q1}
E_2(z_0,z_1)\overset{\textrm{def}}{=}[z_0^{2},z_1^{2}]\ \  \text{и}\ \ E_n(z_0,\dots,z_{n-1}) \overset{\textrm{def}}{=}\bigl[E_{n-1}(z_0,\dots,z_{n-2})^2,z_{n-1}^2\bigr]\ \ \text{при}\ \ n\ge3.
\end{equation}

Для подкортежей\footnote[$\dagger$]{См. определение подкортежа на стр. \pageref{pagewith1}.} $\underline{k}[j,0] = (k_0,\dots,k_{j-1})$, $j=1,\ \dots,\ n$ кортежа $\underline{k}=\underline{k}[n,0]=(k_0,\dots,k_{n-1})\in \mathbb{N}^n$, $n\ge 1$ определим слова $E_{\underline{k}[j,0]}(z_0,\dots,z_{j-1})\in F_n(z_0,\dots,z_{n-1})'$ следующим образом:
\begin{equation}\label{AddIII_NTWE_Q2}
E_{\underline{k}[1,0]}(z_0)\overset{\textrm{def}}{=}z_0^{k_0}\ \ \text{и}\ \ E_{\underline{k}[j,0]}(z_0,\dots,z_{j-1}) \overset{\textrm{def}}{=}\bigl[E_{\underline{k}[j-1,0]}(z_0,\dots,z_{j-2})^2,z_{j-1}^{2k_{j-1}}\bigr]\ \ \text{при}\ \ j\ge2.
\end{equation}

%Для кортежа $\underline{k}=(k_0,\dots,k_{n-1})\in \mathbb{N}^n$, $n\ge 2$ определим слово $E_{\underline{k}}(z_0,\dots,z_{n-1})$ как:
В частности, из определений~(\ref{AddIIIE_NTWE_Q1}) и~(\ref{AddIII_NTWE_Q2}) следует, что для подкортежа $\underline{k}[j,0]$ кортежа $\underline{k}=(k_0,\dots,k_{n-1})\in \mathbb{N}^n$ при $j\ge 2$ мы имеем:
$$E_{\underline{k}[j,0]}(z_0,\dots,z_{j-1}) = E_j(z_0^{k_0},\dots,z_{j-1}^{k_{j-1}}).$$
%В частности, из (\ref{AddIIIE_NTWE_Q1}) и (\ref{AddIII_NTWE_Q2}) при $n\ge 2$ следует, что:
%$$E_{\underline{k}}(z_0,\dots,z_{n-1})=[E_{\underline{k}[n-1,0]}(z_0,\dots,z_{n-2})^2,z_{n-1}^{2k_{n-1}}],$$
%где $E_{\underline{k}[1,0]}(z_0,\dots,z_{n-2})\overset{\textrm{def}}{=}z_0^{k_0}$ и $\underline{k}[j,0] = (k_0,k_1,\dots,k_{j-1})$ есть подкортеж кортежа $\underline{k}$ (подробнее %см. стр. \pageref{pagewith1}).

\begin{Lemma}\label{AddIIILem1} Пусть группа $H=\zvezda\limits_{\mathfrak{a}\in \mathfrak{A}}\!H_\mathfrak{a}$ есть нетривиальное свободное произведение групп $H_\mathfrak{a}$ и пусть $\underline{k}=\underline{k}[n,0]=(k_0,\dots,k_{n-1})\in \mathbb{N}^n$, $n\ge2$ есть произвольный кортеж. Тогда, если кортеж $(y_0,\dots,y_{n-1})$ элементов из $H$ таков, что $E_{\underline{k}}(y_0,\dots,y_{n-1})\in H_\mathfrak{a}^g\!\setminus\!\{1\}$ для некоторого $\mathfrak{a}\in \mathfrak{A}$ и $g\in H$, то $y_0,\ \dots,\ y_{n-1}\in H_\mathfrak{a}^g\!\setminus\!\{1\}$.
\end{Lemma}
\begin{proof} Индукция по длине кортежа $\underline{k}$. При $n=2$, в силу определения~(\ref{AddIII_NTWE_Q2}) и условия леммы, мы имеем $E_{\underline{k}}(y_0,y_{1})=[y_0^{2k_0},y_1^{2k_1}]\in H_\mathfrak{a}^g\!\setminus\!\{1\}$. Значит, в силу пункта~$C1)$ леммы~\ref{MCL}, $y_0^{2k_0},y_{1}^{2k_1}\in H_\mathfrak{a}^g\!\setminus\!\{1\},$ т.е. $y_0,y_{1}\in H_\mathfrak{a}^g\!\setminus\!\{1\}.$

Если $n\ge 3$, то, в силу определения~(\ref{AddIII_NTWE_Q2}) и условия леммы, мы имеем:
$$E_{\underline{k}}(y_0,\dots,y_{n-1}) = \bigl[E_{\underline{k}[n-1,0]}(y_0,\dots,y_{n-2})^2,y_{n-1}^{2k_{n-1}}\bigr]\in H_\mathfrak{a}^g\!\setminus\!\{1\}.$$
Из пункта~$C1)$ леммы~\ref{MCL} легко следует, что в таком случае мы имеем $E_{\underline{k}[n-1,0]}(y_0,\dots,y_{n-2}),\ y_{n-1}\in H_\mathfrak{a}^g\!\setminus\!\{1\}$. Осталось применить к $E_{\underline{k}[n-1,0]}(y_0,\dots,y_{n-2})\in H_\mathfrak{a}^g\!\setminus\!\{1\}$ предположение индукции.
%Значит, в силу пункта $C1)$ леммы \ref{MCL}, $E_{\underline{k}[n-1,0]}(y_1,\dots,y_{n-2})^2,\ y_{n-1}^{2k_{n-1}}\in H_\mathfrak{a}^g\!\setminus\!\{1\}$, следовательно, %$E_{\underline{k}[n-1,0]}(y_1,\dots,y_{n-2}),\ y_{n-1}\in H_\mathfrak{a}^g\!\setminus\!\{1\}$ и осталось применить к $E_{\underline{k}[n-1,0]}(y_1,\dots,y_{n-2})$ предположение %индукции.
\end{proof}

\begin{Lemma}\label{AddIIILem2mod} Пусть группа $H=\zvezda\limits_{\mathfrak{a}\in \mathfrak{A}}\!H_\mathfrak{a}$ есть нетривиальное свободное произведение групп $H_\mathfrak{a}$ и пусть $(x_0,\dots,x_{n-1})$, $n\ge 2$ есть такой кортеж гиперболических элементов группы $H$, что $[x_i,x_{i+1}]\ne 1$ при $i=0,\ \dots,\ n-2$. Тогда константы $k_0,\ \dots,\ k_{n-2}\in\mathbb{N}$ могут быть выбраны так, что:
\begin{equation}\label{AddIIILem2EQ1}
\bigl|E_n(x_0^{k_0},\dots,x_{n-2}^{k_{n-2}},x_{n-1}^{m})\bigr|_r\to\infty\ \ \text{при}\ \ m\to \infty.
\end{equation}
\end{Lemma}
\begin{proof} Индукция по $n$. При $n=2$ положим $k_0=1$, тогда, в силу определения~(\ref{AddIIIE_NTWE_Q1}), мы имеем $E_2(x_0,x_1^{m})=[x_0^2,x_1^{2m}]$ и утверждение немедленно следует из пункта~$C5)$ леммы~\ref{MCL}.

Пусть $n\ge3$. В соответствии с определением~(\ref{AddIIIE_NTWE_Q1}), мы имеем:
\begin{equation}\label{AddIIILem2EQ3}
E_n(x_0^{k_0},\dots,x_{n-2}^{k_{n-2}},x_{n-1}^{m}) = \bigl[E_{n-1}(x_0^{k_0},\dots,x_{n-2}^{k_{n-2}})^2,x_{n-1}^{2m}\bigr].
\end{equation}
Опираясь на индуктивное предположение для $n-1,$ выберем константы $k_0,\ \dots,\ k_{n-3}$ так, чтобы выполнялось:
\begin{equation}\label{AddIIILem2EQ2}
\bigl|E_{n-1}(x_0^{k_0},\dots,x_{n-3}^{k_{n-3}},x_{n-2}^{k_{n-2}})\bigr|_r\to\infty\ \ \text{при}\ \ k_{n-2}\to \infty.
\end{equation}
Из~(\ref{AddIIILem2EQ2}) следует, что мы можем выбрать такое $k_{n-2}$, что $\bigl[E_{n-1}(x_0^{k_0},\dots,x_{n-3}^{k_{n-3}},x_{n-2}^{k_{n-2}}),x_{n-1}\bigr]\ne\!~1$  (действительно, в силу пункта~$C2)$ леммы~\ref{MCL}, для этого достаточно выбрать $k_{n-2}$ так, чтобы $\bigl|E_{n-1}(x_0^{k_0},\dots,x_{n-3}^{k_{n-3}},x_{n-2}^{k_{n-2}})\bigr|_r>|x_{n-1}|_r$). Теперь нужное нам утверждение непосредственно следует из пункта~$C5)$ леммы~\ref{MCL}.
\end{proof}

Для кортежа $\underline{k}=(k_0,\dots,k_{n-2},m_1,m_2)\in \mathbb{N}^{n+1}$, $n\ge 2$ определим слово $J_{\underline{k}}(z_0,\dots,z_{n-1},z_n) \in F_{n+1}(z_0,\dots,z_{n})'$ следующим образом:
\begin{equation}\label{AddIIILemNWEQ3}
J_{\underline{k}}(z_0,\dots,z_{n-1},z_{n})\overset{\textrm{def}}{=}\bigl[E_{\underline{k}[n,0]}(z_0,\dots,z_{n-1})^{2m_2},z_n^2\bigr].
\end{equation}

%Введем следующие слова из группы $F_\infty(z_0,z_1,\dots)'$:
%\begin{align}
%T_i(z_0,\dots,z_{n-1})&=[z_i^3,E_n(z_0^{k_0},z_1^{k_1},\dots,z_{n-2}^{k_{n-2}},z_{n-1}^{m})^3];\ \ \text{для}\ \ i = 0,\ \dots,\ n-1\label{AddIIINWT1}\\
%\hat{T}_i(z_0,\dots,z_{n-1})&=[z_i^{3},E_n(z_0^{k_0},z_1^{k_1},\dots,z_{n-2}^{k_{n-2}},z_{n-1}^{\hat{m}})^3],\ \ \text{для}\ \ i = 0,\ \dots,\ n-1\label{AddIIINWT2}
%\end{align}
%где константы $k_0,\ \dots,\ k_{n-2}$, $m$ и $\hat{m}$ суть некоторые натуральные числа.

\begin{Lemma}\label{AddIIILem2} Пусть группа $H=\zvezda\limits_{\mathfrak{a}\in \mathfrak{A}}\!H_\mathfrak{a}$ есть нетривиальное свободное произведение групп $H_\mathfrak{a}$ и пусть $(x_0,\dots,x_{n-1})$, $n\ge 2$ есть такой кортеж гиперболических элементов группы $H$, что $[x_i,x_{i+1}]\ne 1$ при $i=0,\ \dots,\ n-2$. Тогда найдутся такие кортежи $\underline{k}=(k_0,\dots,k_{n-2},m_1,m_2)\in\mathbb{N}^{n+1}$ и $\underline{\hat{k}}=(k_0,\dots,k_{n-2},\hat{m}_1,\hat{m}_2)\in\mathbb{N}^{n+1}$, что выполнены следующие условия:
 \begin{enumerate}
     \item[\sffamily $\mathrm{a)}$] $\underset{i}{\min}\Big( \big|J_{\underline{k}}(x_0,\dots,x_{n-1},x_i)\big|_r\Big )>\big|E_{\underline{\hat{k}}[n,0]}(x_0,\dots,x_{n-1})\big|_r>\big|E_{\underline{k}[n,0]}(x_0,\dots,x_{n-1})\big|_r$;

     \item[\sffamily $\mathrm{b)}$]
        $\underset{i}{\min}\Big (\big|J_{\underline{\hat{k}}}(x_0,\dots,x_{n-1},x_i)\big|_r\Big)>\underset{i}{\max}\Big (\big|J_{\underline{k}}(x_0,\dots,x_{n-1},x_i)\big|_r \Big ).$
 \end{enumerate}
\end{Lemma}
\begin{proof} Т.к. по условию $\underline{k}[n-1,0]=\underline{\hat{k}}[n-1,0]=(k_0,\dots,k_{n-2})$, то нам необходимо определить $n+3$ константы: $k_0,\ \dots,\ k_{n-2},\ m_1,\ \hat{m}_1,\ m_2\ \text{и}\ \hat{m}_2$.

$a)$ Положим $M=\max\big(|x_0|_r,\dots,|x_{n-1}|_{r}\big)$ и, воспользовавшись леммой~\ref{AddIIILem2mod}, выберем константы $k_0,\ \dots,\ k-2$ так, чтобы для кортежа $(x_0,\dots,x_{n-1})$ было выполнено~(\ref{AddIIILem2EQ1}), т.е. чтобы:
$$\bigl|E_{\underline{k}[n,0]}(x_0,\dots,x_{n-1})\bigr|_r=\bigl|\bigl[E_{\underline{k}[n-1,0]}(x_0,\dots,x_{n-2})^2,x_{n-1}^{2m_1}\bigr]\bigr|_r\to\infty\ \ \text{при}\ \ m_1\to \infty.$$
Фиксируем такое $m_1$, что $|E_{\underline{k}[n,0]}(x_0,\dots,x_{n-1})|_r>M$. После того, как константа $m_1$ фиксирована, вновь опираясь на~(\ref{AddIIILem2EQ1}), выберем константу $\hat{m}_1$ так, чтобы выполнялось неравенство $|E_{\underline{\hat{k}}[n,0]}(x_0,\dots,x_{n-1})|_r>|E_{\underline{k}[n,0]}(x_0,\dots,x_{n-1})|_r$ (мы можем это сделать, т.к. на этом этапе константы $k_0,\ \dots,\ k_{n-2}$,\ $m_1$ уже фиксированы и, следовательно, величина $|E_{\underline{k}[n,0]}(x_0,\dots,x_{n-1})|_r$ так же является фиксированной константой).

Пусть константы $m_1$ и $\hat{m}_1$ фиксированы. Покажем, что для любого $i = 0,\ \dots,\ n-1$ выполнено:
\begin{align}
\bigl|J_{\underline{k}}(x_0,\dots,x_{n-1},x_i)\bigr|_r&\to\infty\ \ \text{при}\ \ m_2\to\infty;\label{AddIIILemNWEQ4_1}\\
\bigl|J_{\underline{\hat{k}}}(x_0,\dots,x_{n-1},x_i)\bigr|_r&\to\infty\ \ \text{при}\ \ \hat{m}_2\to\infty.\label{AddIIILemNWEQ4_2}
\end{align}
В соответствии с определением~(\ref{AddIIILemNWEQ3}), для любого $i$ мы имеем:
$$J_{\underline{k}}(x_0,\dots,x_{n-1},x_i)=\bigl[E_{\underline{k}[n,0]}(x_0,\dots,x_{n-1})^{2m_2},x_i^2\bigr].$$
Т.к. $|E_{\underline{k}[n,0]}(x_0,\dots,x_{n-1})|_r>M\ge|x_i|_r$, то легко видеть, что~(\ref{AddIIILemNWEQ4_1}) непосредственно следует из последнего равенства и пункта~$C5)$ леммы~\ref{MCL} (для~(\ref{AddIIILemNWEQ4_2}) доказательство аналогичное). Значит, опираясь на~(\ref{AddIIILemNWEQ4_1}), мы можем выбрать такое $m_2$, чтобы было выполнено первое неравенство пункта~$a)$.

$b)$ После того как константы $m_1$, $\hat{m}_1$ и $m_2$ фиксированы, используя~(\ref{AddIIILemNWEQ4_2}), мы можем выбрать такое $\hat{m}_2$, чтобы было выполнено неравенство пункта~$b)$.
\end{proof}

\begin{ProofMLIII} Пусть $\underline{z}=(z_0,\dots,z_{n-1})$.
%Пусть группа $H=H_1\ast H_2$ есть свободное произведение нетривиальных групп $H_1$ и $H_2$, $H$ есть вербально замкнутая подгруппа в группе $G$ и $X$ есть порождающее множество %группы $H$, удовлетворяющее удовлетворяющее пунктам $1)$ и $2)$ леммы \ref{ProofMTBL}.
%
%В соответствии с леммой \ref{ProofMT2TL}, достаточно доказать, что произвольная система уравнений вида:
%$$S=\{w_i(\underline{z})=x_i\ |\ w_i(\underline{z})\in F(\underline{z}),\ x_i\in X,\ i=0,\ \dots,\ n-1\},$$ имеющая решение в $G$, имеет решение и в $H$. Далее мы считаем, что $n\ge %2$, т.к. при $n=1$ утверждение непосредственно следует из определения вербальной замкнутости.
%Для удобства мы фиксируем порядок, в котором записаны уравнения системы $S$, т.к. далее мы будем рассматривать кортеж $\underline{x}=(x_0,\dots,x_{n-1})$, составленный из констант, %являющихся правыми частями уравнений этой системы. Покажем, что без ограничения общности мы можем считать, что в кортеже $\underline{x}$ для любого $i=0,\ \dots,\ n-2$ выполнено %условие $[x_i,x_{i+1}]\ne 1$. Действительно, если $[x_i,x_{i+1}]= 1$ для некоторого $i$, то (в силу пункта 2) леммы \ref{ProofMTBL}) найдется такой элемент $x_j\in X$, что %$[x_i,x_j]\ne1$ и $[x_j,x_{i+1}]\ne 1.$ Между уравнениями $w_i(\underline{z})=x_i$ и $w_{i+1}(\underline{z})=x_{i+1}$ системы $S$ мы вставим новое уравнение вида $y_j=x_j$, где %$y_j$ --- новая переменная. Не более чем через $n-1$ таких шагов от исходной системы $S$ мы придем к системе уравнений, удовлетворяющей нужному нам условию и имеющей одинаковое с %$S$ число решений в $G$ (в $H$).
Для кортежа $\underline{x}$ определим кортежи $\underline{k},\underline{\hat{k}}\in\mathbb{N}^{n+1}$, так, как это было сделано в лемме~\ref{AddIIILem2} и, воспользовавшись определениями~(\ref{AddIII_NTWE_Q2}) и~(\ref{AddIIILemNWEQ3}), положим:
\begin{align}
W_i(\underline{z}) &\overset{\textrm{def}}{=} J_{\underline{k}}\big(z_0,\dots,z_{n-1},z_i\big)\ \ \text{при}\ \ i=0,\ \dots,\ n-1;\label{Lem4EQ1}\\
\hat{W}_i(\underline{z}) &\overset{\textrm{def}}{=} J_{\underline{\hat{k}}}\big(z_0,\dots,z_{n-1},z_i\big)\ \ \text{при}\ \ i=0,\ \dots,\ n-1;\label{Lem4EQ2}\\
E(\underline{z})&\overset{\textrm{def}}{=} E_{\underline{k}[n,0]}\big(z_0,\dots,z_{n-1}\big);\label{Lem4EQ3}\\
\hat{E}(\underline{z})&\overset{\textrm{def}}{=} E_{\underline{\hat{k}}[n,0]}\big(z_0,\dots,z_{n-1}\big);\label{Lem4EQ4}\\
X_i &=W_i(\underline{x})= J_{\underline{k}}\big(x_0,\dots,x_{n-1},x_{i}\big)\ \ \text{при}\ \ i=0,\ \dots,\ n-1;\nonumber\\
\hat{X}_i &=\hat{W}_i(\underline{x})= J_{\underline{\hat{k}}}\big(x_0,\dots,x_{n-1},x_{i}\big)\ \ \text{при}\ \ i=0,\ \dots,\ n-1;\nonumber\\
U&=E(\underline{x})= E_{\underline{k}[n,0]}\big(x_0,\dots,x_{n-1}\big);\nonumber\\
V&=\hat{E}(\underline{x})= E_{\underline{\hat{k}}[n,0]}\big(x_0,\dots,x_{n-1}\big).\nonumber
\end{align}
Подчеркнем, что слова $W_i(\underline{z})$, $\hat{W}_i(\underline{z})$, $E(\underline{z})$ и $\hat{E}(\underline{z})$ зависят от кортежа~$\underline{x}$ (т.к. выбор кортежей~$\underline{k}$ и~$\underline{\hat{k}}$ зависит от кортежа~$\underline{x}$).

Покажем, что из леммы~\ref{AddIIILem2} следует, что выполнены следующие условия:
 \begin{enumerate}
    \item[\sffamily $\mathrm{F1)}$] $U,\ V,\ X_i,\ \hat{X}_i$ при $i=0,\ \dots,\ n-1$ суть гиперболические элементы группы $H$;
    \item[\sffamily $\mathrm{F2)}$] $[X_i,\hat{X}_i]\ne 1$ при $i=0,\ \dots,\ n-1$;
    \item[\sffamily $\mathrm{F3)}$] $[\hat{X}_{i},X_{i+1}]\ne 1$ при $i=0,\ \dots,\ n-2$;
    \item[\sffamily $\mathrm{F4)}$] $[\hat{X}_{n-1},U]\ne1$;
    \item[\sffamily $\mathrm{F5)}$] $[U,V]\ne1$.
 \end{enumerate}
Условие $\mathrm{F1)}$ следует из леммы~\ref{AddIIILem2} в силу того, что радикальная длина определена лишь для гиперболических элементов группы $H$.
%Т.к. из пунктов $a)$ и $b)$ леммы \ref{AddIIILem2} следует, что среди элементов $U,\ V,\ X_i,\ \hat{X}_i$, $i=0,\ \dots,\ n-1$ наименьшую радикальную длину имеет элемент $U$, причем $|U|_r\ge 2$, то $F1)$ установлено.
Условия $\mathrm{F2)}$ и $\mathrm{F3)}$ прямо следуют из пункта~$b)$ леммы~\ref{AddIIILem2}, т.к. коммутирующие гиперболические элементы имеют равную радикальную длину. Аналогичные рассуждения показывают, что условия $\mathrm{F4)}$ и $\mathrm{F5)}$ так же выполнены.

Следовательно, кортеж (длины  $2n+2$):
\begin{equation}\label{DefTildex-Tuple}
\underline{\tilde{x}}=\bigl(X_0,\hat{X}_0,X_1,\hat{X}_1,\dots,X_{n-1},\hat{X}_{n-1},U,V\bigr)
\end{equation}
состоит из гиперболических элементов, причем соседние его компоненты не коммутируют.
%Значит для кортежа $\underline{\dot{x}}$, в силу леммы \ref{LemMLPartII}, существует слово $T_{\underline{X}}(z_0,\dots,z_{2n+1})$.

$a)$ Определим слово $P_{\underline{x}}(z_0,\dots,z_{n-1})\in F_n(z_0,\dots,z_{n-1})'$ следующим образом:
\begin{equation}\label{DefP-words}
P_{\underline{x}}(z_0,\dots,z_{n-1})\overset{\textrm{def}}{=}T_{\underline{\tilde{x}}}\Big(W_0(\underline{z}),\hat{W}_0(\underline{z}),\dots,W_{n-1}(\underline{z}),\hat{W}_{n-1}(\underline{z}),E(\underline{z}),\hat{E}(\underline{z})\Big),
\end{equation}
где слова $W_i(\underline{z})$, $\hat{W}_i(\underline{z})$, $E(\underline{z})$ и $\hat{E}(\underline{z})$ определены как~(\ref{Lem4EQ1}), (\ref{Lem4EQ2}), (\ref{Lem4EQ3}) и~(\ref{Lem4EQ4}) соответственно, а слово $T_{\underline{\tilde{x}}}(z_0,\dots,z_{2n+1})$ выбрано для кортежа~(\ref{DefTildex-Tuple}) как в пункте~$a)$ леммы~\ref{LemMLPartII}. Заметим, что данное определение корректно, т.к. кортеж~(\ref{DefTildex-Tuple}) удовлетворяет всем условиям леммы~\ref{LemMLPartII} (т.е. все его компоненты суть гиперболические элементы и его соседние компоненты не коммутируют).

Покажем, что если выполнено равенство~(\ref{LemMLPartIIIMEQ}), т.е. если:
\begin{align}\label{Lem3MainEQtr}
T_{\underline{\tilde{x}}}\Big(X_0,\hat{X}_0,\dots,X_{n-1},\hat{X}_{n-1},U,V\Big)=T_{\underline{\tilde{x}}}\Big(W_0(\underline{y}),\hat{W}_0(\underline{y}),\dots,W_{n-1}(\underline{y}),\hat{W}_{n-1}(\underline{y}),E(\underline{y}),\hat{E}(\underline{y})\Big),
\end{align}
то найдется такой элемент $v\in \left \langle P_{\underline{x}}(x_0,\dots,x_{n-1}) \right \rangle_\infty$, что $x_i = y_i^v$ при $i=0,\ \dots,\ n-1.$

%\begin{align}\label{Lem3MainEQtr}
%T_{\underline{\tilde{x}}}\Big(W_0(\underline{x}),\hat{W}_0(\underline{x}),\dots,W_{n-1}&(\underline{x}),\hat{W}_{n-1}(\underline{x}),E(\underline{x}),\hat{E}(\underline{x})\Big)=\\\nonumber
%&=T_{\underline{\tilde{x}}}\Big(W_0(\underline{y}),\hat{W}_0(\underline{y}),\dots,W_{n-1}(\underline{y}),\hat{W}_{n-1}(\underline{y}),E(\underline{y}),\hat{E}(\underline{y})\Big),
%\end{align}

%\begin{equation}\label{Lem3MainEQ}
%T_{\underline{X}}\big(W_0(\underline{z}),\hat{W}_0(\underline{z}),\dots,W_{n-1}(\underline{z}),\hat{W}_{n-1}(\underline{z}),E(\underline{z}),\hat{E}(\underline{z})\big)=T_{\underline{X}}\big(X_0,\hat{X}_0,\dots,X_{n-1},\hat{X}_{n-1},U,V\big).
%\end{equation}
%Данное уравнение имеет решение в $G$ по построению (действительно, в качестве $\underline{z}$ достаточно взять решение системы $S$), следовательно, в силу вербальной замкнутости %подгруппы $H$ в $G$, это уравнение имеет решение, скажем $\underline{z}=\underline{h}$, и в $H$. Т.е. в $H$ выполнено равенство:
%\begin{equation}\label{Lem3MainEQtr}
%T_{\underline{X}}\big(W_0(\underline{h}),\hat{W}_0(\underline{h}),\dots,W_{n-1}(\underline{h}),\hat{W}_{n-1}(\underline{h}),E(\underline{h}),\hat{E}(\underline{h})\big)=T_{\underline{X}}\big(X_0,\hat{X}_0,\dots,X_{n-1},\hat{X}_{n-1},U,V\big).
%\end{equation}

Чтобы применить пункт~$a)$ леммы~\ref{LemMLPartII} к равенству~(\ref{Lem3MainEQtr}) достаточно показать, что все компоненты кортежа:
$$\underline{W}=\big(W_0(\underline{y}),\hat{W}_0(\underline{y}),\dots,W_{n-1}(\underline{y}),\hat{W}_{n-1}(\underline{y}),E(\underline{y}),\hat{E}(\underline{y})\big)$$
суть гиперболические элементы группы $H$. Покажем, что если хотя бы одна из компонент кортежа $\underline{W}$ не является гиперболическим элементом, то равенство равенство~(\ref{Lem3MainEQtr}) не может быть выполнено.

Для этого прежде всего заметим, что ни одна из компонент кортежа $\underline{W}$ не равна единице, т.к. в противном случае, в силу пункта~$E2)$ леммы~$\ref{AddIIINewLemIGProp}$, правая часть равенства~(\ref{Lem3MainEQtr}) равна единице, а это противоречит пункту~$E1)$ леммы~$\ref{AddIIINewLemIGProp}$ (точнее, здесь мы пользуемся пунктами~$E1)$ и~$E2)$ леммы~$\ref{AddIIINewLemIGProp}$, примененной к кортежу $\underline{\tilde{x}}$, для пары индексов $(j,i)=(2n+2,0)$).

Пусть $E(\underline{y})=E_{\underline{k}[n,0]}\big(y_0,\dots,y_{n-1}\big)\in H_\mathfrak{a}^g\!\setminus\!\{1\}$ для некоторого $\mathfrak{a}\in \mathfrak{A}$ и $g\in H$. Тогда, в силу леммы~\ref{AddIIILem1}, мы имеем $y_0,\dots,y_{n-1}\in H_\mathfrak{a}^g\!\setminus\!\{1\}.$ Значит, в силу определений~(\ref{Lem4EQ1}), (\ref{Lem4EQ2}), (\ref{Lem4EQ3}) и~(\ref{Lem4EQ4}), правая часть равенства~(\ref{Lem3MainEQtr}) есть элемент из $H_\mathfrak{a}^g$, но это противоречит тому, что (в силу пункта~$E1)$ леммы~\ref{AddIIINewLemIGProp}) левая часть равенства~(\ref{Lem3MainEQtr}) есть гиперболический элемент.

Случай $\hat{E}(\underline{y})\in H_\mathfrak{a}^g\!\setminus\!\{1\}$ полностью аналогичен разобранному.

Пусть для некоторого $i=0,\ \dots,\ n-1$ мы имеем $W_i(\underline{y})\in H_\mathfrak{a}^g\!\setminus\!\{1\}$, т.е.:
\begin{align}\label{Lem3MainEQ2Jk}
W_i(\underline{y})&=J_{\underline{k}}\bigl(y_0,\dots,y_{n-1},y_i\bigr);\nonumber\\
&=\bigl[E_{\underline{k}[n,0]}\big(y_0,\dots,y_{n-1}\big)^{2m_2},y_i^2\bigr];\nonumber\\
&=\bigl[E(\underline{y})^{2m_2},y_i^2\bigr]\in H_\mathfrak{a}^g\!\setminus\!\{1\},
\end{align}
где $m_2$ есть $(n+1)$-я компонента кортежа $\underline{k}$, а при первом, втором и третьем переходе мы воспользовались определениями~(\ref{Lem4EQ1}), (\ref{AddIIILemNWEQ3}) и~(\ref{Lem4EQ3}) соответственно. Из~(\ref{Lem3MainEQ2Jk}), в силу пункта~$C1)$ леммы~\ref{MCL}, легко следует, что $E(\underline{y})\in H_\mathfrak{a}^g\!\setminus\!\{1\}$ и мы приходим к уже разобранному случаю.

Случай $\hat{W}_i(\underline{y})\in H_\mathfrak{a}^g\!\setminus\!\{1\}$ полностью аналогичен случаю $W_i(\underline{y})\in H_\mathfrak{a}^g\!\setminus\!\{1\}$.

Итак, мы показали, что все компоненты кортежа $\underline{W}$ суть гиперболические элементы. Значит к равенству~(\ref{Lem3MainEQtr}) применим пункт~$a)$ леммы~\ref{LemMLPartII}, т.е. найдется такой элемент $v\in H$, что:
\begin{align}
\big[U^{2m_2},x_i^2\big]=X_i&=W_i(\underline{y})^v=\big[E(\underline{y})^{2m_2},y_i^2\big]^v\ \ \text{при}\ \ i=0,\ \dots,\ n-1;\label{Lem3MainEQAlthe1}\\
\big[V^{2\hat{m}_2},x_i^2\big]=\hat{X}_i&=\hat{W}_i(\underline{y})^v=\big[\hat{E}(\underline{y})^{2\hat{m}_2},y_i^2\big]^v\ \ \text{при}\ \ i=0,\ \dots,\ n-1;\label{Lem3MainEQAlthe2}\\
U&=E(\underline{y})^v;\label{Lem3MainEQAlthe3}\\
V&=\hat{E}(\underline{y})^v\label{Lem3MainEQAlthe4}
\end{align}
%\begin{align}
%\big[E(\underline{y})^{2m_2},y_i^2\big]=W_i(\underline{y})&=X_i^v=\big[U^{2m_2},x_i^2\big]^v\ \ \text{при}\ \ i=0,\ \dots,\ n-1;\label{Lem3MainEQAlthe1}\\
%\big[\hat{E}(\underline{y})^{2\hat{m}_2},y_i^2\big]=\hat{W}_i(\underline{y})&=\hat{X}_i^v=\big[V^{2\hat{m}_2},x_i^2\big]^v\ \ \text{при}\ \ i=0,\ \dots,\ %n-1;\label{Lem3MainEQAlthe2}\\
%E(\underline{y})&=U^v;\label{Lem3MainEQAlthe3}\\
%\hat{E}(\underline{y})&=V^v\label{Lem3MainEQAlthe4}
%\end{align}
(чтобы легче понять левые и правые равенства в~(\ref{Lem3MainEQAlthe1}) и~(\ref{Lem3MainEQAlthe2}) см. цепочку преобразований~(\ref{Lem3MainEQ2Jk})).

Подставляя~(\ref{Lem3MainEQAlthe3}) в~(\ref{Lem3MainEQAlthe1}) и~(\ref{Lem3MainEQAlthe4}) в~(\ref{Lem3MainEQAlthe2}) мы получаем:
\begin{align}
\big[U^{2m_2},x_i^2\big]&=\big[U^{2m_2},y_i^{2v}\big]\ \ \text{при}\ \ i=0,\ \dots,\ n-1;\label{Lem3MainEQAAlthe1}\\
\big[V^{2\hat{m}_2},x_i^2\big]&=\big[V^{2\hat{m}_2},y_i^{2v}\big]\ \ \text{при}\ \ i=0,\ \dots,\ n-1.\label{Lem3MainEQAAlthe2}
\end{align}
Следовательно:
\begin{align*}
x_i^{2}y_i^{-2v}\in C_H(U^{2m_2})\ \ \text{при}\ \ i=0,\ \dots,\ n-1;\\
x_i^{2}y_i^{-2v}\in C_H(V^{2\hat{m}_2})\ \ \text{при}\ \ i=0,\ \dots,\ n-1.
\end{align*}
Т.к. $U$ и $V$ суть некоммутирующие гиперболические элементы (см. $\mathrm{F1)}$ и $\mathrm{F5)}$), то из теоремы Куроша легко вывести, что $C_H(U^{2m_2})\cap C_H(V^{2\hat{m}_2})=\{1\}$ для любых $m_2,\hat{m}_2\in\mathbb{N}$. Значит мы имеем:
\begin{equation}
x_i^{2}=y_i^{2v}\ \ \text{при}\ \ i=0,\ \dots,\ n-1,
\end{equation}
откуда, учитывая однозначность извлечения в $H$ корня из гиперболического элемента, мы получаем нужные нам равенства.

Осталось заметить, что, в силу пункта~$a)$ леммы~\ref{LemMLPartII}, мы имеем $v\in \big\langle T_{\underline{\tilde{x}}}(\underline{\tilde{x}}) \big\rangle_\infty$ и что $T_{\underline{\tilde{x}}}(\underline{\tilde{x}})=P_{\underline{x}}(\underline{x})$ в силу определения~(\ref{DefP-words}).

$b)$ Определим слова $P_{\underline{x}}'(z_0,\dots,z_{n-1}),P_{\underline{x}}''(z_0,\dots,z_{n-1})\in F_n(z_0,\dots,z_{n-1})'$ следующим образом:
\begin{align*}\label{DefP-wordspartb}
P_{\underline{x}}'(z_0,\dots,z_{n-1})&\overset{\textrm{def}}{=}T_{\underline{\tilde{x}}}'\Big(W_0(\underline{z}),\hat{W}_0(\underline{z}),\dots,W_{n-1}(\underline{z}),\hat{W}_{n-1}(\underline{z}),E(\underline{z}),\hat{E}(\underline{z})\Big);\\
P_{\underline{x}}''(z_0,\dots,z_{n-1})&\overset{\textrm{def}}{=}T_{\underline{\tilde{x}}}''\Big(W_0(\underline{z}),\hat{W}_0(\underline{z}),\dots,W_{n-1}(\underline{z}),\hat{W}_{n-1}(\underline{z}),E(\underline{z}),\hat{E}(\underline{z})\Big),
\end{align*}
где слова $W_i(\underline{z})$, $\hat{W}_i(\underline{z})$, $E(\underline{z})$ и $\hat{E}(\underline{z})$ определены как~(\ref{Lem4EQ1}), (\ref{Lem4EQ2}), (\ref{Lem4EQ3}) и~(\ref{Lem4EQ4}) соответственно, а слова $T_{\underline{\tilde{x}}}'(z_0,\dots,z_{2n+1})$ и $T_{\underline{\tilde{x}}}''(z_0,\dots,z_{2n+1})$ выбраны для кортежа~(\ref{DefTildex-Tuple}) как в пункте~$b)$ леммы~\ref{LemMLPartII}. Заметим, что данное определение корректно, т.к. кортеж~(\ref{DefTildex-Tuple}) удовлетворяет всем условиям леммы~\ref{LemMLPartII} (т.е. все его компоненты суть гиперболические элементы и его соседние компоненты не коммутируют).

 Мы опускаем дальнейшие рассуждения, т.к. они полностью аналогичны рассуждениям, приведенным при доказательстве пункта~$a).$ Отметим лишь, что когда мы применим пункт~$b)$ леммы~\ref{LemMLPartII} к равенствам~(\ref{LemMLPartIIIMEQpartb}), то мы получим равенства вида~(\ref{Lem3MainEQAlthe1}), (\ref{Lem3MainEQAlthe2}), (\ref{Lem3MainEQAlthe3}) и~(\ref{Lem3MainEQAlthe4}), но с дополнительным условием $v=1$.
\end{ProofMLIII}

\section{Завершение доказательства Леммы о словах, доказательства Основной теоремы и Утверждения о гомоморфизмах}

%\begin{Lemma}\label{TFML} Пусть группа $H=\zvezda\limits_{\mathfrak{a}\in \mathfrak{A}}\!H_\mathfrak{a}$ есть свободное произведение нетривиальных групп $H_\mathfrak{a}$ и пусть %$\underline{x}=(x_0,\dots,x_{n-1})$ есть такой кортеж элементов из $H$, что подгруппа $\left \langle x_0,\dots,x_{n-1} \right \rangle$ содержит хотя бы два некоммутирующих %гиперболических элемента группы $H$. Тогда существует слово $M_{\underline{x}}(z_0,\dots,z_{n-1})\in F_n(z_0,\dots,z_{n-1})'$, вообще говоря, зависящее от кортежа $\underline{x}$, %такое, что если для некоторого кортежа $(y_0,\dots,y_{n-1})$ элементов $y_i\in H$ выполнено равенство:
%\begin{equation}\label{AddPartIVEQ1}
%M_{\underline{x}}(x_0,\dots,x_{n-1})=M_{\underline{x}}(y_0,\dots,y_{n-1}),
%\end{equation}
%то найдется такой элемент $v\in \left \langle M_{\underline{x}}(x_0,\dots,x_{n-1}) \right \rangle_\infty$, что $x_i = v^{-1}y_iv$ при $i=0,\ \dots,\ n-1.$
%\end{Lemma}
Перейдем непосредственно к доказательству Леммы о словах.

\begin{ProofWL} Пусть $M_H$ есть множество всех гиперболических элементов группы $H$. Прежде всего заметим, что если подгруппа $S$ группы $H$ содержит хотя бы два некоммутирующих гиперболических элемента группы $H$, то для любого конечного подмножества $X\subset S\!\setminus\!\{1\}$ найдется элемент из множества $S\cap M_H$, не коммутирующий ни с одним из элементов множества $X$. Действительно, разложение Куроша такой подгруппы $S$ должно содержать хотя бы два свободных сомножителя, причем если свободных сомножителей ровно два, то хотя бы один из них должен содержать не меньше трех элементов. Из этого легко следует, что в $S\cap M_H$ есть слова сколь угодно большой радикальной длины. В качестве искомого элемента можно взять любой элемент из $S\cap M_H$, радикальная длина которого больше радикальной длины любого гиперболического элемента из $X$. Построение кортежа $\underline{\hat{x}}$, приведенное ниже, опирается на существование таких элементов.

Положим $S=\left \langle x_0,\dots,x_{n-1} \right \rangle$ и преобразуем кортеж $\underline{x}=(x_0,\dots,x_{n-1})$ следующим образом.

1) Выберем такой элемент $w\in S\!\setminus\!\{1\}$, что $x_iw\ne 1$ при $i=0,\ \dots,\ n-1$ и составим кортеж $(x_0w,\dots,x_{n-1}w,w).$ Далее, выберем элементы $u_1,u_2\in S\cap M_H$ так, чтобы $[u_1,u_2]\ne 1$, $[w,u_1]\ne1$, $[w,u_2]\ne1$ и $[x_iw,u_1]\ne1$, $[x_iw,u_2]\ne1$ при $i=0,\ \dots,\ n-1$. Положим $\dot{x}_{i,1}=[x_iw,u_1]$ и $\dot{x}_{i,2}=[x_iw,u_2]$ при $i=0,\ \dots,\ n-1$ и составим кортеж (длины $2n+4$):
$$\underline{\dot{x}}=(\dot{x}_{0,1},\dot{x}_{0,2},\dots,\dot{x}_{n-1,1},\dot{x}_{n-1,2},[w,u_1],[w,u_2],u_1,u_2).$$
Заметим, что все компоненты кортежа $\underline{\dot{x}}$ суть элементы множества $S\cap M_H$ (действительно, $u_1$ и $u_2$ суть гиперболические элементы группы $H$ по построению, а остальные компоненты являются гиперболическими элементами группы $H$ в силу пункта~$C1)$ леммы~\ref{MCL}).

2) Выберем элемент $s\in S\cap M_H$ так, чтобы $s$ не коммутировал ни с одной из компонент кортежа $\underline{\dot{x}}$ и составим кортеж (длины $4n+6$):
\begin{equation}\label{finTupleX}
\underline{\hat{x}}=(\dot{x}_{0,1},s,\dot{x}_{0,2},s,\dots,\dot{x}_{n-1,1},s,\dot{x}_{n-1,2},s,[w,u_1],s,[w,u_2],s,u_1,u_2),
\end{equation}
полученный из кортежа $\underline{\dot{x}}$ посредством <<вставки>> элемента $s$ между каждыми соседними его компонентами, кроме двух последних.
Заметим, что все компоненты кортежа $\underline{\hat{x}}$ суть элементы множества $S\cap M_H$, причем его соседние компоненты не коммутируют.

Положим $\underline{z}= (z_0,\dots,z_{n-1})$. Пусть $f_w(\underline{z}),\ f_{u_1}(\underline{z}),\ f_{u_2}(\underline{z})$ и $f_{s}(\underline{z})$ суть такие слова свободной группы $F_n(\underline{z})$, что в $S$ мы имеем:
\begin{align}\label{AddPartIVEQ3}
f_w(x_0,\dots,x_{n-1})=w,\ \ f_{u_1}(x_0,\dots,x_{n-1})=u_1,\ \  f_{u_2}(x_0,\dots,x_{n-1})=u_2,\ \  f_{s}(x_0,\dots,x_{n-1})=s.
\end{align}
Докажем пункт~$a)$. Определим слово $M_{\underline{x}}(z_0,\dots,z_{n-1})\in F_n(z_0,\dots,z_{n-1})'$ следующим образом:
\begin{equation}\label{you-know-what}
M_{\underline{x}}(\underline{z})\overset{\textrm{def}}{=}P_{\underline{\hat{x}}}\Big(f_{\dot{x}_{0,1}}(\underline{z}),f_s(\underline{z}),\dots,f_{\dot{x}_{n-1,2}}(\underline{z}),f_s(\underline{z}),f_{[w,u_1]}(\underline{z}),f_s(\underline{z}),f_{[w,u_2]}(\underline{z}),f_s(\underline{z}),f_{u_1}(\underline{z}),f_{u_2}(\underline{z})\Big),
\end{equation}
где $f_{\dot{x}_{i,j}}(\underline{z})\overset{\textrm{def}}{=}\big[z_if_w(\underline{z}),f_{u_j}(\underline{z})\big]$ и $f_{[w,u_j]}(\underline{z})\overset{\textrm{def}}{=}\big[f_w(\underline{z}),f_{u_j}(\underline{z})\big]$ при $i=0,\ \dots,\ n-1$ и $j=1,2,$ а слово $P_{\underline{\hat{x}}}(z_0,\dots,z_{4n+5})$ выбрано для кортежа~(\ref{finTupleX}) как в пункте~$a)$ леммы~\ref{LemMLPartIII}. Заметим, что данное определение корректно, т.к. кортеж $\underline{\hat{x}}$ удовлетворяет всем условиям леммы~\ref{LemMLPartIII} (т.е. все его компоненты суть гиперболические элементы и его соседние компоненты не коммутируют).
%Заметим, что данное определение корректно, т.к. слово $P_{\underline{\hat{x}}}(z_0,\dots,z_{4n+5})$ может быть построено для кортежа $\underline{\hat{x}}$ в силу леммы %\ref{LemMLPartIII}.

Покажем, что слово $M_{\underline{x}}(\underline{z})$ действительно обладает описанным в пункте~$a)$ Леммы о словах свойством. Пусть для некоторого кортежа $\underline{y}= (y_0,\dots,y_{n-1})$ элементов $y_i\in H$ выполнено равенство~(\ref{LemMLPartIIIMEQff}). Тогда, в силу определения~(\ref{you-know-what}) и леммы~\ref{LemMLPartIII}, найдется такой элемент $v\in \left \langle M_{\underline{x}}(x_0,\dots,x_{n-1}) \right \rangle_\infty=\left \langle P_{\underline{\hat{x}}}(\underline{\hat{x}}) \right \rangle_{\infty}$, что одновременно выполнены равенства:
\begin{align*}
\big[x_if_w(\underline{x}),f_{u_1}(\underline{x})\big]&=\big[y_if_w(\underline{y}),f_{u_1}(\underline{y})\big]^v\ \ \text{при}\ \ i=0,\ \dots,\ n-1;\\
\big[x_if_w(\underline{x}),f_{u_2}(\underline{x})\big]&=\big[y_if_w(\underline{y}),f_{u_2}(\underline{y})\big]^v\ \ \text{при}\ \ i=0,\ \dots,\ n-1;\\
f_{u_j}(\underline{x})&=f_{u_j}(\underline{y})^v\ \ \text{при}\ \ j=1,2;\\
\big[f_w(\underline{x}),f_{u_j}(\underline{x})\big]&=\big[f_w(\underline{y}),f_{u_j}(\underline{y})\big]^v\ \ \text{при}\ \ j=1,2;\\
f_{s}(\underline{x})&=f_{s}(\underline{y})^v.
\end{align*}
Учитывая~(\ref{AddPartIVEQ3}), из первых четырех наборов равенств мы получаем\footnote[$\dagger$]{Отметим, что равенство $f_{s}(\underline{x})=f_{s}(\underline{y})^v$ мы никак не используем, т.к. оно является <<побочным продуктом>> введения в кортеж $\underline{\hat{x}}$ компонент $s$, необходимых для того, чтобы гарантировать то, что стоящие рядом компоненты кортежа $\underline{\hat{x}}$ не коммутируют.}:
\begin{align*}
\bigl( y_if_w(\underline{y})\bigr)^v\big(x_iw\big)^{-1}&\in C_H(u_1)\ \ \text{при}\ \ i=0,\ \dots,\ n-1;\\
\big(y_if_w(\underline{y})\big)^v\big(x_iw\big)^{-1}&\in C_H(u_2)\ \ \text{при}\ \ i=0,\ \dots,\ n-1;\\
f_w(\underline{y})^vw^{-1}&\in C_H(u_1);\\
f_w(\underline{y})^vw^{-1}&\in C_H(u_2).
\end{align*}
Т.к. гиперболические элементы $u_1$ и $u_2$ группы $H$ не коммутируют (по построению), то из теоремы Куроша легко вывести, что  $C_H(u_1)\cap  C_H(u_2)=\{1\}.$ Следовательно, мы имеем:
\begin{align*}
x_iw&=y_i^vf_w(\underline{y})^v\ \ \text{при}\ \ i=0,\ \dots,\ n-1;\\
w&=f_w(\underline{y})^v,
\end{align*}
откуда немедленно следуют необходимые нам равенства.

Докажем пункт~$b)$. Определим слова $M_{\underline{x}}'(z_0,\dots,z_{n-1}),M_{\underline{x}}''(z_0,\dots,z_{n-1})\in F_n(z_0,\dots,z_{n-1})'$ следующим образом:
\begin{align*}
M_{\underline{x}}'(\underline{z})&\overset{\textrm{def}}{=}P_{\underline{\hat{x}}}'\Big(f_{\dot{x}_{0,1}}(\underline{z}),f_s(\underline{z}),\dots,f_{\dot{x}_{n-1,2}}(\underline{z}),f_s(\underline{z}),f_{[w,u_1]}(\underline{z}),f_s(\underline{z}),f_{[w,u_2]}(\underline{z}),f_s(\underline{z}),f_{u_1}(\underline{z}),f_{u_2}(\underline{z})\Big);\\
M_{\underline{x}}''(\underline{z})&\overset{\textrm{def}}{=}P_{\underline{\hat{x}}}''\Big(f_{\dot{x}_{0,1}}(\underline{z}),f_s(\underline{z}),\dots,f_{\dot{x}_{n-1,2}}(\underline{z}),f_s(\underline{z}),f_{[w,u_1]}(\underline{z}),f_s(\underline{z}),f_{[w,u_2]}(\underline{z}),f_s(\underline{z}),f_{u_1}(\underline{z}),f_{u_2}(\underline{z})\Big),
\end{align*}
где $f_{\dot{x}_{i,j}}(\underline{z})\overset{\textrm{def}}{=}\big[z_if_w(\underline{z}),f_{u_j}(\underline{z})\big]$ и $f_{[w,u_j]}(\underline{z})\overset{\textrm{def}}{=}\big[f_w(\underline{z}),f_{u_j}(\underline{z})\big]$ при $i=0,\ \dots,\ n-1$ и $j=1,2,$ а слова $P_{\underline{\hat{x}}}'(z_0,\dots,z_{4n+5})$ и $P_{\underline{\hat{x}}}''(z_0,\dots,z_{4n+5})$ выбраны для кортежа~(\ref{finTupleX}) как в пункте~$b)$ леммы~\ref{LemMLPartIII}. Заметим, что данное определение корректно, т.к. кортеж $\underline{\hat{x}}$ удовлетворяет всем условиям леммы~\ref{LemMLPartIII} (т.е. все его компоненты суть гиперболические элементы и его соседние компоненты не коммутируют).

Мы опускаем дальнейшие рассуждения, т.к. они полностью аналогичны рассуждениям, приведенным при доказательстве пункта~$a)$ (единственное отличие --- это дополнительное условие $v=1$).
\end{ProofWL}

Прежде чем перейти к доказательству Основной теоремы, напомним следующую хорошо известную лемму.

\begin{Lemma}{\rm\cite{[KM17]}}\label{ProofMT1TL} Если подгруппа $H$ группы $G$ такова, что любая система уравнений вида:
\begin{equation}\label{ProofMTEq1}
\left \{ w_i(\underline{z})=h_i\ |\ w_i(\underline{z})\in F(\underline{z}),\ h_i\in H,\ i=0,\ \dots,\ n-1 \right \},
\end{equation}
имеющая решение в $G$, имеет решение в $H$, то $H$ алгебраически замкнута в $G$.
\end{Lemma}

\begin{ProofMT} Пусть группа $H=\zvezda\limits_{\mathfrak{a}\in \mathfrak{A}}\!H_\mathfrak{a}$ есть нетривиальное свободное произведение групп $H_\mathfrak{a}$, причем $H\not\simeq \left \langle a \right \rangle_2\ast\left \langle b \right \rangle_2$ (как уже отмечалось во введении, сильная вербальная замкнутость группы $D_\infty\simeq \left \langle a \right \rangle_2\ast\left \langle b \right \rangle_2$ была установлена в работе~\cite{[KlMaMi]}). В таком случае ясно, что в $H$ найдутся два некоммутирующих гиперболических элемента. В силу леммы~\ref{ProofMT1TL} достаточно доказать, что любая система уравнений вида~(\ref{ProofMTEq1}), имеющая решение в $G$, имеет решение в $H$.

Пусть система~(\ref{ProofMTEq1}) имеет решение в $G$ и пусть $u_1$ и $u_2$ суть некоторые некоммутирующие гиперболические элементы группы $H$. Тогда кортеж:
$$\underline{x}=(h_0,\dots,h_{n-1},u_1,u_2)$$
удовлетворяет условию Леммы о словах. Рассмотрим уравнение:
\begin{equation}\label{PartIVMEQ}
M_{\underline{x}}\bigl(w_0(\underline{z}),\dots,w_{n-1}(\underline{z}),\dot{z}_1,\ddot{z}_{2}\bigr)=M_{\underline{x}}\bigl(h_0,\dots,h_{n-1},u_1,u_{2}\bigr),
\end{equation}
где слово $M_{\underline{x}}(z_0,\dots,z_{n+1})$ выбрано для кортежа $\underline{x}$ как в пункте~$a)$ Леммы о словах, а $\dot{z}_1$ и $\ddot{z}_{2}$ суть новые (т.е. не присутствующие в кортеже $\underline{z}$) независимые переменные.

Уравнение~(\ref{PartIVMEQ}) имеет решение в $G$ по построению (действительно, в качестве $\underline{z}$ достаточно взять решение системы~(\ref{ProofMTEq1}) и положить $\dot{z}_1=u_1$ и $\ddot{z}_{2}=u_2$). Следовательно, в силу вербальной замкнутости подгруппы $H$ в группе $G$, это уравнение имеет некоторое решение, скажем $\underline{z}=\underline{y}$, $\dot{z}_1=\dot{y}_1$ и $\ddot{z}_2=\ddot{y}_2,$ в $H$. Значит, в силу Леммы о словах, найдется такой элемент $v\in H$, что $w_i(\underline{y})^v=h_i$ при $i=0,\ \dots,\ n-1,$ т.е. $\underline{z}=\underline{y}^{v}$ есть искомое решение системы~(\ref{ProofMTEq1}) в $H$.
\end{ProofMT}

\begin{ProofHomStatement} Пусть $S=\{s_0,\dots,s_{n-1}\}$ есть порождающее множество группы $G$. По условию подгруппа $\big\langle\psi(s_0),\dots,\psi(s_{n-1})\big\rangle=\textrm{Im}\bigl(\psi(G)\bigr)$ группы $H$ содержит хотя бы два некоммутирующих гиперболических элемента группы $H$, т.е. к кортежу $\underline{\psi(S)}=\big(\psi(s_0),\dots,\psi(s_{n-1})\big)$ применим пункт~$b)$ Леммы о словах. Положим:
$$w_1\overset{\textrm{def}}{=}M_{\underline{\psi(S)}}'(s_0,\dots,s_{n-1})\ \ \text{и}\ \  w_2\overset{\textrm{def}}{=}M_{\underline{\psi(S)}}''(s_0,\dots,s_{n-1}),$$
где слова $M_{\underline{\psi(S)}}'(z_0,\dots,z_{n-1})$ и $M_{\underline{\psi(S)}}''(z_0,\dots,z_{n-1})$ выбраны для кортежа $\underline{\psi(S)}$ как в пункте~$b)$ Леммы о словах.

Если для некоторого гомоморфизма $\varphi \in \text{\textrm{Hom}}(G,H)$ выполнены равенства $\psi(w_1)=\varphi(w_1)$ и $\psi(w_2)=\varphi(w_2)$, т.е. если:
%то для кортежа $\underline{\varphi(S)}=\big(\varphi(s_0),\dots,\varphi(s_{n-1})\varphi)$ мы имеем равенства:
\begin{align*}
M_{\underline{\psi(S)}}'\big(\psi(s_0),\dots,\psi(s_{n-1})\big)&=M_{\underline{\psi(S)}}'\big(\varphi(s_0),\dots,\varphi(s_{n-1})\big);\\
M_{\underline{\psi(S)}}''\big(\psi(s_0),\dots,\psi(s_{n-1})\big)&=M_{\underline{\psi(S)}}''\big(\varphi(s_0),\dots,\varphi(s_{n-1})\big),
\end{align*}
то, в силу пункта~$b)$ Леммы о словах, мы имеем $\psi(s_i)=\varphi(s_i)$ при $i=0,\ \dots,\ n-1$. Т.к. $S$ есть порождающее множество группы $G$, то это означает, что $\psi=\varphi$.
\end{ProofHomStatement}

\section{Благодарности}

%Работа выполнена при поддержке Российского фонда фундаментальных
%исследований, грант 15-01-05823.

Автор выражает благодарность своему научному руководителю А.А.~Клячко за постоянное внимание к работе и ценные замечания.

\end{document}